\documentclass[10pt]{amsart}
\usepackage{amsmath,amsthm}
\usepackage{comment}
\usepackage{amssymb}
\usepackage{stmaryrd}
\usepackage{color}
\usepackage{enumitem}
\usepackage{pdfsync}
\usepackage{tikz-cd}
\usepackage{graphics}
\usepackage{tikz}
\usetikzlibrary{decorations.markings}
\usepackage{hyperref}
\usepackage[capitalise,noabbrev]{cleveref}

\newcommand{\defterm}[1]{\emph{#1}}

\numberwithin{equation}{subsection}

\theoremstyle{plain}
	\newtheorem{thm}[equation]{Theorem}
	\newtheorem*{thm*}{Theorem}	
	\newtheorem*{question*}{Question}

	\newtheorem{cor}[equation]{Corollary}
	\newtheorem*{cor*}{Corollary}
	\newtheorem{prop}[equation]{Proposition}
	\newtheorem*{prop*}{Proposition}
	\newtheorem{lem}[equation]{Lemma}
	\newtheorem*{lem*}{Lemma}
	
	\newtheorem*{ex*}{Exercise}
	
	\newtheorem*{claim*}{Claim}
	
	\newtheorem*{conj*}{Conjecture}
	
	\newtheorem*{fact*}{Fact}
	\newtheorem{notation}[equation]{Notation}
	\newtheorem*{notation*}{Notation}
\theoremstyle{definition}
	\newtheorem{Def}[equation]{Definition}
	\newtheorem*{Def*}{Definition}
	
	\newtheorem{soln*}{Solution}
	
	\newtheorem*{note*}{Note}
	\newtheorem{eg}[equation]{Example}
	\newtheorem*{eg*}{Example}
	\newtheorem{rmk}[equation]{Remark}
	\newtheorem*{rmk*}{Remark}
	
	\newtheorem*{obs*}{Observation}
	\newtheorem{warning}[equation]{Warning}
	\newtheorem*{warning*}{Warning}
	
\newcommand{\op}{\mathrm{op}}

\newcommand{\id}{\mathrm{id}}

\newcommand{\Mor}{\operatorname{Mor}}

\newcommand{\Hom}{\mathrm{Hom}}

\newcommand{\Aut}{\mathrm{Aut}}
\newcommand{\inv}{^{-1}}

\newcommand{\dual}{\vee}

\newcommand{\Set}{\mathsf{Set}}

\newcommand{\Cat}{\mathsf{Cat}}
\newcommand{\Top}{\mathsf{Top}}
\newcommand{\Fun}{\mathrm{Fun}}

\DeclareMathOperator{\Ind}{Ind}

\newcommand{\Sm}{\mathsf{Sm}}
\newcommand{\SMC}{\mathsf{SMC}}
\newcommand{\SMD}{\mathsf{SMD}}
\newcommand{\SMP}{\mathsf{SMP}}
\newcommand{\FinBij}{\mathsf{Fin}^{\mathrm{iso}}}
\newcommand{\Bord}{\mathsf{Bord}}
\newcommand{\Spt}{\mathsf{Spt}}
\newcommand{\Span}{\mathsf{Span}}
\newcommand{\Fin}{\mathsf{Fin}}
\newcommand{\Vect}{\mathsf{Vect}}

\newcommand{\ints}{\mathbb{Z}}

\newcommand{\nats}{\mathbb{N}}
\newcommand{\pjv}{\mathbb{P}}
\newcommand{\aff}{\mathbb{A}}
\newcommand{\GL}{\operatorname{GL}}
\newcommand{\field}{\mathbb{F}}

\newcommand{\calC}{\mathcal{C}}
\newcommand{\calD}{\mathcal{D}}
\newcommand{\calE}{\mathcal{E}}

\newcommand{\calI}{\mathcal{I}}
\newcommand{\calK}{\mathcal{K}}
\newcommand{\calL}{\mathcal{L}}

\newcommand{\calO}{\mathcal{O}}
\newcommand{\calP}{\mathcal{P}}
\newcommand{\calR}{\mathcal{R}}

\newcommand{\calU}{\mathcal{U}}

\newcommand{\Th}{\operatorname{Th}}
\newcommand{\pt}{\mathrm{pt}}
\newcommand{\freeduals}{\mathbb{D}}

\newcommand{\swapmatrix}{B}
\newcommand{\eventuallyconstantcat}{{(\prod_p \Vect^{\mathrm{f.d.}}_{\field_p})_\mathrm{ev.const.}}}

\newcommand{\cof}{\mathsf{cof}}
\newcommand{\cocart}{\mathrm{cocart}}
\newcommand{\fr}{\mathrm{fr}}
\newcommand{\fin}{\mathrm{fin}}
\newcommand{\fd}{\mathrm{f.d.}}
\newcommand{\stab}{\mathrm{stab}}
\newcommand{\gptrivsusp}{{\mathrm{gp},\Sigma\mathrm{-triv}}}
\newcommand{\trivsusp}{{\Sigma\mathrm{-triv}}}
\newcommand{\notgp}{{\neg \mathrm{gp}}}

\newcommand{\rex}{\mathrm{rex}}
\newcommand{\add}{\mathrm{add}}

\newcommand{\ho}{\mathsf{ho}}
\newcommand{\gp}{\mathrm{gp}}
\newcommand{\Psh}{\mathsf{Psh}}
\newcommand{\bbS}{\mathbb{S}}
\newcommand{\Idem}{\mathsf{Idem}}
\newcommand{\Split}{\mathrm{split}}
\newcommand{\divides}{\mid}
\newcommand{\notdivides}{\nmid}

\newcommand{\catdeloop}{\mathbb{B}}
\newcommand{\Tope}[1]{{#1\mathsf{Top}}}
\newcommand{\Spte}[1]{{#1\mathsf{Spt}}}

\tikzset{
            downward/.style={out=-90, in=90},
            upward/.style={out=90, in=-90},
            downcrossnw/.style={out=-90,in=135},
            downcrossne/.style={out=-90,in=45},
            crossswdown/.style={out=-135,in=90},
            crosssedown/.style={out=-45, in=90},
            capwdown/.style={out=180,in=90},
            capedown/.style={out=0,in=90},
            downcupw/.style={out=-90,in=180},
            downcupe/.style={out=-90,in=0},
            sw_w/.style={out=-135,in=180},
            sw_e/.style={out=-135,in=0},
            se_w/.style={out=-45,in=180},
            se_e/.style={out=-45,in=0},
            nw_w/.style={out=135,in=180},
            nw_e/.style={out=135,in=0},
            ne_w/.style={out=45,in=180},
            ne_e/.style={out=45,in=0},
            w_sw/.style={out=180,in=-135},
            e_sw/.style={out=0,in=-135},
            w_se/.style={out=180,in=-45},
            e_se/.style={out=0,in=-45},
            w_nw/.style={out=180,in=135},
            e_nw/.style={out=0,in=135},
            w_ne/.style={out=180,in=45},
            e_ne/.style={out=0,in=45},
            sw_sw/.style={out=-135,in=-135},
            sw_se/.style={out=-135,in=-45},
            sw_ne/.style={out=-135,in=45},
            sw_nw/.style={out=-135,in=135},
            se_sw/.style={out=-45,in=-135},
            se_se/.style={out=-45,in=-45},
            se_ne/.style={out=-45,in=45},
            se_nw/.style={out=-45,in=135},
            nw_sw/.style={out=135,in=-135},
            nw_se/.style={out=135,in=-45},
            nw_ne/.style={out=135,in=45},
            nw_nw/.style={out=135,in=135},
            ne_sw/.style={out=45,in=-135},
            ne_se/.style={out=45,in=-45},
            ne_ne/.style={out=45,in=45},
            ne_nw/.style={out=45,in=135},
            witharrow/.style={postaction = {
            decorate,decoration={markings,mark=at position 0.5 with {\arrow{>}}}}},
            withbackarrow/.style={postaction = {
            decorate,decoration={markings,mark=at position 0.5 with {\arrow{<}}}}},
            circlemor/.style={shape=circle, draw},
            rectmor/.style={shape=rectangle, draw},
            crossing/.style={},
            dummy/.style={},
            mycap/.style={},
            mycup/.style={}}

\begin{document}

\title[Free Duals]{Free Duals \protect\\ and \protect\\ A New Universal Property for Stable Equivariant Homotopy Theory}
\author{Timothy F. Campion}

\maketitle

\begin{abstract}
  We study the left adjoint $\mathbb{D}$ to the forgetful functor from the $\infty$-category of symmetric monoidal $\infty$-categories with duals and finite colimits to the $\infty$-category of symmetric monoidal $\infty$-categories with finite colimits, and related free constructions. The main result is that $\mathbb{D} \mathcal C$ always splits as the product of 3 factors, each characterized by a certain universal property. As an application, we show that, for any compact Lie group $G$, the $\infty$-category of genuine $G$-spectra is obtained from the $\infty$-category of Bredon (\emph{a.k.a} ``naive") $G$-spectra by freely adjoining duals for compact objects, while respecting colimits.
\end{abstract}

\tableofcontents

\section{Introduction}

This thesis\footnote{This is a pre-publication version of my 2021 PhD thesis at the University of Notre Dame, written under the supervision of Chris Schommer-Pries.} is a meditation on the implications that dualizable objects have for the global structure of symmetric monoidal $\infty$-category $(\calC,\wedge,S)$. We expose a number of such implications.\footnote{See \cref{def:dual} for a recollection of the notion of a \defterm{dual} for an object $X$ in a symmetric monoidal category $\calC$, familiar to homotopy theorists from the study of Spanier-Whitehead duality. This is the notion often called \defterm{strong dualizability} in the homotopy theory literature. We follow the convention of the category-theoretic literature in simply saying \defterm{dualizable}, because we do not use any other notion of ``dualizability" in this thesis.}
My ultimate goal is to apply such observations to understand what happens when one \emph{freely} adjoins duals to a symmetric monoidal $\infty$-category.

\subsection{Motivation: equivariant homotopy theory}
Motivating questions in this direction were asked by Charles Rezk (\cite{rezk}). Let $G$ be a compact Lie group (for example, $G$ might be finite). Roughly, Rezk asked:

\begin{enumerate}
    \item What happens when duals are freely adjoined to the $\infty$-category $\Tope{G}^\fin_\ast$ of finite pointed $G$-spaces?
    \item In particular, is the resulting category related to the category $\Spte G$ of genuine $G$-spectra?
\end{enumerate}

As a preliminary note, Rezk is asking about certain \emph{homotopical universal properties}, which immediately suggests that the appropriate framework to answer his question is the language of \emph{$\infty$-categories} (rather than e.g.\ the language of \emph{model categories}). For this reason, although much of the ``action" in this thesis takes place in the relevant homotopy categories, the main results are formulated and proven in the language of $\infty$-categories.

But perhaps Rezk's motivating questions themselves require some motivation. To this end, let us recall a couple of facts about equivariant homotopy theory with respect to a compact Lie group $G$. Unstably, there are, to a first approximation, two forms of $G$-equivariant homotopy theory. On the one hand there is \defterm{Borel $G$-equivariant homotopy theory}, which at the point-set level studies spaces with $G$-action up to \emph{underlying} weak homotopy equivalences. That is, in the Borel setting, a $G$-equivariant map is considered to be an ``equivalence" if it is a weak equivalence of underlying spaces. On the other hand there is \defterm{$G$-equivariant homotopy theory}, where a $G$-equivariant map is only considered to be a weak equivalence if it restricts to a weak homotopy equivalence on $H$-fixed-point sets for every closed subgroup $H \subseteq G$. Many natural equivariant questions originate in the Borel setting, but the $G$-equivariant setting provides more refined information and has better formal properties. Rezk's question concerns $G$-equivariant homotopy theory.

Having agreed to study equivariant homotopy theory rather than Borel equivariant homotopy theory, there is a further distinction to make when passing to the \emph{stable} arena of equivariant homotopy theory. On the one hand we have \defterm{Bredon} equivariant stable homotopy theory,\footnote{Bredon equivariant stable homotopy theory is often (somewhat pejoratively) called \emph{naive} equivariant stable homotopy theory.} where the ``$G$-spectra" involved admit deloopings with respect to $S^1$. On the other hand there is \defterm{genuine} equivariant stable homotopy theory, where the ``$G$-spectra" involved admit deloopings not just with respect to $S^1$, but with respect to every representation sphere $S^V$.\footnote{Here $V$ is a real $G$-representation, and the representation sphere $S^V$ is its one-point compactification.} Similar remarks are applicable here: Bredon equivariant homotopy theory is in some sense easier to define, but genuine equivariant homotopy theory has better formal properties, and it is essential to use this somewhat more refined setting for the important applications of equivariant homotopy theory.

Among the pleasant formal properties enjoyed by genuine equivariant homotopy theory is a good theory of \defterm{Spanier-Whitehead duality}. That is, just as in nonequivariant homotopy theory, it is the case that every finite genuine $G$-spectrum $X$ admits a dual $X^\vee$ (cf.\ \cref{lem:spt-g-duals}), called its ``Spanier-Whitehead dual". Precisely, the homotopy category of finite $G$-spectra, which is symmetric monoidal under smash product, has duals for all objects. By contrast, finite Bredon spectra are not dualizable in general.

Although not directly relevant to the present thesis, it is worth noting that this pattern is repeated elsewhere in homotopy theory. For example, if $S$ is a scheme, then the unstable motivic category $H(S)$ may be stabilized with respect to $S^1$ to obtain a stable category $SH^{S^1}(S)$, but better formal properties are obtained when one additionally stabilizes with respect to $\pjv^1$ and not just $S^1$, to obtain the stable motivic category $SH(S)$. Analogously to before, one of the pleasant formal properties of $SH(S)$ is that every smooth projective scheme over $S$ is dualizable in $SH(S)$, but not in $SH^{S^1}(S)$.

The general pattern appears to be something like the following. Given a ``geometric" $\infty$-category $\calC$, if one wishes to obtain a satisfactory ``stable version" of $\calC$, it is not enough to look for deloopings with respect to $S^1$. Instead, one should also look for any other ``sphere-like" objects in $\calC$, and ask for deloopings with respect to those objects as well. The story so far suggests that determining which objects are relevantly ``sphere-like", and deserve to be delooped, is some kind of art form, which encodes in an essential way certain information about the geometric nature of the ``sphere-like" objects of $\calC$. An attempt to formalize this general procedure which recurs both equivariantly and motivically might start by axiomatizing suitable properties for an object $C \in \calC$ to be regarded as a ``sphere", and then investigate what happens when such an object is inverted under the monoidal product. This paradigm has been studied in some generality by Robalo \cite{robalo}, following Voevodsky \cite{voevodsky}, who attributes it to Jeff Smith's unpublished study of Adams' construction of the smash product of spectra \cite{adams}. There, an object $T$ is ``sphere-like" if it is \defterm{symmetric}, in the sense that the cyclic permutation $T^{\wedge 3} \to T^{\wedge 3}$ (which, if the monoidal product $\wedge$ were cartesian, would be denoted $(x,y,z) \mapsto (z,x,y)$), should be homotopic to the identity (\cite[Definition 2.16]{robalo}). Robalo / Voevodsky / Smith show that if $T$ is symmetric, then the category $\calC[T\inv]$, by definition obtained from $\calC$ by universally inverting $T$ under the monoidal product, may be calculated in the familiar way as the sequential colimit $\calC_T$ of symmetric monoidal $\infty$-categories $\calC_T = \varinjlim(\calC \xrightarrow{T \wedge(-)} \calC \xrightarrow{T \wedge(-)}\cdots)$. The converse holds: if $\calC_T$ correctly computes $\calC[T\inv]$, then $T$ is symmetric (cf. \cite[Proposition 6]{nikolaus}). There is also an infinitary version: if $\calC$ is presentably symmetric monoidal, then let $\calC[T\inv,L]$  be the presentably symmetric monoidal $\infty$-category obtained from $\calC$ by universally inverting $T$ under $\wedge$ while respecting colimits. Then, if $T$ is symmetric, we have the inverse limit formula $\calC[T\inv,L] = \varprojlim(\calC \overset{(-)^T}{\leftarrow} \calC \overset{(-)^T}{\leftarrow} \cdots)$; conversely, if this formula is correct, then $T$ is symmetric. These sorts of theorems should be compared to ``group-completion" theorems in algebraic $K$-theory, such as that of McDuff-Segal \cite{mcduff-segal}.

At this point, one might conclude the following:
\begin{enumerate}
    \item A ``sphere-like" object should simply be defined to be a symmetric object.
    \item A good theory of ``genuine stabilization" should seek to identify such sphere-like objects, and monoidally invert them.
    \item We should consider ourselves lucky that there is a particularly simple formula for inverting an object when it is sphere-like.
\end{enumerate}
However, these conclusions are a bit question-begging: if we believe that this is the end of the story, then in some sense the ultimate reason for inverting sphere-like objects is simply the fact that we happen to know a convenient formula for doing so. We have not, for instance, identified a universal property associated with the inversion of sphere-like objects. We have not connected the inversion of sphere-like objects to the desideratum of good duality properties. And we have little to say about the significance of the nice formulas for these constructions other than that they are straightforward and familiar.

In this thesis, the logic is reversed. we propose understanding the phenomenon of ``genuine stabilization" not as a process of inverting carefully-chosen objects, but rather as a process of imposing a theory of Spanier-Whitehead duality for rather generically-chosen objects. That is, the use of genuine equivariant stable homotopy theory is often justified by pointing to its good theory of Spanier-Whitehead duality, and we propose to take this justification seriously. We show (\cref{cor:freecats}) that there is a well-defined way to \emph{freely adjoin duals} to the compact objects of a symmetric monoidal compactly-generated $\infty$-category. We show moreover, that when one does this to the $\infty$-category ${\Tope G}_\ast$ of pointed $G$-spaces, the resulting category $\freeduals_\omega {\Tope G}_\ast$ is \emph{almost} the $\infty$-category $\Spte G$ of genuine $G$-spectra; in fact $\freeduals_\omega {\Tope G}_\ast$ splits as a product of several subcategories; one factor (characterized as the unique stable factor) is none other than $\Spte G$ (\cref{cor:main-result}).

Thus the answers to Rezk's questions are that one can indeed freely adjoin duals to $G$-spaces in an appropriate sense, and that the resulting category is closely related to the category of $G$-spectra, but has some additional unstable ``cruft" attached to it. This ``cruft" is shaken off if we additionally stabilize our category, for instance, so that when the functor $\freeduals_\omega$ is applied to the symmetric monoidal $\infty$-category ${\Tope G}_\ast[(S^1)\inv]$ of Bredon $G$-spectra, the result is precisely the symmetric monoidal $\infty$-category $\Spte G$ of genuine $G$-spectra.

\subsection{Twisted-trivial braiding}

The proof of \cref{cor:main-result} does shed some light on the nature of ``sphere-like objects", insofar as the result hinges on the behavior of symmetric objects. In fact, for most of the thesis we work with the slightly stronger (\cref{prop:tt-smith}) condition that an object $T$ have \defterm{twisted-trivial braiding} (\cref{def:tt}), meaning that the braiding, or ``swap" map $\beta_{T,T}: T \wedge T \to T \wedge T$ is (homotopic to) a map of the form $\id_T \wedge t: T \wedge T \to T \wedge T$, for some $t : T \to T$.\footnote{If $T$ is additionally dualizable, then $T$ is symmetric iff it has twisted-trivial braiding (\cref{prop:smith-tt}). As we are generally forcing $T$ to become dualizable anyway, the distinction is not too important for us.} For instance (\cref{eg:twist-sphere}), when $T = S^n$ is a sphere in the pointed homotopy category $\ho\Top_\ast$, this condition is satisfied with $t = (-1)^n : S^n \to S^n$. Similarly, representation spheres $S^V$ in the homotopy category $\Tope{G}_\ast$ of pointed $G$-spaces have twisted-trivial braiding (\cref{eg:twist-G-sphere}), as does projective space $\pjv^1_S$ in the unstable pointed motivic category $H(S)$ (\cref{eg:twist-tate-sphere}). The crucial observation is that if $T$ has twisted-trivial braiding and is forced to become dualizable, then under mild conditions, the resulting category splits as a product of the subcategories of ``$T$-stable" objects and ``$T$-torsion objects" (\cref{prop:monoidal-splitting}). In the former factor, $T$ becomes invertible, while in the latter factor, $T$ becomes zero. Thus, if we take the property ``having twisted-trivial braiding" as a crude, approximate explication of what it means to be ``sphere-like", then \cref{prop:monoidal-splitting} may be stated loosely as follows: \emph{sphere-like objects, when dualizable, are either invertible or nil}.

In the case of equivariant homotopy theory, the representation spheres $S^V$ provide an ample supply of objects of $\Tope{G}^\fin_\ast$ with twisted-trivial braiding, each of which induces a splitting of the $\infty$-category $\freeduals^\rex \Tope{G}^\fin_\ast$ obtained by freely adjoining duals to $\Tope{G}^\fin_\ast$.
In the splitting induced by $S^V$, one factor (the ``$S^V$-stable" factor) has $S^V$ becoming invertible, while the other (the ``$S^V$-torsion" factor) has it becoming trivial. In order to arrive at the main result (\cref{cor:main-result}, \cref{cor:maincor}), the real input from equivariant homotopy theory consists in showing that at least if we restrict attention to the stable factor of $\freeduals^\rex \Tope{G}^\fin_\ast$, all of the $S^V$-torsion factors actually vanish, so that the stable factor of $\freeduals^\rex \Tope{G}^\fin_\ast$ has all representation spheres $S^V$ invertible. As $\Spte{G}$ is standardly defined (\cite{robalo}, \cite{gepner-meier}) to be obtained from $\Tope{G}$ by universally inverting the representation spheres $S^V$, we may now conclude that this stable factor of $\freeduals^\rex \Tope{G}^\fin_\ast$ is precisely the $\infty$-category $\Spte{G}^\fin$ of finite $G$-spectra. In other words, genuine equivariant homotopy theory is obtained from Bredon equivariant homotopy theory by universally adjoining duals for compact objects.

\subsection{Freely adjoining duals: finding the appropriate setting}\label{sec:appropriate}


We shall now explain what we mean by ``universally adjoining duals", previewing a bit from \cref{subsec:cob}. Recall that by the 1-dimesional cobordism hypothesis \cite{lurie-cob,harpaz-cob}, there is a free symmetric monoidal $\infty$-category on an object, which is given by the symmetric monoidal $\infty$-category $\Bord_1^\fr$ of oriented 0-manifolds and oriented bordisms between them. Thus if $X$ is an object of a symmetric monoidal $\infty$-category $\calC$, then a dual may be universally adjoined to the object $X$ by passing to the symmetric monoidal $\infty$-category $\calC\cup_{\FinBij} \Bord^\fr_1$, where $\FinBij$ is the symmetric monoidal groupoid of finite sets and the pushout is taken in the $\infty$-category $\SMC$ of symmetric monoidal $\infty$-categories. Moreover, because $\SMC$ is presentable, and because the full subcategory $\SMD \subset \SMC$ of symmetric monoidal $\infty$-categories with duals for objects is precisely the right orthogonal complement of the canonical symmetric monoidal functor $\FinBij \to \Bord^\fr_1$, it follows by the adjoint functor theorem that there is also a left adjoint $\freeduals: \SMC \to \SMD$ to the inclusion functor.

However, the functor $\freeduals$ is a bit too crude for the present purposes. Instead we work with variant categories such as the inclusion $\SMD^\rex \to \SMC^\rex$ of symmetric monoidal $\infty$-categories with duals and finite colimits into symmetric monoidal $\infty$-categories with finite colimits (preserved by the tensor product in each variable), and the reflection $\freeduals^\rex$ onto this full subcategory. After all, one is rarely interested in studying \emph{arbitrary} symmetric monoidal functors out of symmetric monoidal $\infty$-categories such as $\Tope{G}^\fin$, which has interesting finite colimits. The category $\Spte{G}^\fin$ which is supposed to be related to the output of the construction likewise has interesting finite colimits. Thus it is natural to take account of these colimits on both ends of the construction. Moroever, in order to invoke the full strength of the spitting theorem alluded to before (\cref{prop:monoidal-splitting}), it is necessary to know that the universal functor $\Tope{G}^\fin_\ast \to \freeduals^\rex \Tope{G}^\fin_\ast $ preserves cogroup objects, which are defined using coproducts. Furthermore, in order to have some reasonable interpretation of the object $S^1 \in \Tope{G}^\fin_\ast$ and its crucial twisted-trivial braiding, it is desireable for all functors in sight to preserve suspensions. For these reasons, we generally consider the functor $\freeduals^\rex$ and variants respecting other classes of colimits, rather than the plain functor $\freeduals$.


\subsection{Variations}

The main result comes in several versions, both finitary and infinintary (see \cref{cor:maincor}. For instance, $\Spte G$ is the free compactly-generated symmetric monoidal $\infty$-category on $\Tope G_\ast$ which has duals for compact objects. This results from the finitary version by applying the $\Ind$ construction to everything.

A word on basepoints: we believe that their use is inessential, because any symmetric monoidal $\infty$-category with duals for objects and an initial object is in fact pointed (\cref{prop:zero-duals}). So, modulo showing that $\Tope G^\fin_\ast$ is the free symmetric monoidal pointed $\infty$-category with finite colimits on $\Tope G^\fin$, viewed as a symmetric monoidal $\infty$-category with finite colimits, it will result that $\freeduals^\rex \Tope G^\fin = \freeduals^\rex \Tope G^\fin_\ast$, so that it should be possible to drop the basepoints from the above discussion. However, we have not at this point proven that $\Tope G^\fin_\ast$ has the required universal property. In fact, modulo similar considerations, it should be possible to derive a universal property for equivariant stable homotopy theory whose starting data is even simpler -- just the orbit category $\calO_G$ itself.

We hope to apply this methodology to also obtain a universal property for the stable motivic category $SH(S)$ as well. However, there the story is a bit more complicated: not every compact object is dualizable in general (though the smooth projective schemes provide a good supply of dualizable objects). Moreover, in this case, it appears that simply freely adjoining a dual to the object $\pjv^1$ and stabilizing does not cut down precisely to $SH(S)$ -- there seem to be other stable factors. We hope to study these factors in future work.

\subsection{Overview}

\cref{chap:1-cat} and \cref{chap:tt}, unlike the rest of this work, take place entirely in the 1-categorical setting. \cref{chap:1-cat} contains some standard background on (symmetric) (monoidal) 1-categories, none of which is new. The least familiar fact mentioned here may be Houston's theorem, \cref{prop:biproducts}, which says that any braided monoidal category with finite coproducts and duals is semiadditive. \cref{chap:tt} also takes place purely in the 1-categorical setting, but contains the most crucial results about split monoidal localizations, which are 1-categorical in nature. In \cref{def:tt}, we introduce the notion of an \emph{object with twisted-trivial braiding}, which is closely related to the \emph{symmetry} condition (\cref{def:smith}), used by Voevodsky, \cite[Definition 2.16]{robalo}, \cite[Definition 9.2]{hovey}, \cite{nikolaus}, and others to control monoidal localizations. We show (\cref{prop:clopens}) that any dualizable object with twisted-trivial braiding gives rise to a smashing-cosmashing localization. We review in \cref{sec:comp-smash} how this often implies that the entire category splits along lines defined by this object.

Next, these 1-categorical results are applied in the $\infty$-categorical setting. \cref{chap:oo-cat} develops the $\infty$-categorical infrastructure necessary to apply these results in the $\infty$-categorical setting, primarily recalling material from \cite{ha}. \cref{chap:dual-colim} explores certain splittings which occur in great generality when duals and colimits interact. \cref{chap:app} applies the theory to the example of equivariant homotopy theory, leading up to the main theorem, \cref{cor:maincor}.

\subsection{Acknowledgements}
This document is a pre-publication version of my PhD thesis, completed at the University of Notre Dame between 2016 and 2021. Minor additions and corrections have been made since then. I would like to thank my advisor Chris Schomer-Pries, for helping me to grow and mature as a mathematician, and for supporting me in challenging times. Thanks also to the members of my dissertation committe: Chris Schommer-Pries, Mark Behrens, Stephan Stolz, and Pavel Mnev -- each of whom have been generous with their time and support during my time at Notre Dame. On that note, I would like to thank all the wonderful people I have worked with at Notre Dame during my time there, particularly the members of the Geometry and Topology group, I would like to thank my family for supporting me through challenging times. I am grateful for the support of NSF grant DMS-1547292, which supported some of this work, as well as the support of the ARO under MURI Grant W911NF-20-1-0082.

\section{1-Categorical Preliminaries}\label{chap:1-cat}



This chapter discusses several elementary structures encountered in 1-categories and monoidal 1-categories. In contrast to the rest of the text, no $\infty$-categories appear in this chapter. Most of the material comprises a review of some basic concepts. However, the central idea of this thesis also appears in \cref{sec:tt} and \cref{sec:tt-dual}, namely the attention given to objects with \defterm{twisted-trivial braiding} (\cref{def:tt}) in the context of dualizability. The main theorem on these objects is \cref{prop:clopens}, which in combination with \cref{prop:ordinary-split} shows that dualizable objects with twisted-trivial braidings often lead to splittings of the entire category. The implications of these observations, when applied to various ``sphere-like" examples of such objects (see \cref{sec:tt}) , are the main topic of this thesis, and will be explored further in \cref{chap:dual-colim} and \cref{chap:app}.

One perhaps well-known but under-publicized fact is discussed in \cref{sec:semiadd}, where in \cref{prop:biproducts} a result of Robin Houston is recalled \cite{houston}: any symmetric monoidal category with duals and finite coproducts is semiadditive.

Throughout this chapter, we work usually with braided monoidal categories. This is not because we have particular examples in mind which are braided but not symmetric, but rather because we find the concomitant string diagrams easier to understand in the braided setting than in the symmetric setting.

\subsection{Monoidal categories}\label{sec:moncat}

We assume that the reader is familiar with basic concepts of category theory, as well as monoidal categories, braided monoidal categories, and symmetric monoidal categories. We will often take advantage of the Mac Lane coherence theorem to assume that monoidal categories are \defterm{strict}, i.e. that the associators and unitors are identities. We shall also freely use string diagrams to reason in monoidal categories and braided or symmetric monoidal categories. We assume the reader is familiar with their usage. My conventions are that string diagrams are read from top to bottom, and that the braiding isomorphism is depicted as follows:

\begin{equation*}
    \beta_{X,Y} = 
    \begin{tikzpicture}[baseline=(C.south)]
        
        \node (Xtop) at (0,0) {$X$};
        \node (Ytop) at (1,0) {$Y$};
        \node (C) at (.5,-1) {};

        \node (Ybot) at (0,-2) {$Y$};
        \node (Xbot) at (1,-2) {$X$};
        

        \draw (Xtop) [downcrossnw] to (C.center);
        \draw (Ytop) [downcrossne] to (C.north east);
        \draw (C.center) [crosssedown] to (Xbot);
        \draw (C.south west) [crossswdown] to (Ybot);
    \end{tikzpicture}
\end{equation*}

As part of my use of string diagrams, we assume the reader is familiar with the use of ``cups" and ``caps" associated with dual objects (the definition of a dual object is recalled in \cref{def:dual} below). When considering morphisms involving both an object and its dual in a string diagram, we include an arrow indicating the orientation of the string. If the orientation is downward for an object $X$, then it is upward for $X^\dual$, etc.

\begin{notation}
There is an equivalence between monoidal categories and strong mon\-oid\-al functors on the one hand, and bicategories with a unique object on the other hand. If $\calC$ is a monoidal category, denote by $\catdeloop \calC$ the corresponding 1-object bicategory.
\end{notation}

\begin{Def}
Let $\calK$ be a class of small categories. A \defterm{monoidal category with $\calK$-colimits} is a monoidal category $(\calC,\wedge,S)$ such that $\calC$ has $\calK$-colimits and $\wedge$ preserves them separately in each variable. A braided (resp. symmetric) monoidal category with $\calK$-colimits is a braided (resp. symmetric) monoidal category whose underlying monoidal category is a monoidal category with $\calK$-colimits.
\end{Def}

\begin{rmk}
See \cref{sec:colims} for a more systematic account of the interaction of monoidal structure and colimits in the $\infty$-categorical setting.
\end{rmk}

\begin{Def}
A \defterm{localization} is an adjunction $L: \calC^\to_\leftarrow \calD: i$ whose counit $Li \Rightarrow \id_\calD$ is an isomorphism. A \defterm{monoidal localization} is a localization $L \dashv i$ where $L$ is strong monoidal. Braided (resp. symmetric) monoidal localizations are defined in the obvious way.

A \defterm{smashing localization} is a braided monoidal localization $L: \calC^\to_\leftarrow \calD: i$ such that the map $i(S) \wedge X \to iL(X)$ is an isomorphism for all $X \in \calC$ (here $S$ is the unit of $\calD$ and $\wedge$ is the tensor in $\calC$). 
\end{Def}

We will have more to say about smashing localizations in \cref{sec:idem-smash} below.

\subsection{Duals}\label{sec:dual}
In this section we review the concept of dualizability in monoidal and particularly in braided monoidal categories.

\begin{Def}\label{def:inv}
Let $(\calC,S,\wedge)$ be a monoidal category. An object $C \in \calC$ is said to be \defterm{invertible} if $C$ is an equivalence in $\catdeloop \calC$. Explicitly, this means that there exists an object $C\inv \in \calC$ such that $C\inv \wedge C \cong S$ and $C \wedge C\inv \cong S$. We call $C\inv$ the \defterm{inverse} of $C$.
\end{Def}

\begin{rmk}
If $\calC$ is braided, then the two isomorphisms $C\inv \wedge C \cong S$ and $C \wedge C\inv \cong S$ are equivalent, so only one of them need be checked.
\end{rmk}

\begin{rmk}
Let $(\calC,S,\wedge)$ be a monoidal category. An inverse to $C \in \calC$, if it exists, is unique up to isomorphism, justifying the notation $C\inv$ used in \cref{def:inv}.
\end{rmk}

\begin{rmk}
Strong monoidal functors preserve invertible objects and inverses to invertible objects.
\end{rmk}

\begin{eg}
The unit object of a monoidal category is inverse to itself.
\end{eg}

\begin{Def}\label{def:dual}
Let $(\calC,S,\wedge)$ be a monoidal category. A \defterm{duality datum} in $\calC$ consists of an adjunction in the associated bicategory $\catdeloop \calC$. Explicitly, this means we have a pair of objects $R,L$ and morphisms $\eta: S \to  R \wedge L$ (the \defterm{unit}), $\varepsilon: L \wedge R \to S$ (the \defterm{counit}) satisfying the \defterm{triangle equations}: $(\varepsilon \wedge L)(L\wedge \eta) = \id_L$ and $ \id_R = (R \wedge \varepsilon)(\eta \wedge R)$. We say that $L$ is the \defterm{left dual} and $R$ is the \defterm{right dual}, or that \defterm{$L$ is left dual to $R$} and \defterm{$R$ is right dual to $L$}. If $\calC$ is braided then left and right duals coincide, so we drop the handedness from the terminology and write $L = R^\dual$, $R = L^\dual$.
\end{Def}

\begin{rmk}
Let $L$ be an object of a monoidal category $\calC$. Then the groupoid of extensions of $L$ to a duality datum where $L$ is the left dual is either empty or contractible. Thus admitting a right dual,like admitting a right adjoint, is a property of an object rather than a structure. This fact justifies the notation $X^\dual$ introduced in \cref{def:dual}.
\end{rmk}

\begin{rmk}
Strong monoidal functors preserve duals in the sense that if $(L,R,\eta,\varepsilon)$ is a duality datum and $F$ is a strong monoidal functor, then $(F(L),F(R),F(\eta),F(\varepsilon))$ is also a dualtiy datum.
\end{rmk}

\begin{rmk}
In some of the algebraic topology literature, the term ``strongly dualizable" is used for what we have called ``dualizable" in \cref{def:dual}. There is some justification for this terminology: if $(\calC,\wedge,S)$ is monoidal with internal homs, then the dual of $X$ must, by uniqueness of adjoints, coincide with the internal hom $[X,S]$, so there is some justification for referring to $[X,S]$ as a ``dual" which always exists in the presence of internal homs, and reserving ``strong dual" for the case when a unit map exists to complete the duality data. However, as we do not consider here the internal hom $[X,S]$ in cases where $X$ is not (strongly) dualizable, we have opted to use the plain term ``dualizable" for the situation of \cref{def:dual}.
\end{rmk}

\begin{eg}
Any invertible object is dualizable. This is a special case of the fact that any equivalence in a 2-category may be upgraded to an adjoint equivalence.
\end{eg}

\subsection{Pointed categories}\label{sec:pointed}
In this section we review basic concepts about categories with a zero object.

\begin{Def}
A category $\calC$ is \defterm{pointed} if it has an object $0 = 0_\calC \in \calC$ which is both initial and terminal. A \defterm{pointed functor} between pointed categories is a functor preserving the zero object. A \defterm{monoidal pointed category} is a monoidal category with initial object whose underlying category is pointed. A braided (resp. symmetric) monoidal pointed categories is a braided (resp. symmetric) monoidal category whose underlying monoidal category is a monoidal pointed category.
\end{Def}


\begin{rmk}
Let $\calC$ be an $\infty$-category with an initial object $\emptyset$ and a terminal object $1$. Then $\calC$ is pointed if and only if the unique map $\emptyset \to 1$ is invertible.
\end{rmk}

A folklore result is that
\begin{prop}\label{prop:zero-duals}
Let $(\calC,\wedge,S)$ be a monoidal category with an initial object $\emptyset$ preserved by $X \wedge(-)$ for any object $X$. Then $\emptyset$ has a right dual $\emptyset^\dual$ iff $\emptyset = \emptyset^\dual$ is a zero object.
\end{prop}
\begin{proof}
In one direction, it is straightforward to check that a zero object is canonically self-dual. Conversely, suppose that $\emptyset$ has a dual $\emptyset^\dual$ and $\wedge$ preserves the initial object. One of the triangle equations says that $\emptyset^\dual$ is a retract of $\emptyset \wedge \emptyset^\dual \wedge \emptyset = \emptyset$ and hence $\emptyset^\dual \cong \emptyset$. Then maps $X \to \emptyset$ correspond by adjunction to maps $\emptyset = X \wedge \emptyset^\dual \to S$, so that $\emptyset$ is terminal as well as initial as desired.
\end{proof}

\subsection{Semiadditive categories}\label{sec:semiadd}
In this section we review basic concepts about categories enriched in commutative monoids.

\begin{Def}
Let $X,Y,Z$ be objects of a category $\calC$. We say that $Z$ is the (binary) \defterm{biproduct} or \defterm{direct sum} of $X$ and $Y$, and write $Z = X \oplus Y$, if there is a diagram $X^\leftarrow_\to Z^\to_\leftarrow Y$ which exhibits $X,Y$ as splittings of commuting idempotents, while simultaneously exhibiting $Z$ as both the product and coproduct of $X$ and $Y$. We say that $X,Y$ are \defterm{complementary} retracts of $Z$. If every $X,Y$ fit into a biproduct diagram, we say that \defterm{$\calC$ has binary biproducts}, and if in addition $\calC$ has a zero object, we say that \defterm{$\calC$ has finite biproducts} or \defterm{$\calC$ is semiadditive}.

If $\calC$ is monoidal, then a biproduct is \defterm{monoidal} if it is preserved by $X \wedge (-)$ and $(-) \wedge X$ for each $X \in \calC$.
A \defterm{monoidal semiadditive category} is a monoidal category with finite coproducts whose underlying category is semiadditive. Braided (resp. symmetric) monoidal semiadditive categories are defined in the obvious way.
\end{Def}

\begin{rmk}
Suppose that $Z = X \oplus Y$ in a category $\calC$. Let $e,f$ be the idempotents on $Z$ split by $X,Y$ respectively. Then $ef$ is also an idempotent, and its splitting $W$ (which exists in the idempotent completion of $\calC$) is both subterminal (in the sense that any object admits at most one map to $W$) and co-subterminal (any object admits at most one map from $W$). Thus if $\calC$ has a zero object, then $W$ is the zero object.
\end{rmk}

\begin{rmk}
Suppose that $X$ is a retract of $Z$, and that $Y,Y'$ are retracts of $Z$ which are both complementary to $X$. Then $Y,Y'$ are equal as subobjects (or quotient objects) of $Z$. Thus we may speak of \emph{the} complement of a retract of $Z$, if it exists.
\end{rmk}

\begin{rmk}
Let $\calC$ be a category with finite biproducts.
Then $\calC$ is uniquely enriched in commutative monoids. Conversely, if $\calC$ is enriched in commutative monoids and has finite products or finite coproducts, then $\calC$ is semiadditive.
\end{rmk}

\begin{rmk}
If $\calC$ is pointed and $X,Y \in \calC$, then the biproduct of $X$ and $Y$ exists if and only if the coproduct $X \amalg Y$ and the product $X \times Y$ exist and the canonical map $X \amalg Y \to X \times Y$ is an isomorphism.
\end{rmk}

In analogy to \cref{prop:zero-duals}, we have another well-known result:

\begin{prop}[\cite{houston}, see also \cite{garner-schappi}]\label{prop:biproducts}
Let $(\calC,\wedge,S)$ be a braided monoidal category with zero object. Suppose that the product $S \times S$ exists and let $X,Y \in \calC$ be objects such that the coproduct $X \amalg  Y$ exists and moreover $X \wedge(-)$, $Y \wedge(-)$, and $(X \amalg  Y) \wedge(-)$ all preserve the product $S \times S$, while $(-) \wedge (S \times S)$ preserves the coproduct $X \amalg  Y$. Then $X \amalg  Y$ is a biproduct $X \oplus Y$.
\end{prop}
\begin{rmk}
In particular, this result and its dual imply that if $(\calC, \wedge, S)$ is a braided monoidal category with finite coproducts preserved by $\wedge$ in each variable, then the biproduct $X\oplus Y$ exists in any of the following cases:
\begin{enumerate}
    \item $S \times S$ exists and $X,Y$ have duals.
    \item $S \amalg  S$ has a dual.
    \item Every object has a dual.
\end{enumerate}
\end{rmk}
\begin{proof}[Proof of Proposition \ref{prop:biproducts}]
By the hypotheses, we have
\begin{align*}
    (X \times X)\amalg (Y \times Y) &= ((X \wedge S) \times (X \wedge S))\amalg ((Y \wedge S) \times (Y \wedge S)) \\
                              &= ( X \wedge (S \times S)) \amalg  (Y \wedge (S \times S)) \\
                              &= (X \amalg  Y) \wedge (S \times S) \\
                              &= ((X\amalg Y) \wedge S) \times ((X\amalg Y) \wedge S) \\
                              &= ((X \wedge S)\amalg (Y \wedge S)) \times ((X \wedge S)\amalg (Y \wedge S)) \\
                              &= (X\amalg Y)\times (X\amalg Y)
\end{align*}
(To be careful, we have that $(X\amalg Y) \wedge (S \times S)$ exists and we derive the existence of the other objects from the hypotheses). The canonical map $X \amalg  Y \to X \times Y$ whose invertibility defines semiadditivity is a retract of this canonical isomorphism, and is thus also an isomorphism.
\end{proof}

\subsection{(Co)group Objects}\label{sec:cogp}
In this section we review basic concepts about group objects and cogroup objects.

\begin{Def}
A \defterm{group object} in a category $\calC$ with finite products is an object $X$ with maps $1 \xrightarrow u X \overset{m}{\leftarrow} X \times X$ (unit and multiplication) and a \defterm{negation} map $(-1): X \to X$ satisfying the group equations. Dually, we define a \defterm{cogroup object} in a category with finite coproducts.

An \defterm{additive} category is a semiadditive category where every object is a group object, or (equivalently -- see \cref{rmk:semiadd-cogp}) a cogroup object. A \defterm{monoidal additive category} is a monoidal semiadditive category whose underlying category is additive. Braided (resp. symmetric) monoidal additive categories are defined in the obvious way.
\end{Def}

\begin{rmk}\label{rmk:semiadd-cogp}
If $\calC$ is semiadditive, then group objects and cogroup objects are the same thing, and an object can be a (co)group in at most one way. An object $X$ is a (co)group object if and only if the commutative monoid $\Hom(X,Y)$ is an abelian group for all $Y$, if and only if the commutative monoid $\Hom(Y,X)$ is an abelian group for all $Y$, if and only if the commutative monoid $\Hom(X,X)$ is an abelian group. Another way to say this is that the ``shear" map $\begin{pmatrix} \id_X & \id_X \\ 0 & \id_X \end{pmatrix}: X \oplus X \to X \oplus X$ is an isomorphism.
\end{rmk}

\begin{rmk}
(Co)group objects are preserved by any functor which preserves finite (co)products.
\end{rmk}

\begin{rmk}\label{rmk:cogp-cof-of-diag}
Suppose that $X$ is a cogroup object, and moreover that the cofiber $C$ of the comultiplication $X \to X \amalg  X$ exists. Then $C$ corepresents invertible elements of $\Hom(X,-)$.
\end{rmk}

\begin{eg}\label{eg:cogp-sphere}
The sphere $S^1$ is a cogroup object in $\ho\Top_\ast$, the homotopy category of pointed spaces.
\end{eg}

\begin{eg}\label{eg:cogp-susp}
Let $\calC$ be a pointed $\infty$-category with suspensions and finite coproducts. Then it follows from \cref{eg:cogp-sphere} that the suspension $\Sigma X$ of any object $X$ is a cogroup object in $\ho\calC$.
\end{eg}

Cogroup structure provides one way to ensure that complements exist:
\begin{lem}\label{lem:cogroup}
Let $\calC$ be a pointed category, and suppose that $X^{\overset{r}{\leftarrow}}_{\underset{i}{\to}} Z$ is a retract. If a complement $Z^{\overset{q}{\to}}_{\underset{j}{\leftarrow}}Y$ to $X$ exists, then $q$ exhibits $Y$ as the cokernel of $i$ (and dually, $j$ exhibits $Y$ as the kernel of $p$). Conversely, if $\calC$ is semiadditive and $X$ is a group object, then any cokernel of $i$ (or kernel of $p$) is a complement to $X$. Alternatively, an idempotent splitting for $1 + i(-1)r$ yields a complementary retract of $Z$.
\end{lem}
\begin{proof}
Only the ``Conversely" part is not entirely standard. For this, use the negation on $X$ to construct the idempotent $1 + i(-1)p$ on $Z$ and then proceed as usual.
\end{proof}

\section{Structures in 1-Category Theory}\label{chap:tt}

\subsection{Objects with (twisted) trivial braiding and the symmetry condition}\label{sec:tt}
In this section we discuss objects with \defterm{twisted trivial braiding} (\cref{def:tt}), and the relationship to the \defterm{symmetry condition} (\cref{def:smith}). The latter is well-known from the study of monoidal localization; the former is stronger in general (\cref{prop:tt-smith}), but sometimes easier to check. Moreover, the two conditions become equivalent as soon as the object in question is dualizable (\cref{prop:smith-tt}). We shall see later, for instance in \cref{sec:tt-dual}, that these conditions continue to be relevant when one is adjoining a $\wedge$-dual to an object rather than a $\wedge$-inverse.

\begin{prop}\label{prop:twisted-equiv}
Let $(\calC,\wedge,S)$ be a braided monoidal category. Let $T$ be an object of $\calC$, and let $s,t: T \rightrightarrows T$ be endomorphisms of $T$. Consider the braiding morphism $\beta_{T,T}: T \wedge T \to T \wedge T$. The following are equivalent:
\begin{enumerate}
    \item $\beta_{T,T} = s \wedge t$
    \item $\beta_{T,T} = t \wedge s$
    \item $\beta_{T,T} = ts \wedge \id_T$
    \item $\beta_{T,T} = st \wedge \id_T$
    \item $\beta_{T,T} = \id_T \wedge ts$
    \item $\beta_{T,T} = \id_T \wedge st$
\end{enumerate}
\end{prop}
\begin{proof}
First we show that $(1) \Rightarrow (2)$. Assuming $(1)$, we calculate:

\begin{equation*}
    \begin{tikzpicture}[baseline=(base.south)]
        \node (base) at (0,-2) {};
    
        \node (T0top) at (0,0) {$T$};
        \node (T1top) at (1,0) {$T$};

        \node (T0bot) at (0,-4) {$T$};
        \node (T1bot) at (1,-4) {$T$};
        
        
        \draw (T0top) [downward] to (T0bot);
        \draw (T1top) [downward] to (T1bot);
    \end{tikzpicture}
    =
    \begin{tikzpicture}[baseline=(base.south)]
        \node (base) at (0,-2) {};
        
        \node (T0top) at (0,0) {$T$};
        \node (T1top) at (1,0) {$T$};
        \node (C0) at (.5,-1) {};
        \node (C1) at (.5,-3) {};

        \node (T0bot) at (0,-4) {$T$};
        \node (T1bot) at (1,-4) {$T$};
        
        
        \draw (T0top) [downcrossnw] to (C0.center);
        \draw (T1top) [downcrossne] to (C0.north east);
        \draw (C0.center) [se_ne] to (C1.center);
        \draw (C0.south west) [sw_nw] to (C1.north west);
        \draw (C1.center) [crossswdown] to (T0bot);
        \draw (C1.south east) [crosssedown] to (T1bot);
    \end{tikzpicture}
    =
    \begin{tikzpicture}[baseline=(base.south)]
        \node (base) at (0,-2) {};
        
        \node (T0top) at (0,0) {$T$};
        \node (T1top) at (1,0) {$T$};
        \node [circlemor] (S) at (0,-1.5) {$s$};
        \node [circlemor] (T) at (1,-1.5) {$t$};
        \node (C1) at (.5,-3) {};

        \node (T0bot) at (0,-4) {$T$};
        \node (T1bot) at (1,-4) {$T$};
        
        
        \draw (T0top) [downward] to (S);
        \draw (T1top) [downward] to (T);
        \draw (T) [downcrossne] to (C1.center);
        \draw (S) [downcrossnw] to (C1.north west);
        \draw (C1.center) [crossswdown] to (T0bot);
        \draw (C1.south east) [crosssedown] to (T1bot);
    \end{tikzpicture}
=
\begin{tikzpicture}[baseline=(base.south)]
        \node (base) at (0,-2) {};
        
        \node (T0top) at (0,0) {$T$};
        \node (T1top) at (1,0) {$T$};
        \node [circlemor] (S) at (1,-2.5) {$s$};
        \node [circlemor] (T) at (0,-2.5) {$t$};
        \node (C1) at (.5,-1) {};

        \node (T0bot) at (0,-4) {$T$};
        \node (T1bot) at (1,-4) {$T$};
        
        
        \draw (T) [downward] to (T0bot);
        \draw (S) [downward] to (T1bot);
        \draw (T1top) [downcrossne] to (C1.center);
        \draw (T0top) [downcrossnw] to (C1.north west);
        \draw (C1.center) [crossswdown] to (T);
        \draw (C1.south east) [crosssedown] to (S);
    \end{tikzpicture}
\end{equation*}
Here in the first equation we have isotoped, so that in the second equation we may apply the identity $\beta_{T,T} = s \wedge t$ from $(1)$. In the third equation we have isotoped. This shows that $\id_{T \wedge T} = (t \wedge s) \beta_{T,T}\inv$. Because the braiding is invertible, this implies that $t \wedge s = \beta_{T,T}$, i.e. that $(2)$ holds. Of course, $(2) \Rightarrow (1)$ follows by symmetry. This also implies that $(3) \Leftrightarrow (5)$ and $(2) \Leftrightarrow (4)$.

We now show that $(1) \Rightarrow (5)$. Assuming $(1)$, we calculate:
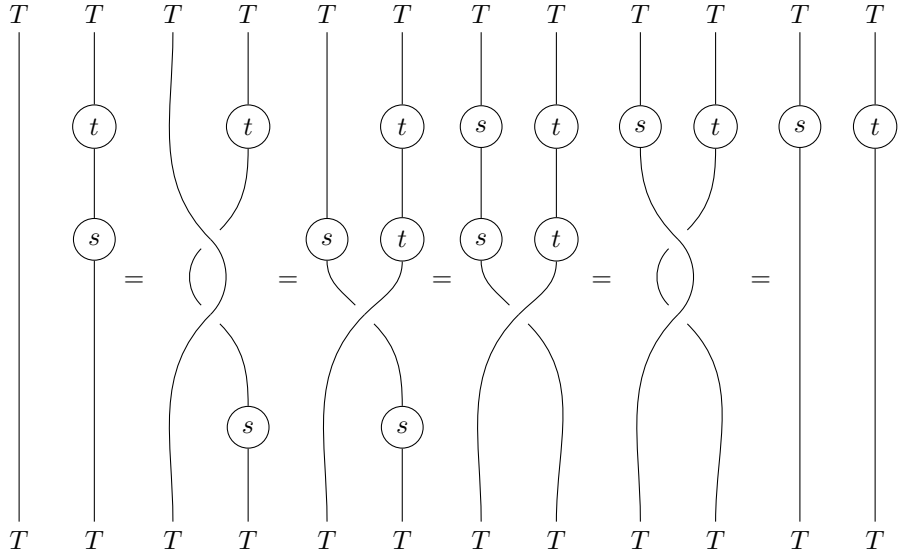
\begin{equation*}
    \begin{tikzpicture}[baseline=(base.south)]
        \node (base) at (0,-3.5) {};
        
        \node (T0top) at (0,0) {$T$};
        \node (T1top) at (1,0) {$T$};
        \node [circlemor] (T) at (1,-1.5) {$t$};
        \node [circlemor] (S) at (1,-3) {$s$};

        \node (T0bot) at (0,-7) {$T$};
        \node (T1bot) at (1,-7) {$T$};
        
        
        \draw (T0top) [downward] to (T0bot);
        \draw (T1top) [downward] to (T);
        \draw (T) [downward] to (S);
        \draw (S) [downward] to (T1bot);
    \end{tikzpicture}
    =
    \begin{tikzpicture}[baseline=(base.south)]
        \node (base) at (0,-3.5) {};
        
        \node (T0top) at (0,0) {$T$};
        \node (T1top) at (1,0) {$T$};
        \node [circlemor] (T) at (1,-1.5) {$t$};
        \node (C0) at (.5,-3) {};
        \node (C1) at (.5,-4) {};
        \node [circlemor] (S) at (1,-5.5) {$s$};

        \node (T0bot) at (0,-7) {$T$};
        \node (T1bot) at (1,-7) {$T$};
        
        
        \draw (T0top) [downcrossnw] to (C0.center);
        \draw (T1top) [downward] to (T);
        \draw (T) [downcrossne] to (C0.north east);
        \draw (C0.center) [se_ne] to (C1.center);
        \draw (C0.south west) [sw_nw] to (C1.north west);
        \draw (C1.center) [crossswdown] to (T0bot);
        \draw (C1.south east) [crosssedown] to (S);
        \draw (S) [downward] to (T1bot);
    \end{tikzpicture}
    =
    \begin{tikzpicture}[baseline=(base.south)]
        \node (base) at (0,-3.5) {};
        
        \node (T0top) at (0,0) {$T$};
        \node (T1top) at (1,0) {$T$};
        \node [circlemor] (T) at (1,-1.5) {$t$};
        \node (SS) [circlemor] at (0,-3) {$s$};
        \node (TT) [circlemor] at (1,-3) {$t$};
        \node (C1) at (.5,-4) {};
        \node [circlemor] (S) at (1,-5.5) {$s$};

        \node (T0bot) at (0,-7) {$T$};
        \node (T1bot) at (1,-7) {$T$};
        
        
        \draw (T0top) [downward] to (SS);
        \draw (T1top) [downward] to (T);
        \draw (T) [downward] to (TT);
        \draw (SS) [downcrossnw] to (C1.north west);
        \draw (TT) [downcrossne] to (C1.center);
        \draw (C1.center) [crossswdown] to (T0bot);
        \draw (C1.south east) [crosssedown] to (S);
        \draw (S) [downward] to (T1bot);
    \end{tikzpicture}
    =
    \begin{tikzpicture}[baseline=(base.south)]
        \node (base) at (0,-3.5) {};
        
        \node (T0top) at (0,0) {$T$};
        \node (T1top) at (1,0) {$T$};
        \node [circlemor] (T) at (1,-1.5) {$t$};
        \node (SS) [circlemor] at (0,-1.5) {$s$};
        \node (TT) [circlemor] at (1,-3) {$t$};
        \node [circlemor] (S) at (0,-3) {$s$};
        \node (C1) at (.5,-4) {};

        \node (T0bot) at (0,-7) {$T$};
        \node (T1bot) at (1,-7) {$T$};
        
        
        \draw (T0top) [downward] to (SS);
        \draw (T1top) [downward] to (T);
        \draw (T) [downward] to (TT);
        \draw (SS) [downward] to (S);
        \draw (S) [downcrossnw] to (C1.north west);
        \draw (TT) [downcrossne] to (C1.center);
        \draw (C1.center) [crossswdown] to (T0bot);
        \draw (C1.south east) [crosssedown] to (T1bot);
    \end{tikzpicture}
    =
    \begin{tikzpicture}[baseline=(base.south)]
        \node (base) at (0,-3.5) {};
        
        \node (T0top) at (0,0) {$T$};
        \node (T1top) at (1,0) {$T$};
        \node [circlemor] (T) at (1,-1.5) {$t$};
        \node (SS) [circlemor] at (0,-1.5) {$s$};
        \node (C0) at (.5,-3) {};
        \node (C1) at (.5,-4) {};

        \node (T0bot) at (0,-7) {$T$};
        \node (T1bot) at (1,-7) {$T$};
        
        
        \draw (T0top) [downward] to (SS);
        \draw (T1top) [downward] to (T);
        \draw (T) [downcrossne] to (C0.north east);
        \draw (SS) [downcrossnw] to (C0.center);
        \draw (C0.south west) [sw_nw] to (C1.north west);
        \draw (C0.center) [se_ne] to (C1.center);
        \draw (C1.center) [crossswdown] to (T0bot);
        \draw (C1.south east) [crosssedown] to (T1bot);
    \end{tikzpicture}
    =
    \begin{tikzpicture}[baseline=(base.south)]
        \node (base) at (0,-3.5) {};
        
        \node (T0top) at (0,0) {$T$};
        \node (T1top) at (1,0) {$T$};
        \node [circlemor] (T) at (1,-1.5) {$t$};
        \node (SS) [circlemor] at (0,-1.5) {$s$};

        \node (T0bot) at (0,-7) {$T$};
        \node (T1bot) at (1,-7) {$T$};
        
        
        \draw (T0top) [downward] to (SS);
        \draw (T1top) [downward] to (T);
        \draw (T) [downward] to (T1bot);
        \draw (SS) [downward] to (T0bot);
    \end{tikzpicture}
\end{equation*}
Here in the first equation we have isotoped, so that in the second equation we may apply the identity $\beta_{T,T} = s \wedge t$ of $(1)$. In the third equation we have isotoped so that in the fourth equation we may apply the identity $\beta_{T,T} = s \wedge t$ again. In the fifth equation we have isotoped. We are left with $\id_T \wedge ts = s \wedge t$, which by hypothesis is the same as $\beta_{T,T}$, so that $(5)$ holds as claimed.

We now show that $(5) \Rightarrow (1)$. Assuming $(5)$, we calculate:
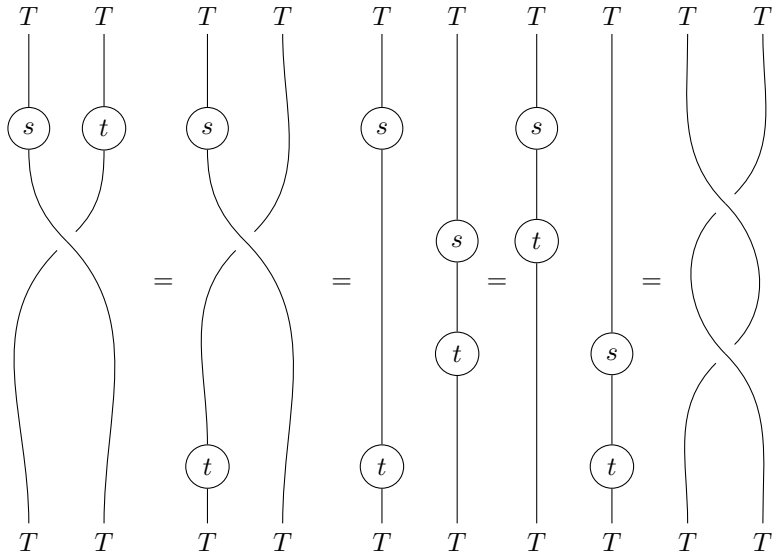
\begin{equation*}
    \begin{tikzpicture}[baseline=(base.south)]
        \node (base) at (0,-3.5) {};
        
        \node (T0top) at (0,0) {$T$};
        \node (T1top) at (1,0) {$T$};
        \node [circlemor] (T) at (1,-1.5) {$t$};
        \node (S) [circlemor] at (0,-1.5) {$s$};
        \node (C) at (.5,-3) {};

        \node (T0bot) at (0,-7) {$T$};
        \node (T1bot) at (1,-7) {$T$};
        
        
        \draw (T0top) [downward] to (S);
        \draw (T1top) [downward] to (T);
        \draw (T) [downcrossne] to (C.north east);
        \draw (S) [downcrossnw] to (C.center);
        \draw (C.center) [crosssedown] to (T1bot);
        \draw (C.south west) [crossswdown] to (T0bot);
    \end{tikzpicture}
    =
    \begin{tikzpicture}[baseline=(base.south)]
        \node (base) at (0,-3.5) {};
        
        \node (T0top) at (0,0) {$T$};
        \node (T1top) at (1,0) {$T$};
        \node [circlemor] (T) at (0,-6) {$t$};
        \node (S) [circlemor] at (0,-1.5) {$s$};
        \node (C) at (.5,-3) {};

        \node (T0bot) at (0,-7) {$T$};
        \node (T1bot) at (1,-7) {$T$};
        
        
        \draw (T0top) [downward] to (S);
        \draw (T1top) [downcrossne] to (C.north east);
        \draw (S) [downcrossnw] to (C.center);
        \draw (C.center) [crosssedown] to (T1bot);
        \draw (C.south west) [crossswdown] to (T);
        \draw (T) [downward] to (T0bot);
    \end{tikzpicture}
    =
    \begin{tikzpicture}[baseline=(base.south)]
        \node (base) at (0,-3.5) {};
        
        \node (T0top) at (0,0) {$T$};
        \node (T1top) at (1,0) {$T$};
        \node [circlemor] (T) at (0,-6) {$t$};
        \node (S) [circlemor] at (0,-1.5) {$s$};
        \node (SS) [circlemor] at (1,-3) {$s$};
        \node (TT) [circlemor] at (1,-4.5) {$t$};

        \node (T0bot) at (0,-7) {$T$};
        \node (T1bot) at (1,-7) {$T$};
        
        
        \draw (T0top) [downward] to (S);
        \draw (T1top) [downward] to (SS);
        \draw (S) [downward] to (T);
        \draw (SS) [downward] to (TT);
        \draw (TT) [downward] to (T1bot);
        \draw (T) [downward] to (T0bot);
    \end{tikzpicture}
    =
    \begin{tikzpicture}[baseline=(base.south)]
        \node (base) at (0,-3.5) {};
        
        \node (T0top) at (0,0) {$T$};
        \node (T1top) at (1,0) {$T$};
        \node [circlemor] (T) at (0,-3) {$t$};
        \node (S) [circlemor] at (0,-1.5) {$s$};
        \node (SS) [circlemor] at (1,-4.5) {$s$};
        \node (TT) [circlemor] at (1,-6) {$t$};

        \node (T0bot) at (0,-7) {$T$};
        \node (T1bot) at (1,-7) {$T$};
        
        
        \draw (T0top) [downward] to (S);
        \draw (T1top) [downward] to (SS);
        \draw (S) [downward] to (T);
        \draw (SS) [downward] to (TT);
        \draw (TT) [downward] to (T1bot);
        \draw (T) [downward] to (T0bot);
    \end{tikzpicture}
    =
    \begin{tikzpicture}[baseline=(base.south)]
        \node (base) at (0,-3.5) {};
        
        \node (T0top) at (0,0) {$T$};
        \node (T1top) at (1,0) {$T$};
        \node (C0) at (.5,-2.5) {};
        \node (C1) at (.5,-4.5) {};

        \node (T0bot) at (0,-7) {$T$};
        \node (T1bot) at (1,-7) {$T$};
        
        
        \draw (T0top) [downcrossnw] to (C0.center);
        \draw (T1top) [downcrossne] to (C0.north east);
        \draw (C0.center) [se_ne] to (C1.north east);
        \draw (C0.south west) [sw_nw] to (C1.center);
        \draw (C1.center) [crosssedown] to (T1bot);
        \draw (C1.south west) [crossswdown] to (T0bot);
    \end{tikzpicture}
\end{equation*}
Here, in the first equation, we have isotoped so that in the second equation we can use the equation $\beta_{T,T} = \id_T \wedge ts$ from $(5)$. In the third equation, we have isotoped so that in the fourth equation we may use the equation $\beta_{T,T} = \id_T \wedge ts$ from $(5)$ again, as well as the equation $\beta_{T,T} = ts \wedge \id_T$ from $(3)$ which we have seen is equivalent to $(5)$. The resulting equation $\beta_{T,T}(s \wedge t) = \beta_{T,T}^2$ implies that $s \wedge t = \beta_{T,T}$ because $\beta_{T,T}$ is invertible. Thus $(5) \Rightarrow (1)$ as desired.

Of course, $(1) \Leftrightarrow (5)$ implies that $(2) \Leftrightarrow (4)$ by symmetry. Combining with the earlier equivalences, we are done.
\end{proof}

\begin{Def}\label{def:tt}
Let $T$ be an object in a braided monoidal category and let $t$ be an endomorphism of $T$. We say that $T$ has \defterm{$t$-twisted trivial braiding} if the braiding morphism $\beta_{T,T}: T \wedge T \to T \wedge T$ is equal to $\id_T \wedge t$ (or equivalently by \cref{prop:twisted-equiv}, to $t \wedge \id_T$). If $t$ is the identity, we say that $T$ has a \defterm{trivial braiding}. We say that $T$ has a \defterm{twisted-trivial braiding} to mean that it has a $t$-twisted trivial braiding for some $t$.
\end{Def} 

\begin{rmk}\label{rmk:twisted-triv-def}
Twisted trivial braidings are preserved by braided monoidal functors: if $X \in \calC$ has twisted-trivial braiding and $F: \calC \to \calD$ is a braided monoidal functor, then $F(X)$ also has twisted-trivial braiding. 
\end{rmk}

Objects with twisted-trivial braiding show up in a variety of contexts as ``pre-$\wedge$-invertible" objects.

\begin{eg}\label{eg:twist-sphere}
The sphere $S^1$ has $(-1)$-trivial braiding as an object of $\ho\Top_\ast$ pointed at $1 = e^0$. Here $(-1)$ is the automorphism of $S^1$ sending $e^{i\theta}$ to $e^{-i\theta}$. This can be seen because the matrix $\swapmatrix = \begin{pmatrix} 0 & 1 \\ 1 & 0 \end{pmatrix}$ is homotopic through invertible matrices to the matrix $1 \oplus (-1) = \begin{pmatrix} 1 & 0 \\ 0 & -1 \end{pmatrix}$, for instance by taking the homotopy $\begin{pmatrix} \sin t & \cos t \\ \cos t & -\sin t \end{pmatrix}$.
\end{eg}

\begin{eg}\label{eg:twist-susp-unit}
Let $(\calC,\wedge,S)$ be an $E_n$-monoidal pointed $\infty$-category ($n\geq 2$) with suspension (the suspension of $X$ is denoted $\Sigma X$) over which $\wedge$ distributes. Then $\Sigma S$ has $(-1)_\calC$-twisted trivial braiding, where $(-1)_\calC: \Sigma S \to \Sigma S$ is induced by the automorphism $(-1)$ of $S^1$ from \cref{eg:twist-sphere}.
\end{eg}

\begin{eg}\label{eg:twist-G-sphere}
Let $G$ be a compact Lie group, and $V$ a representation of dimension $n$. Let $S^V$ be the corresponding representation sphere, i.e. the one-point compactification of $V$. Then $S^V$ is naturally an object of the homotopy category $\ho G\Top_\ast$ of based $G$-spaces, with basepoint either at $0$ or $\infty$. Either basepoint is fixed by the $G$-equivariant automorphism $(-1)^n: S^V \to S^V$ induced by $-\id_V$. The braiding on $S^V$ is $(-1)^n$-trivial. This can be seen by tensoring the linear homotopy in \cref{eg:twist-sphere} with $V$ to get an equivariant homotopy through automorphisms of $V$, inducing an equivariant homotopy on spheres $S^V$. \cref{eg:twist-sphere} is the case where $G$ is the trivial group.

This example leads to the main result of this thesis, \cref{cor:main-result}.
\end{eg}

\begin{eg}\label{eg:twist-tate-sphere}
Let $S$ be a scheme. Then $\pjv^1_S$, the Tate sphere over $S$, admits a basepoint at either $0$ or $\infty$; either way it is preserved by the automorphism $(-1)$ induced by the $(-1)$ map on $\aff^1_S$, the trivial bundle of rank 1. \footnote{Note that even after passing to the stable category, this map is not the additive inverse of the identity, so $(-1)$ is perhaps a misnomer. For instance, its Betti realization is homotopic to the identity, not the additive inverse of the identity.} The braiding on $\pjv^1_S$ is $(-1)$-trivial in the homotopy category $\ho\Sm_{S,\ast}$ of pointed smooth schemes over $S$ localized at $\aff^1_S$. To see this, we use the Thom space structure: $\pjv^1_S \wedge_S \pjv^1_S = \Th(\aff^1_S) \wedge_S \Th(\aff^1_S) = \Th(\aff^1_S \times_S \aff^1_S) = \Th(\aff^2_S)$. Hence an $\aff^1$-homotopy through linear automorphisms of $\aff^2_S$ gives rise to an $\aff^1$-homotopy between the induced self-maps of $\pjv^1_S \wedge_S \pjv^1_S$, and so it suffices to show that the matrices $\swapmatrix$ and $1 \oplus (-1)$ from (\cref{eg:twist-sphere}) are $\aff^1$-homotopic in $\GL_2$. It's not clear that this can be done through a single $\aff^1$-homotopy, but it may be done in two steps, as the linear homotopies $\begin{pmatrix} t & 1-t \\ 1+t & -t\end{pmatrix}$ and $\begin{pmatrix}1 & 0 \\ 2t & -1 \end{pmatrix}$ (note that both have constant determinant $-1$) exhibit $\swapmatrix$ and $1 \oplus (-1)$ respectively as $\aff^1$-homotopic to $\begin{pmatrix}1 & 0 \\ 2 & -1 \end{pmatrix}$.

We plan to study this example further in future work.
\end{eg}

In the literature, a slightly weaker condition is often used instead:

\begin{Def}[Smith, \cite{voevodsky}]\label{def:smith}
Let $T$ be an object in a braided monoidal category. We say that $T$ is \defterm{symmetric} if the map $(\id_T \wedge \beta_{T,T}\inv)(\beta_{T,T} \wedge \id_T)$ is equal to $\id_{T^3}$.
\begin{equation*}
    \begin{tikzpicture}[baseline=(base.south)]
        \node (base) at (0,-1.75) {};
        
        \node (T0top) at (0,0) {$T$};
        \node (T1top) at (1,0) {$T$};
        \node (T2top) at (2,0) {$T$};
        \node (C0) at (.5,-1) {};
        \node (C1) at (1.5,-2.5) {};

        \node (T0bot) at (0,-3.5) {$T$};
        \node (T1bot) at (1,-3.5) {$T$};
        \node (T2bot) at (2,-3.5) {$T$};
        
        
        \draw (T0top) [downcrossnw] to (C0.center);
        \draw (T1top) [downcrossne] to (C0.north east);
        \draw (T2top) [downcrossne] to (C1.center);
        \draw (C0.center) [se_nw] to (C1.north west);
        \draw (C0.south west) [crossswdown] to (T0bot);
        \draw (C1.center) [crossswdown] to (T1bot);
        \draw (C1.south east) [crosssedown] to (T2bot);
    \end{tikzpicture}
    =
    \begin{tikzpicture}[baseline=(base.south)]
        \node (base) at (0,-1.75) {};
        
        \node (T0top) at (0,0) {$T$};
        \node (T1top) at (1,0) {$T$};
        \node (T2top) at (2,0) {$T$};
        \node (C0) at (.5,-1) {};
        \node (C1) at (1.5,-2.5) {};

        \node (T0bot) at (0,-3.5) {$T$};
        \node (T1bot) at (1,-3.5) {$T$};
        \node (T2bot) at (2,-3.5) {$T$};
        
        
        \draw (T0top) [downward] to (T0bot);
        \draw (T1top) [downward] to (T1bot);
        \draw (T2top) [downward] to (T2bot);
    \end{tikzpicture}
\end{equation*}    
\end{Def}

\begin{rmk}
We do not know whether \cref{def:smith} has appeared in the literature in the generality of a braided monoidal category (rather than just a symmetric monoidal category). The choice of handedness for the crossings is made to ensure that \cref{prop:tt-smith} below is true. In a symmetric monoidal category, if $T$ is symmetric then for any $n \in \nats$, the homomorphism $\Sigma_n \to \Aut(T^{\wedge n})$ factors through the sign homomorphism $\operatorname{sgn} : \Sigma_n \to C_2$. Similarly, in the braided setting, if $T$ is symmetric then for every $n \in \nats$, the canonical map $B_n \to \Aut(T^{\wedge n})$ factors through the abelianization homomorphism $B_n \to \ints$ (where $B_n$ is the braid group on $n$ strands). This is easily seen from the usual generators-and-relations description of $B_n$, where generator $b_i$ is the positive crossing of the $i$th strand over the $(i+1)$st strand, since the symmetry condition forces that $b_i b_{i+1}\inv \equiv 1$.
\end{rmk}

\begin{lem}\label{prop:tt-smith}
Let $T$ be an object in a braided monoidal category with twisted-trivial braiding. Then $T$ is symmetric.
\end{lem}
\begin{proof}
Let $t$ be the twist. We have (using \cref{prop:twisted-equiv}):
\begin{equation*}
    \begin{tikzpicture}[baseline=(base.south)]
        \node (base) at (0,-1.75) {};
        
        \node (T0top) at (0,0) {$T$};
        \node (T1top) at (1,0) {$T$};
        \node (T2top) at (2,0) {$T$};
        \node (C0) at (.5,-1) {};
        \node (C1) at (1.5,-2.5) {};

        \node (T0bot) at (0,-3.5) {$T$};
        \node (T1bot) at (1,-3.5) {$T$};
        \node (T2bot) at (2,-3.5) {$T$};
        
        
        \draw (T0top) [downcrossnw] to (C0.center);
        \draw (T1top) [downcrossne] to (C0.north east);
        \draw (T2top) [downcrossne] to (C1.center);
        \draw (C0.center) [se_nw] to (C1.north west);
        \draw (C0.south west) [crossswdown] to (T0bot);
        \draw (C1.center) [crossswdown] to (T1bot);
        \draw (C1.south east) [crosssedown] to (T2bot);
    \end{tikzpicture}
    =
    \begin{tikzpicture}[baseline=(base.south)]
        \node (base) at (0,-1.75) {};
        
        \node (T0top) at (0,0) {$T$};
        \node (T1top) at (1,0) {$T$};
        \node (T2top) at (2,0) {$T$};
        \node (T) [circlemor] at (1,-1) {$t$};
        \node (C1) at (1.5,-2.5) {};

        \node (T0bot) at (0,-3.5) {$T$};
        \node (T1bot) at (1,-3.5) {$T$};
        \node (T2bot) at (2,-3.5) {$T$};
        
        
        \draw (T0top) [downward] to (T0bot);
        \draw (T1top) [downward] to (T.north);
        \draw (T2top) [downcrossne] to (C1.center);
        \draw (T.south) [downcrossnw] to (C1.north west);
        \draw (C1.center) [crossswdown] to (T1bot);
        \draw (C1.south east) [crosssedown] to (T2bot);
    \end{tikzpicture}
    =
    \begin{tikzpicture}[baseline=(base.south)]
        \node (base) at (0,-1.75) {};
        
        \node (T0top) at (0,0) {$T$};
        \node (T1top) at (1,0) {$T$};
        \node (T2top) at (2,0) {$T$};
        \node (C0) at (1.5,-1) {};
        \node (C1) at (1.5,-2.5) {};

        \node (T0bot) at (0,-3.5) {$T$};
        \node (T1bot) at (1,-3.5) {$T$};
        \node (T2bot) at (2,-3.5) {$T$};
        
        
        \draw (T0top) [downward] to (T0bot);
        \draw (T1top) [downcrossnw] to (C0.center);
        \draw (T2top) [downcrossne] to (C0.north east);
        \draw (C0.center) [se_ne] to (C1.center);
        \draw (C0.south west) [sw_nw] to (C1.north west);
        \draw (C1.center) [crossswdown] to (T1bot);
        \draw (C1.south east) [crosssedown] to (T2bot);
    \end{tikzpicture}
    =
    \begin{tikzpicture}[baseline=(base.south)]
        \node (base) at (0,-1.75) {};
        
        \node (T0top) at (0,0) {$T$};
        \node (T1top) at (1,0) {$T$};
        \node (T2top) at (2,0) {$T$};
        \node (C0) at (.5,-1) {};
        \node (C1) at (1.5,-2.5) {};

        \node (T0bot) at (0,-3.5) {$T$};
        \node (T1bot) at (1,-3.5) {$T$};
        \node (T2bot) at (2,-3.5) {$T$};
        
        
        \draw (T0top) [downward] to (T0bot);
        \draw (T1top) [downward] to (T1bot);
        \draw (T2top) [downward] to (T2bot);
    \end{tikzpicture}
\end{equation*}
\end{proof}

The converse holds if $T$ is dualizable -- a condition we are interested in forcing anyway.

\begin{lem}\label{prop:smith-tt}
Let $T$ be an object in a braided monoidal category. Assume that $T$ is symmetric and dualizable. Then $T$ has $t$-twisted trivial braiding, where $t = \id_T \wedge (\varepsilon \beta_{T,T^\vee}\inv \eta)$ is mulitplication by the Euler characteristic of $T$.
\end{lem}
\begin{proof}
First, we introduce some cups, caps and crossings:
\begin{equation*}
    \begin{tikzpicture}[baseline=(base.south)]
        \node (base) at (0,-1.75) {};
        
        \node (T0top) at (0,0) {$T$};
        \node (T1top) at (1,0) {$T$};
        \node (C0) at (.5,-1) {};

        \node (T0bot) at (0,-3.5) {$T$};
        \node (T1bot) at (1,-3.5) {$T$};
        
        
        \draw (T0top) [downcrossnw] to (C0.center);
        \draw (T1top) [downcrossne] to (C0.north east);
        \draw (C0.center) [crosssedown, witharrow] to (T1bot);
        \draw (C0.south west) [crossswdown, witharrow] to (T0bot);
    \end{tikzpicture}
    = \begin{tikzpicture}[baseline=(base.south)]
        \node (base) at (0,-1.75) {};
        
        \node (T0top) at (0,0) {$T$};
        \node (T1top) at (1,0) {$T$};
        \node (Cap) [mycap] at (2.5,0) {};
        \node (C0) at (.5,-1) {};
        \node (C1) at (2.5,-1) {};
        \node (C2) at (2.5,-2) {};
        \node (Cup) [mycup] at (1.5,-3) {};

        \node (T0bot) at (0,-3.5) {$T$};
        \node (T1bot) at (3,-3.5) {$T$};

        
        \draw (T0top) [downcrossnw] to (C0.center);
        \draw (T1top) [downcrossne] to (C0.north east);
        \draw (Cap.center) [w_nw] to (C1.north west);
        \draw (Cap.center) [e_ne] to (C1.center);
        \draw (C0.center) [se_w, witharrow] to (Cup.center);
        \draw (C0.south west) [crossswdown, witharrow] to (T0bot);
        \draw (C1.center) [sw_nw, witharrow] to (C2.center);
        \draw (C1.south east) [se_ne,withbackarrow] to (C2.north east);
        \draw (C2.south west) [sw_e] to (Cup.center);
        \draw (C2.center) [crosssedown] to (T1bot);
    \end{tikzpicture}
\end{equation*}
Now we introduce a few more crossings, then apply the symmetry condition, and finally simplify:
\begin{equation*}
    = 
    \begin{tikzpicture}[baseline=(base.south)]
        \node (base) at (0,-1.75) {};
        
        \node (T0top) at (0,0) {$T$};
        \node (T1top) at (1,0) {$T$};
        \node (Cap) [mycap] at (2.5,0) {};
        \node (C0) at (.5,-1) {};
        \node (C1) at (2.5,-1) {};
        \node (C2) at (2.5,-4) {};
        \node (Cup) [mycup] at (1.5,-5) {};
        \node (C3) at (1.5,-2) {};
        \node (C4) at (1.5,-3) {};

        \node (T0bot) at (0,-5.5) {$T$};
        \node (T1bot) at (3,-5.5) {$T$};

        
        \draw (T0top) [downcrossnw] to (C0.center);
        \draw (T1top) [downcrossne] to (C0.north east);
        \draw (Cap.center) [w_nw] to (C1.north west);
        \draw (Cap.center) [e_ne] to (C1.center);
        \draw (C0.center) [se_nw, witharrow] to (C3.north west);
        \draw (C0.south west) [crossswdown, witharrow] to (T0bot);
        \draw (C1.center) [sw_ne, witharrow] to (C3.center);
        \draw (C3.center) [sw_nw] to (C4.center);
        \draw (C3.south east) [se_ne] to (C4.north east);
        \draw (C4.south west) [sw_w] to (Cup.center);
        \draw (C4.center) [se_nw] to (C2.center);
        \draw (C1.south east) [se_ne, withbackarrow] to (C2.north east);
        \draw (C2.south west) [sw_e] to (Cup.center);
        \draw (C2.center) [crosssedown] to (T1bot);
    \end{tikzpicture}
    = 
    \begin{tikzpicture}[baseline=(base.south)]
        \node (base) at (0,-1.75) {};
        
        \node (T0top) at (0,0) {$T$};
        \node (T1top) at (1,0) {$T$};
        \node (Cap) [mycap] at (2.5,0) {$T$};
        \node (C1) at (2.5,-1) {};
        \node (C2) at (2.5,-4) {};
        \node (Cup) [mycup] at (1.5,-5) {};
        \node (C4) at (1.5,-3) {};

        \node (T0bot) at (0,-5.5) {$T$};
        \node (T1bot) at (3,-5.5) {$T$};

        
        \draw (T0top) [downward, witharrow] to (T0bot);
        \draw (T1top) [downcrossnw, witharrow] to (C4.center);
        \draw (Cap.center) [w_nw] to (C1.north west);
        \draw (Cap.center) [e_ne] to (C1.center);
        \draw (C1.center) [sw_ne, witharrow] to (C4.north east);
        \draw (C4.south west) [sw_w] to (Cup.center);
        \draw (C4.center) [se_nw] to (C2.center);
        \draw (C1.south east) [se_ne, withbackarrow] to (C2.north east);
        \draw (C2.south west) [sw_e] to (Cup.center);
        \draw (C2.center) [crosssedown] to (T1bot);
    \end{tikzpicture}
    = 
    \begin{tikzpicture}[baseline=(base.south)]
        \node (base) at (0,-1.75) {};
        
        \node (T0top) at (0,0) {$T$};
        \node (T1top) at (1,0) {$T$};
        \node (Cap) [mycap] at (2.5,0) {$T$};
        \node (C1) at (2.5,-1) {};
        \node (Cup) [mycup] at (2.5,-2) {};

        \node (T0bot) at (0,-5.5) {$T$};
        \node (T1bot) at (3,-5.5) {$T$};

        
        \draw (T0top) [downward, witharrow] to (T0bot);
        \draw (T1top) [downward, witharrow] to (T1bot);
        \draw (Cap.center) [w_nw] to (C1.north west);
        \draw (Cap.center) [e_ne] to (C1.center);
        \draw (C1.center) [sw_w] to (Cup.center);
        \draw (C1.south east) [se_e, withbackarrow] to (Cup.center);
    \end{tikzpicture}
\end{equation*}
\end{proof}

\subsection{Idempotent objects and smashing localizations}\label{sec:idem-smash}

In this section we review some basic concepts of smashing and co-smashing localizations.
Much of the material of this section may be found in the preprint \cite{boyarchenko-drinfeld}, which was previously available on Drinfeld's website.

\begin{Def}[\cite{boyarchenko-drinfeld}]
Let $(\calC,S,\wedge)$ be a braided monoidal category.
\begin{itemize}
    \item A \defterm{closed idempotent} in $\calC$ is an object $E$ equipped with a morphism $r: S \to E$ such that the induced morphism $r \wedge \id_E: S \wedge E \to E \wedge E$ is an isomorphism (equivalently, $\id_E \wedge r: E \wedge S \to E \wedge E$ is an isomorphism).
    \item Dually, an \defterm{open idempotent} is a closed idempotent in $\calC^\op$, i.e. an object $E$ equipped with a map $i: E \to S$ such that $i \wedge \id_E$ is an isomorphism.
    \item A \defterm{clopen idempotent} is an object $E$ with maps $r: S \to E$, $i: E \to S$ satisfying the equations \[ri = \id_E \quad \text{(``splitting")} \quad \text{and} \quad ir \wedge \id_E = \id_{E \wedge E} \quad \text{(``stability")}\]
    \item If $r: S \to E$ (resp. $i: E \to S$) is a map and $X \in \calC$, say that $X$ is \defterm{$r$-stable} (resp. \defterm{$i$-stable}) if $r \wedge \id_X: X \to E \wedge X$ (resp. $i \wedge \id_X: E \wedge X \to X$) is an isomorphism. We may say \defterm{$E$-stable} if $r$ (resp. $i$) is understood. We write $\calC_E \subseteq \calC$ for the full subcategory of $E$-stable objects.
    \item If $E$ is an object and $\calC$ is pointed, say that an object $X$ is \defterm{$E$-torsion} if $E \wedge X = 0$.
\end{itemize}
\end{Def}

\begin{rmk}
Any braided monoidal functor $F: \calC \to \calD$ preserves open, closed, and clopen idempotents, along with stability in the sense that if $X$ is $E$-stable, then $F(X)$ is $F(E)$-stable. If moreover $\calC$ and $\calD$ are pointed and $F$ preserves the zero object, then if $X$ is $E$-torsion, then $F(X)$ is $F(E)$-torsion.
\end{rmk}

\begin{prop}\label{prop:smashing-localization}
Let $(\calC, \wedge, S)$ be a braided monoidal category and $(E, r: S \to E)$ a closed idempotent. Then the functor $E \wedge (-): \calC \to \calC$ is a localization, with essential image $\calC_E$. The unit at $X \in \calC$ is given by $r \wedge \id_X: X \to E \wedge X$, and the equation $r \wedge \id_E = \id_E \wedge r$ holds.

Moreover, $\calC_E$ is braided monoidal under $\wedge$ with unit $E$, and the localization $\calC \to \calC_E$ is braided monoidal.
\end{prop}
\begin{proof}
As a first step, we check that $E \wedge (-)$ is faithful when restricted to its image. Let $f,g: E \wedge X {}^\to_\to E \wedge Y$, and suppose that $\id_E \wedge f = \id_E \wedge g$. To this equation postcompose $(r \wedge \id_E)\inv \wedge \id_Y$ and precompose $r \wedge \id_E \wedge \id_X$; slide the $r$ down and cancel the inverses to obtain $f = g$.

\begin{equation}
    \begin{tikzpicture}[baseline=(base.south)]
        \node (base) at (0,2) {};

        \node (E0top) at (0,4) {$E$};
        \node (Etop) at (1,4) {$E$};
        \node (Xtop) at (2,4) {$X$};
        
        \node (E0bot) at (0,0) {$E$};
        \node (Ebot) at (1,0) {$E$};
        \node (Xbot) at (2,0) {$X$};
        
        \node [shape=rectangle, draw] (F) at (1.5, 2) {$f$};

        \draw (E0top) to [downward] (E0bot);
        \draw (Etop) to [downward] (F.100);
        \draw (Xtop) to [downward] (F.80);
        \draw (Ebot) to [upward] (F.260);
        \draw (Xbot) to [upward] (F.280);
    \end{tikzpicture}
    =
    \begin{tikzpicture}[baseline=(base.south)]
        \node (base) at (0,2) {};

        \node (E0top) at (0,4) {$E$};
        \node (Etop) at (1,4) {$E$};
        \node (Xtop) at (2,4) {$X$};
        
        \node (E0bot) at (0,0) {$E$};
        \node (Ebot) at (1,0) {$E$};
        \node (Xbot) at (2,0) {$X$};
        
        \node [shape=rectangle, draw] (G) at (1.5, 2) {$g$};

        \draw (E0top) to [downward] (E0bot);
        \draw (Etop) to [downward] (G.100);
        \draw (Xtop) to [downward] (G.80);
        \draw (Ebot) to [upward] (G.260);
        \draw (Xbot) to [upward] (G.280);
    \end{tikzpicture}
\end{equation}

\begin{equation}
    \begin{tikzpicture}[baseline=(base.south)]
        \node (base) at (0,1.5) {};
        
        \node (Etop) at (1,4) {$E$};
        \node (Xtop) at (2,4) {$X$};
        
        \node (Ebot) at (1,-1) {$E$};
        \node (Xbot) at (2,-1) {$X$};
        
        \node [shape=circle, draw] (R) at (0,3) {$r$};
        \node [shape=rectangle, draw] (F) at (1.5, 2) {$f$};
        \node [shape=rectangle, draw] (RI) at (0.25,0) {$(r \wedge \id_E)\inv$};
        
        \draw (R) to [downward] node [left] {$E$} (RI.115);
        \draw (Etop) to [downward] (F.100);
        \draw (Xtop) to [downward] (F.80);
        \draw (RI.65) to [upward] node [right] {$E$} node [right] {$E$} (F.260);
        \draw (Xbot) to [upward] (F.280);
        \draw (Ebot) to [upward] (RI.-90);
    \end{tikzpicture}
    =
    \begin{tikzpicture}[baseline=(base.south)]
        \node (base) at (0,1.5) {};
            
        \node (Etop) at (1,4) {$E$};
        \node (Xtop) at (2,4) {$X$};
        
        \node (Ebot) at (1,-1) {$E$};
        \node (Xbot) at (2,-1) {$X$};
        
        \node [shape=circle, draw] (R) at (0,3) {$r$};
        \node [shape=rectangle, draw] (G) at (1.5, 2) {$g$};
        \node [shape=rectangle, draw] (RI) at (0.25,0) {$(r \wedge \id_E)\inv$};
        
        \draw (R) to [downward] node [left] {$E$} (RI.115);
        \draw (Etop) to [downward] (G.100);
        \draw (Xtop) to [downward] (G.80);
        \draw (RI.65) to [upward] node [right] {$E$} (G.260);
        \draw (Xbot) to [upward] (G.280);
        \draw (Ebot) to [upward] (RI.-90);
    \end{tikzpicture}
\end{equation}

\begin{equation}
    \begin{tikzpicture}[baseline=(base.south)]
        \node (base) at (0,1.5) {};
            
        \node (Etop) at (1,4) {$E$};
        \node (Xtop) at (2,4) {$X$};
        
        \node (Ebot) at (1,-1) {$E$};
        \node (Xbot) at (2,-1) {$X$};
        
        \node [shape=circle, draw] (R) at (0,1.5) {$r$};
        \node [shape=rectangle, draw] (F) at (1.5, 2.5) {$f$};
        \node [shape=rectangle, draw] (RI) at (0.25,0) {$(r \wedge \id_E)\inv$};
        
        \draw (R) to [downward] node [left] {$E$} (RI.115);
        \draw (Etop) to [downward] (F.100);
        \draw (Xtop) to [downward] (F.80);
        \draw (RI.65) to [upward] node [right] {$E$} (F.260);
        \draw (Xbot) to [upward] (F.280);
        \draw (Ebot) to [upward] (RI.-90);
    \end{tikzpicture}
    =
    \begin{tikzpicture}[baseline=(base.south)]
        \node (base) at (0,1.5) {};
            
        \node (Etop) at (1,4) {$E$};
        \node (Xtop) at (2,4) {$X$};
        
        \node (Ebot) at (1,-1) {$E$};
        \node (Xbot) at (2,-1) {$X$};
        
        \node [shape=circle, draw] (R) at (0,1.5) {$r$};
        \node [shape=rectangle, draw] (G) at (1.5, 2.5) {$g$};
        \node [shape=rectangle, draw] (RI) at (0.25,0) {$(r \wedge \id_E)\inv$};
        
        \draw (R) to [downward] node [left] {$E$} (RI.115);
        \draw (Etop) to [downward] (G.100);
        \draw (Xtop) to [downward] (G.80);
        \draw (RI.65) to [upward] node [right] {$E$} (G.260);
        \draw (Xbot) to [upward] (G.280);
        \draw (Ebot) to [upward] (RI.-90);
    \end{tikzpicture}
\end{equation}

\begin{equation}
    \begin{tikzpicture}[baseline=(base.south)]
        \node (base) at (0,2) {};
            
        \node (Etop) at (1,4) {$E$};
        \node (Xtop) at (2,4) {$X$};
        
        \node (Ebot) at (1,0) {$E$};
        \node (Xbot) at (2,0) {$X$};
        
        \node [shape=rectangle, draw] (F) at (1.5, 2) {$f$};

        \draw (Etop) to [downward] (F.100);
        \draw (Xtop) to [downward] (F.80);
        \draw (Ebot) to [upward] (F.260);
        \draw (Xbot) to [upward] (F.280);
    \end{tikzpicture}
    =
    \begin{tikzpicture}[baseline=(base.south)]
        \node (base) at (1,2) {};
            
        \node (Etop) at (1,4) {$E$};
        \node (Xtop) at (2,4) {$X$};
        
        \node (Ebot) at (1,0) {$E$};
        \node (Xbot) at (2,0) {$X$};
        
        \node [shape=rectangle, draw] (G) at (1.5, 2) {$g$};

        \draw (Etop) to [downward] (G.100);
        \draw (Xtop) to [downward] (G.80);
        \draw (Ebot) to [upward] (G.260);
        \draw (Xbot) to [upward] (G.280);
    \end{tikzpicture}
\end{equation}

Now we claim that $r \wedge \id_E = \id_E \wedge r$. By the first step, it suffices to verify that $\id_E \wedge r \wedge \id_E = \id_E \wedge \id_E \wedge r$. Precomposing the isomorphism $\id_E \wedge r$, it suffices to verify that $\id_E \wedge r \wedge r = \id_E \wedge r \wedge r$, which is true.

We now have an endofunctor $E \wedge(-): \calC \to \calC$ and a natural transformation $r \wedge (-): \id_\calC \Rightarrow E \wedge(-)$ such that $r\wedge E \wedge (-) = E \wedge r \wedge (-)$ is invertible. This is precisely the data of a localization. It is then a general fact that the essential image of $E \wedge (-)$ consists precisely of those objects $Y$ with $r \wedge Y: Y \to E \wedge Y$ an isomorphism, i.e. the $E$-stable objects.

The associator and braiding for $\calC_E$ are as in $\calC$. For the unit, note that if $X \in \calC_E$, then $r \wedge \id_X: X \to E \wedge X$ is invertible; the left unit is defined to be the inverse to this morphism, and similarly the right unit. Only the unit equation must be checked, and it is immediate upon precomposition with $r$. To make the localization functor monoidal, use copies of $r$ for the constraints. It is straightforward to check that this makes $E \wedge (-)$ into a monoidal functor; the fact that it is braided follows from the equation $r \wedge \id_E = \id_E \wedge r$.
%
\end{proof}


\subsection{Twisted trivial braiding, duality, and idempotents}\label{sec:tt-dual}

In this section, we show that any dualizable object with twisted-trivial braiding gives rise to a smashing-cosmashing localization \cref{prop:clopens}. As we shall review in \cref{sec:comp-smash}, this often implies that the entire category splits along lines defined by this object.

\begin{prop}\label{prop:clopens}
Let $(\calC, \wedge, S)$ be a braided monoidal category and $E \in \calC$ an object.
\begin{enumerate}
    \item\label{prop:clopens.item:selfdual-trivial} If $E$ is a closed or open idempotent, then $E$ has trivial braiding.
    \item\label{prop:clopens.item:equiv} The following are equivalent:
    \begin{enumerate}
        \item $E$ admits the structure of a clopen idempotent.
        \item $E$ admits the structure of a closed idempotent and the structure of an open idempotent.
        \item $E$ is dualizable and admits the structure of a closed idempotent.
    \end{enumerate}
    Moreover, in this case $E$ is self-dual. For any closed structure on $E$, there is at most one clopen structure compatible with it.
    \item\label{prop:clopens.item:untwist} If $T \in \calC$ has twisted-trivial braiding and is dualizable, then $T^\dual \wedge T$ is a clopen idempotent (and in particular has trivial braiding by (\ref{prop:clopens.item:selfdual-trivial})).
    \item\label{prop:clopens.item:stable-t} If $E$ is a closed idempotent, then an object $X$ is $E$-stable iff $X$ is of the form $X = E \wedge Y$. If $T \in \calC$ has twisted-trivial braiding, then $X$ is $T^\dual \wedge T$-stable iff $X$ is of the form $X = T \wedge Z$.
\end{enumerate}
\end{prop}

\begin{proof}
(1): Let $r: S \to E$ be a closed idempotent. We want to show that $\beta_{E,E} = \id_{E\wedge E}$. Composing with the isomorphism $r \wedge \id_E$, it suffices to show that $r \wedge \id_E = \id_E \wedge r$. This is part of Proposition \ref{prop:smashing-localization}. The open case is dual.
    
(2): $(a) \Rightarrow (b)$: (This implication is superfluous in light of the others we will prove, but we point it out because it is the most straightforward.) If $(E,r,i)$ is a clopen idempotent, then the splitting and stability equations yield explicit inverses for $r \wedge \id_E$ and $i \wedge \id_E$, so that $(E,r)$ is a closed idempotent and $(E,i)$ is an open idempotent.
    
$(b) \Rightarrow (a)$: Let $(E,r)$ be a closed idempotent and $(E,i)$ be an open idempotent. Set 
$i' = i (\id_E \wedge r)\inv (\id_E \wedge i)\inv$. Then we show that $(E,r,i')$ is a clopen idempotent, as follows. First we verify the stability equation $i'r \wedge \id_E = \id_E$. 

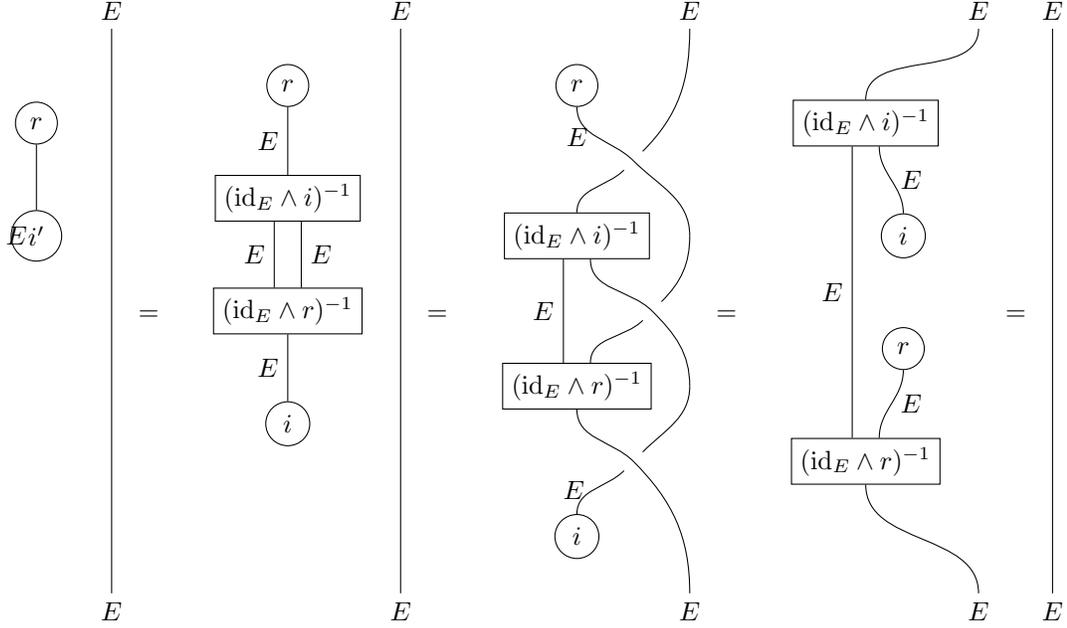
\begin{equation*}
    \begin{tikzpicture}[baseline=(base.south)]
        \node (base) at (0,-4) {};
            
        \node (Etop) at (1,0) {$E$};

        \node (Ebot) at (1,-8) {$E$};

        \node [shape=circle, draw] (R) at (0, -1.5) {$r$};
        \node [shape=circle, draw] (Ip) at (0,-3) {$i'$};

        \draw (Etop) to [downward] (Ebot);
        \draw (R) to [downward] (Ip) node [left] {$E$};
    \end{tikzpicture}
=
    \begin{tikzpicture}[baseline=(base.south)]
        \node (base) at (0,-4) {};
            
        \node (Etop) at (3,0) {$E$};

        \node (Ebot) at (3,-8) {$E$};

        \node [circlemor] (R) at (1.5, -1) {$r$};
        \node [rectmor] (II) at (1.5,-2.5) {$(\id_E \wedge i)\inv$};
        \node [rectmor] (RI) at (1.5,-4) {$(\id_E \wedge r)\inv$};
        \node [circlemor] (I) at (1.5,-5.5) {$i$};

        \draw (Etop) to [downward] (Ebot);
        \draw (R) to [downward] node [left] {$E$} (II);
        \draw (II.-120) to [downward]  node [left] {$E$} (RI.120);
        \draw (II.-60) to [downward] node [right] {$E$} (RI.60);
        \draw (RI) to [downward] node [left] {$E$} (I);
    \end{tikzpicture}
=
    \begin{tikzpicture}[baseline=(base.south)]
        \node (base) at (0,-4) {};
            
        \node (Etop) at (3,0) {$E$};

        \node (Ebot) at (3,-8) {$E$};

        \node [circlemor] (R) at (1.5, -1) {$r$};
        \node [crossing] (C1) at (2.25,-2) {};
        \node [rectmor] (II) at (1.5,-3) {$(\id_E \wedge i)\inv$};
        \node [crossing] (E1) at (3,-3) {};
        \node [crossing] (C2) at (2.5,-4) {};
        \node [rectmor] (RI) at (1.5,-5) {$(\id_E \wedge r)\inv$};
        \node [crossing] (E2) at (3,-5) {};
        \node [crossing] (C3) at (2.25,-6) {};
        \node [circlemor] (I) at (1.5,-7) {$i$};

        \draw (Etop) to [downcrossne] (C1.north east);
        \draw (R) to [downcrossnw] node [left] {$E$} (C1.center);
        \draw (C1.south west) to [crossswdown] (II);
        \draw (C1.center) to [crosssedown] (E1.center);
        \draw (II.-120) to [downward]  node [left] {$E$} (RI.120);
        \draw (II.-60) to [downcrossnw] (C2.center);
        \draw (E1.center) to [downcrossne] (C2.north east);
        \draw (C2.center) to [crosssedown] (E2.center);
        \draw (C2.south west) to [crossswdown] (RI.60);
        \draw (RI) to [downcrossnw] (C3.center);
        \draw (E2.center) to [downcrossne] (C3.north east);
        \draw (C3.center) to [crosssedown] (Ebot);
        \draw (C3.south west) to [crossswdown] node [left] {$E$} (I);
    \end{tikzpicture}
=
    \begin{tikzpicture}[baseline=(base.south)]
        \node (base) at (0,-4) {};
            
        \node (Etop) at (3,0) {$E$};

        \node (Ebot) at (3,-8) {$E$};

        \node [rectmor] (II) at (1.5,-1.5) {$(\id_E \wedge i)\inv$};
        \node [circlemor] (I) at (2, -3) {$i$};
        \node [circlemor] (R) at (2,-4.5) {$r$};
        \node [rectmor] (RI) at (1.5,-6) {$(\id_E \wedge r)\inv$};

        \draw (Etop) to [downward] (II);
        \draw (II.-120) to [downward]  node [left] {$E$} (RI.120);
        \draw (II.-60) to [downward] node [right] {$E$} (I);
        \draw (R) to [downward] node [right] {$E$} (RI.60);
        \draw (RI) to [downward] (Ebot);
    \end{tikzpicture}
    =
    \begin{tikzpicture}[baseline=(base.south)]
        \node (base) at (0,-4) {};
        \node (Etop) at (0,0) {$E$};

        \node (Ebot) at (0,-8) {$E$};
        
        \draw (Etop) to [downward] (Ebot);
    \end{tikzpicture}
\end{equation*}
Here, in the first equality we expand the definition of $i'$. In the second equality, we use the fact we showed in (1): $E$ has trivial braiding. In fact, we use this three times, introducing three crossings into the diagram. In the third equality we ``yank" the free ends of this string diagram through. In the fourth equality we cancel the two pairs of inverse morphisms, leaving us with $\id_E$, as desired.

Now we verify the splitting equation $ri' = \id_E$.
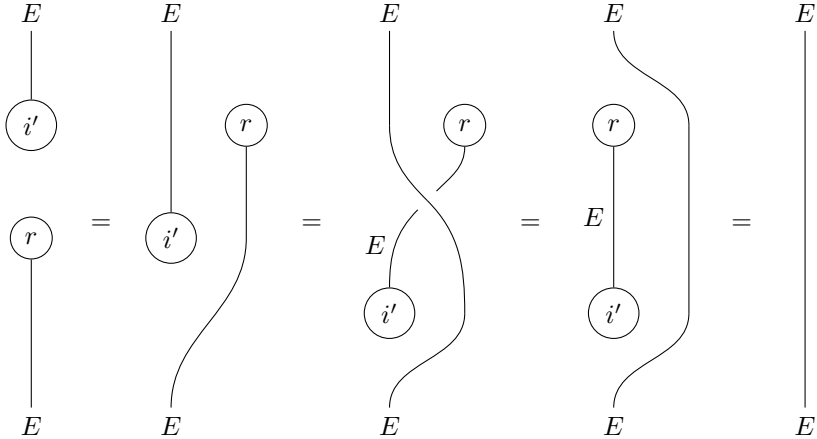
\begin{equation*}
    \begin{tikzpicture}[baseline=(base.south)]
        \node (base) at (0,-2.75) {};
        
        \node (Etop) at (0,0) {$E$};

        \node (Ebot) at (0,-5.5) {$E$};
        
        \node (Ip) [circlemor] at (0,-1.5) {$i'$};
        \node (R) [circlemor] at (0,-3) {$r$};
        
        \draw (Etop) to [downward] (Ip);
        \draw (R) to [downward] (Ebot);
    \end{tikzpicture}
    \quad = \quad
    \begin{tikzpicture}[baseline=(base.south)]
        \node (base) at (0,-2.75) {};
        
        \node (Etop) at (0,0) {$E$};

        \node (Ebot) at (0,-5.5) {$E$};
        
        \node (R) [circlemor] at (1,-1.5) {$r$};
        \node (Ip) [circlemor] at (0,-3) {$i'$};
        \node (E1) [crossing] at (1,-3) {};
        
        \draw (Etop) to [downward] (Ip);
        \draw (R) to [downward] (E1.center);
        \draw (E1.center) to [downward] (Ebot);
    \end{tikzpicture}
    \quad = \quad
    \begin{tikzpicture}[baseline=(base.south)]
        \node (base) at (0,-2.75) {};
        
        \node (Etop) at (0,0) {$E$};

        \node (Ebot) at (0,-5.5) {$E$};
        
        \node (R) [circlemor] at (1,-1.5) {$r$};
        \node (E2) [crossing] at (0,-1.5) {};
        \node (C) [crossing] at (0.5,-2.5) {};
        \node (Ip) [circlemor] at (0,-4) {$i'$};
        \node (E1) [crossing] at (1,-4) {};
        
        \draw (Etop) to [downward] (E2.center);
        \draw (E2.center) to [downcrossnw] (C.center);
        \draw (R) to [downcrossne] (C.north east);
        \draw (C.south west) [crossswdown] to node [left] {$E$} (Ip);
        \draw (C.center) to [crosssedown] (E1.center);
        \draw (E1.center) to [downward] (Ebot);
    \end{tikzpicture}
   \quad = \quad
    \begin{tikzpicture}[baseline=(base.south)]
        \node (base) at (0,-2.75) {};
        
        \node (Etop) at (0,0) {$E$};

        \node (Ebot) at (0,-5.5) {$E$};
        
        \node (E2) [crossing] at (1,-1.5) {};
        \node (R) [circlemor] at (0,-1.5) {$r$};
        \node (Ip) [circlemor] at (0,-4) {$i'$};
        \node (E1) [crossing] at (1,-4) {};
        
        \draw (Etop) to [downward] (E2.center);
        \draw (R) [downward] to node [left] {$E$} (Ip);
        \draw (E2.center) to [downward] (E1.center);
        \draw (E1.center) to [downward] (Ebot);
    \end{tikzpicture}
    \quad = \quad
    \begin{tikzpicture}[baseline=(base.south)]
        \node (base) at (0,-2.75) {};
        
        \node (Etop) at (0,0) {$E$};

        \node (Ebot) at (0,-5.5) {$E$};
        
        
        \draw (Etop) to [downward] (Ebot);
    \end{tikzpicture}
\end{equation*}

Here in the first equation we have simply slid strands past each other. In the second equation, we have used the triviality of $\beta_{E,E}$ from (1) to introduce a crossing. In the third equation, we have slid through the crossing, and the fourth equation then follows from the stability equation which we have already proven, leaving us with $\id_E$ as desired. Note that this argument did not require us to unpack the definition of $i'$ -- so we have in fact proven that the splitting equation follows from the stability equation plus the triviality of the braiding.

$(a) \Rightarrow (c)$: As remarked, $E$ is self-dual, with unit $\eta = r \wedge r$ and counit $\varepsilon = i \wedge i$. We verify one triangle equation:
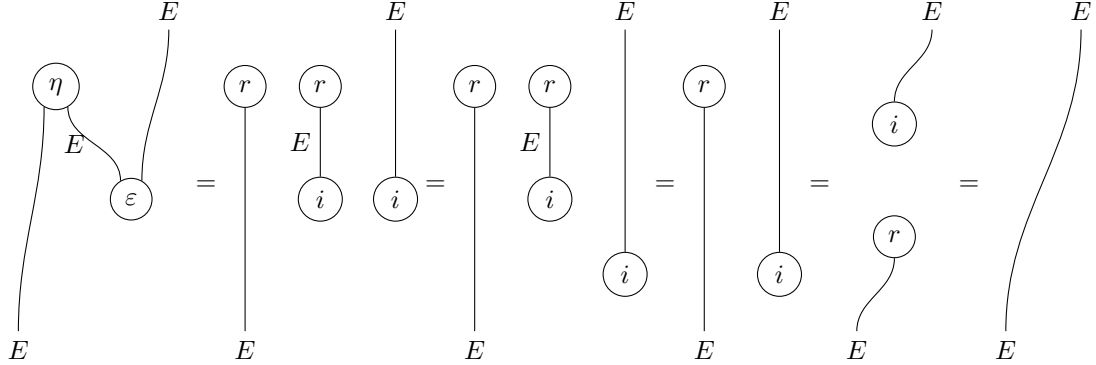
\begin{equation*}
    \begin{tikzpicture}[baseline=(base.south)]
        \node (base) at (0,-2.25) {};
        
        \node (Etop) at (2,0) {$E$};

        \node (Ebot) at (0,-4.5) {$E$};
        
        \node (eta) [circlemor] at (0.5,-1) {$\eta$};
        \node (eps) [circlemor] at (1.5,-2.5) {$\varepsilon$};
        
        \draw (Etop) to [downward] (eps.60);
        \draw (eta.-60) [downward] to node [left] {$E$} (eps.120);
        \draw (eta.-120) to [downward] (Ebot);
    \end{tikzpicture}
    =
    \begin{tikzpicture}[baseline=(base.south)]
        \node (base) at (0,-2.25) {};
        
        \node (Etop) at (2,0) {$E$};

        \node (Ebot) at (0,-4.5) {$E$};
        
        \node (R1) [circlemor] at (0,-1) {$r$};
        \node (R2) [circlemor] at (1,-1) {$r$};
        \node (I1) [circlemor] at (1,-2.5) {$i$};
        \node (I2) [circlemor] at (2,-2.5) {$i$};
        
        \draw (Etop) to [downward] (I2);
        \draw (R2) [downward] to node [left] {$E$} (I1);
        \draw (R1) to [downward] (Ebot);
    \end{tikzpicture}
    =
    \begin{tikzpicture}[baseline=(base.south)]
        \node (base) at (0,-2.25) {};
        
        \node (Etop) at (2,0) {$E$};

        \node (Ebot) at (0,-4.5) {$E$};
        
        \node (R1) [circlemor] at (0,-1) {$r$};
        \node (R2) [circlemor] at (1,-1) {$r$};
        \node (I1) [circlemor] at (1,-2.5) {$i$};
        \node (I2) [circlemor] at (2,-3.5) {$i$};
        
        \draw (Etop) to [downward] (I2);
        \draw (R2) [downward] to node [left] {$E$} (I1);
        \draw (R1) to [downward] (Ebot);
    \end{tikzpicture}
    =
    \begin{tikzpicture}[baseline=(base.south)]
        \node (base) at (0,-2.25) {};
        
        \node (Etop) at (1,0) {$E$};

        \node (Ebot) at (0,-4.5) {$E$};
        
        \node (R1) [circlemor] at (0,-1) {$r$};
        \node (I2) [circlemor] at (1,-3.5) {$i$};
        
        \draw (Etop) to [downward] (I2);
        \draw (R1) to [downward] (Ebot);
    \end{tikzpicture}
    =
    \begin{tikzpicture}[baseline=(base.south)]
        \node (base) at (0,-2.25) {};
        
        \node (Etop) at (1,0) {$E$};

        \node (Ebot) at (0,-4.5) {$E$};
        
        \node (I2) [circlemor] at (0.5,-1.5) {$i$};
        \node (R1) [circlemor] at (0.5,-3) {$r$};
        
        \draw (Etop) to [downward] (I2);
        \draw (R1) to [downward] (Ebot);
    \end{tikzpicture}
    =
    \begin{tikzpicture}[baseline=(base.south)]
        \node (base) at (0,-2.25) {};
        
        \node (Etop) at (1,0) {$E$};

        \node (Ebot) at (0,-4.5) {$E$};
        
        
        \draw (Etop) to [downward] (Ebot);
    \end{tikzpicture}
\end{equation*}

Here in the first equation we have expanded the definitions of $\eta$ and $\varepsilon$. In the second equation, we have isotoped, preparing for the third equation where we apply the stability equation. In the fourth equation we have isotoped again, preparing for the fifth equation where we apply the splitting equation. We are left with $\id_E$ as desired. The proof of the other triangle identity is similar.
    
$(c) \Rightarrow (b)$: Let $\eta: S \to E \wedge E^\dual$ and $\varepsilon: E^\dual \wedge E \to S$ be the unit and counit of the duality between $E$ and $E^\dual$. Define $i: E \to S$ to be $(\id_E \wedge \varepsilon\beta_{E,E^\dual})((r \wedge \id_E)\inv \wedge \id_{E^\dual})(\id_E \wedge \eta)$. We claim that $(E,i)$ is an open idempotent, i.e. that $\id_E \wedge i$ is an isomorphism. Precomposing the isomorphism $\id_E \wedge r$, it will suffice to show that $\id_E \wedge ir$ is an isomorphism:

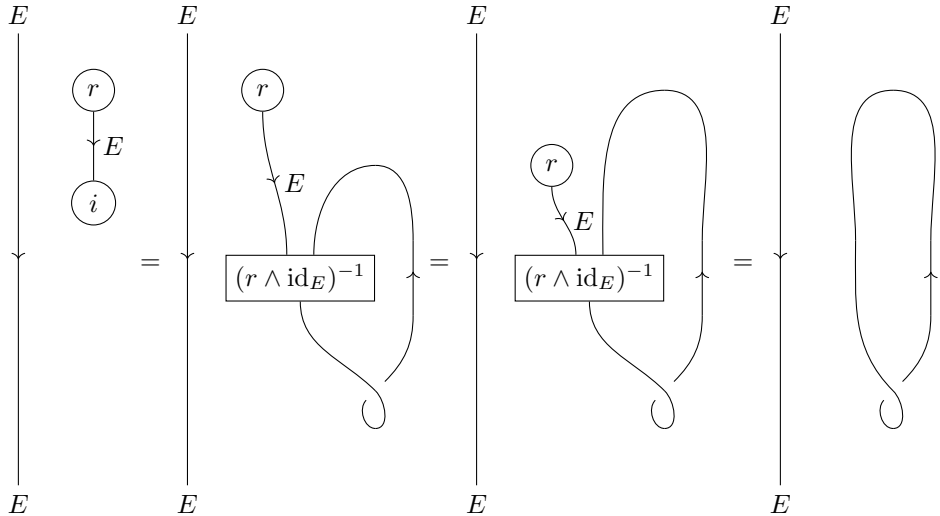
\begin{equation*}
    \begin{tikzpicture}[baseline=(base.south)]
        \node (base) at (0,-3.25) {};
        
        \node (Etop) at (0,0) {$E$};

        \node (Ebot) at (0,-6.5) {$E$};
        
        \node (R) [circlemor] at (1,-1) {$r$};
        \node (I) [circlemor] at (1,-2.5) {$i$};
        \draw (Etop) [downward,witharrow] to (Ebot);
        \draw (R) [downward,witharrow] to node [right] {$E$}  (I);
    \end{tikzpicture}
    =
    \begin{tikzpicture}[baseline=(base.south)]
        \node (base) at (0,-3.25) {};
        
        \node (Etop) at (0,0) {$E$};

        \node (Ebot) at (0,-6.5) {$E$};
        
        \node (R) [circlemor] at (1,-1) {$r$};
        \node (Cap) [mycap] at (2.5,-2) {};
        \node (RI) [rectmor] at (1.5,-3.5) {$(r \wedge \id_E)\inv$};
        \node (E1) [dummy] at (3,-3) {};
        \node (E2) [dummy] at (3,-4) {};
        \node (C) [crossing] at (2.5,-5) {};
        \node (Cup) [mycup] at (2.5,-5.5) {};
        
        \draw (Etop) [downward,witharrow] to (Ebot);
        \draw (R) [downward,witharrow] to node [right] {$E$}  (RI.120);
        \draw (Cap.center) [capwdown] to (RI.60);
        \draw (Cap.center) [capedown] to (E1.center);
        \draw (E1.center) [downward,withbackarrow] to (E2.center);
        \draw (RI) [downcrossnw] to (C.center);
        \draw (E2.center) [downcrossne] to (C.north east);
        \draw (C.south west) [sw_w] to (Cup.center);
        \draw (C.center) [se_e] to (Cup.center);
    \end{tikzpicture}
    =
    \begin{tikzpicture}[baseline=(base.south)]
        \node (base) at (0,-3.25) {};
        
        \node (Etop) at (0,0) {$E$};

        \node (Ebot) at (0,-6.5) {$E$};
        
        \node (R) [circlemor] at (1,-2) {$r$};
        \node (Cap) [mycap] at (2.5,-1) {};
        \node (RI) [rectmor] at (1.5,-3.5) {$(r \wedge \id_E)\inv$};
        \node (E1) [dummy] at (3,-3) {};
        \node (E2) [dummy] at (3,-4) {};
        \node (C) [crossing] at (2.5,-5) {};
        \node (Cup) [mycup] at (2.5,-5.5) {};
        
        \draw (Etop) [downward,witharrow] to (Ebot);
        \draw (R) [downward,witharrow] to node [right] {$E$}  (RI.120);
        \draw (Cap.center) [capwdown] to (RI.60);
        \draw (Cap.center) [capedown] to (E1.center);
        \draw (E1.center) [downward,withbackarrow] to (E2.center);
        \draw (RI) [downcrossnw] to (C.center);
        \draw (E2.center) [downcrossne] to (C.north east);
        \draw (C.south west) [sw_w] to (Cup.center);
        \draw (C.center) [se_e] to (Cup.center);
    \end{tikzpicture}
    =
    \begin{tikzpicture}[baseline=(base.south)]
        \node (base) at (0,-3.25) {};
        
        \node (Etop) at (0,0) {$E$};

        \node (Ebot) at (0,-6.5) {$E$};
        
        \node (Cap) [mycap] at (1.5,-1) {};
        \node (E1) [dummy] at (2,-3) {};
        \node (E2) [dummy] at (2,-4) {};
        \node (E3) [dummy] at (1,-3) {};
        \node (C) [crossing] at (1.5,-5) {};
        \node (Cup) [mycup] at (1.5,-5.5) {};
        
        \draw (Etop) [downward,witharrow] to (Ebot);
        \draw (Cap.center) [capwdown] to (E3.center);
        \draw (Cap.center) [capedown] to (E1.center);
        \draw (E1.center) [downward,withbackarrow] to (E2.center);
        \draw (E3.center) [downcrossnw] to (C.center);
        \draw (E2.center) [downcrossne] to (C.north east);
        \draw (C.south west) [sw_w] to (Cup.center);
        \draw (C.center) [se_e] to (Cup.center);
    \end{tikzpicture}
\end{equation*}

Here in the first equation we have expanded the definition of $i$. In the second equation we have isotoped, setting up the third equation where we cancel inverse morphisms. We continue:

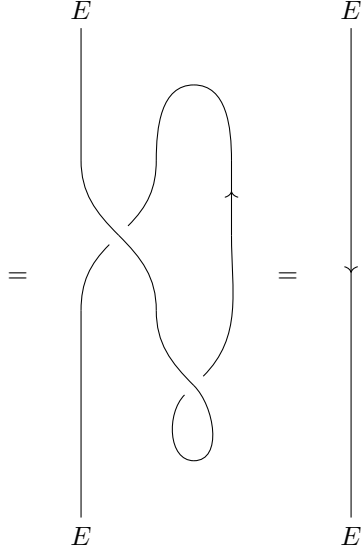
\begin{equation*}
   = \quad
    \begin{tikzpicture}[baseline=(base.south)]
        \node (base) at (0,-3.5) {};
        
        \node (Etop) at (0,0) {$E$};

        \node (Ebot) at (0,-7) {$E$};
        
        \node (Cap) [mycap] at (1.5,-1) {};
        \node (E5) [dummy] at (0,-2) {};
        \node (E6) [dummy] at (0,-4) {};
        \node (E1) [dummy] at (2,-2) {};
        \node (E2) [dummy] at (2,-3) {};
        \node (C2) [crossing] at (0.5,-3) {};
        \node (E3) [dummy] at (1,-2) {};
        \node (E4) [dummy] at (1,-4) {};
        \node (C1) [crossing] at (1.5,-5) {};
        \node (Cup) [mycup] at (1.5,-6) {};
        
        \draw (Etop) [downward] to (E5.center);
        \draw (Cap.center) [capwdown] to (E3.center);
        \draw (Cap.center) [capedown] to (E1.center);
        \draw (E5.center) [downcrossnw] to (C2.center);
        \draw (E3.center) [downcrossne] to (C2.north east);
        \draw (E1.center) [downward,withbackarrow] to (E2.center);
        \draw (C2.south west) [crossswdown] to (E6.center);
        \draw (C2.center) [crosssedown] to (E4.center);
        \draw (E2.center) [downcrossne] to (C.north east);
        \draw (E6.center) [downward] to (Ebot);
        \draw (E4.center) [downcrossnw] to (C1.center);
        \draw (C1.south west) [sw_w] to (Cup.center);
        \draw (C1.center) [se_e] to (Cup.center);
    \end{tikzpicture}
    \quad = \quad
    \begin{tikzpicture}[baseline=(base.south)]
        \node (base) at (0,-3.5) {};
        
        \node (Etop) at (0,0) {$E$};

        \node (Ebot) at (0,-7) {$E$};
        

        \draw (Etop) [downward,witharrow] to (Ebot);
    \end{tikzpicture}
\end{equation*}

Here in the first equation, we have used the fact that $E$ has trivial braiding (because it is a closed idempotent, applying (1)). In the second equation, we have isotoped, leaving us with $\id_E$, which is an isomorphism as desired.

This leaves the last two statements of (2). We have already shown that $E$ is self-dual. For the last statement, suppose that $(E,r,i)$ and $(E,r,i')$ are both clopen. Then we calculate:

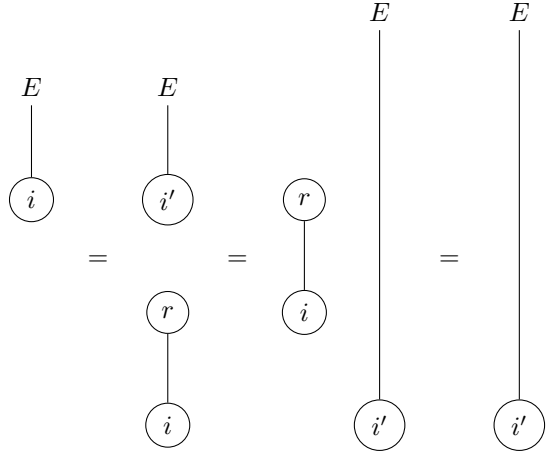
\begin{equation*}
    \begin{tikzpicture}[baseline=(base.south)]
        \node (base) at (0,-2.25) {};
        
        \node (Etop) at (0,0) {$E$};
        \node (I) [circlemor] at (0,-1.5) {$i$};



        \draw (Etop) [downward] to (I);
    \end{tikzpicture}
   \quad = \quad
    \begin{tikzpicture}[baseline=(base.south)]
        \node (base) at (0,-2.25) {};
        
        \node (Etop) at (0,0) {$E$};
        \node (Ip) [circlemor] at (0,-1.5) {$i'$};
        \node (R) [circlemor] at (0,-3) {$r$};
        \node (I) [circlemor] at (0,-4.5) {$i$};



        \draw (Etop) [downward] to (Ip);
        \draw (R) [downward] to (I);
    \end{tikzpicture}
   \quad = \quad
    \begin{tikzpicture}[baseline=(base.south)]
        \node (base) at (0,-2.25) {};
        
        \node (Etop) at (1,1) {$E$};
        \node (Ip) [circlemor] at (1,-4.5) {$i'$};
        \node (R) [circlemor] at (0,-1.5) {$r$};
        \node (I) [circlemor] at (0,-3) {$i$};



        \draw (Etop) [downward] to (Ip);
        \draw (R) [downward] to (I);
    \end{tikzpicture}
    \quad = \quad
    \begin{tikzpicture}[baseline=(base.south)]
        \node (base) at (0,-2.25) {};
        
        \node (Etop) at (0,1) {$E$};
        \node (Ip) [circlemor] at (0,-4.5) {$i'$};



        \draw (Etop) [downward] to (Ip);
    \end{tikzpicture}
\end{equation*}
    Here in the first equation we have used the splitting equation $ri' = \id_E$. In the second equation, we have isotoped, so that in the third equation we may use the stability equation $ir \wedge \id_E = \id_E$. The equation $i=i'$ results.
    
(3): Let $\eta: S \to T^\dual \wedge T$ be the unit and and $\varepsilon: T \wedge T^\dual \to S$ the counit of the duality between $T$ and $T^\dual$; let $t: T \to T$ be the twist. We define $r: S \to T^\dual \wedge T$ to be $(\id_{T^\dual} \wedge t)\eta$ and $i: T^\dual \wedge T \to S$ to be $\varepsilon \beta_{T^\dual,T}$. We first verify the splitting equation $ri = \id_{T^\dual \wedge T}$:

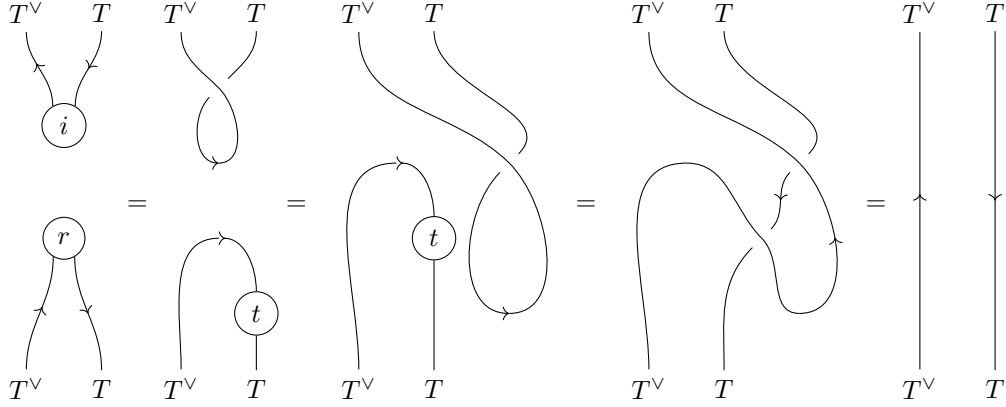
\begin{equation*}
    \begin{tikzpicture}[baseline=(base.south)]
        \node (base) at (0,-2.5) {};
        
        \node (Tdtop) at (0,0) {$T^\dual$};
        \node (Ttop) at (1,0) {$T$};

        \node (Tdbot) at (0,-5) {$T^\dual$};
        \node (Tbot) at (1,-5) {$T$};
        
        \node (I) [circlemor] at (.5,-1.5) {$i$};
        \node (R) [circlemor] at (.5,-3) {$r$};
        
        \draw (Tdtop) [downward,withbackarrow] to (I.120);
        \draw (Ttop) [downward,witharrow] to (I.60);
        \draw (R.-120) [downward,withbackarrow] to (Tdbot);
        \draw (R.-60) [downward,witharrow] to (Tbot);
    \end{tikzpicture}
    =
    \begin{tikzpicture}[baseline=(base.south)]
        \node (base) at (0,-2.5) {};
        
        \node (Tdtop) at (0,0) {$T^\dual$};
        \node (Ttop) at (1,0) {$T$};

        \node (Tdbot) at (0,-5) {$T^\dual$};
        \node (Tbot) at (1,-5) {$T$};
        
        \node (C0) [crossing] at (.5,-1) {};
        \node (Cup) [mycup] at (.5,-2) {};
        \node (Cap) [mycap] at (.5,-3) {};
        \node (Twist) [circlemor] at (1,-4) {$t$};

        \draw (Tdtop) [downcrossnw] to (C0.center);
        \draw (Ttop) [downcrossne] to (C0.north east);
        \draw (C0.south west) [sw_w,->] to (Cup.center);
        \draw (C0.center) [se_e] to (Cup.center);
        \draw (Cap.center) [capwdown] to (Tdbot);
        \draw (Cap.center) [capedown,>-] to (Twist);
        \draw (Twist) [downward] to (Tbot);
    \end{tikzpicture}
    =
    \begin{tikzpicture}[baseline=(base.south)]
        \node (base) at (0,-2.5) {};
        
        \node (Tdtop) at (0,0) {$T^\dual$};
        \node (Ttop) at (1,0) {$T$};

        \node (Tdbot) at (0,-5) {$T^\dual$};
        \node (Tbot) at (1,-5) {$T$};
        
        \node (C0) [crossing] at (2,-2) {};
        \node (Cup) [mycup] at (2,-4) {};
        \node (Cap) [mycap] at (.5,-2) {};
        \node (Twist) [circlemor] at (1,-3) {$t$};

        \draw (Tdtop) [downcrossnw] to (C0.center);
        \draw (Ttop) [downcrossne] to (C0.north east);
        \draw (C0.south west) [sw_w,->] to (Cup.center);
        \draw (C0.center) [se_e] to (Cup.center);
        \draw (Cap.center) [capwdown] to (Tdbot);
        \draw (Cap.center) [capedown,>-] to (Twist);
        \draw (Twist) [downward] to (Tbot);
    \end{tikzpicture}
    =
    \begin{tikzpicture}[baseline=(base.south)]
        \node (base) at (0,-2.5) {};
        
        \node (Tdtop) at (0,0) {$T^\dual$};
        \node (Ttop) at (1,0) {$T$};

        \node (Tdbot) at (0,-5) {$T^\dual$};
        \node (Tbot) at (1,-5) {$T$};
        
        \node (C0) [crossing] at (2,-2) {};
        \node (Cup) [mycup] at (2,-4) {};
        \node (Cap) [mycap] at (.5,-2) {};
        \node (C1) [crossing] at (1.5,-3) {};

        \draw (Tdtop) [downcrossnw] to (C0.center);
        \draw (Ttop) [downcrossne] to (C0.north east);
        \draw (C0.south west) [sw_ne,witharrow] to (C1.north east);
        \draw (C0.center) [se_e,withbackarrow] to (Cup.center);
        \draw (Cap.center) [capwdown] to (Tdbot);
        \draw (Cap.center) [e_nw] to (C1.center);
        \draw (C1.south west) [crossswdown] to (Tbot);
        \draw (C1.center) [se_w] to (Cup.center);
    \end{tikzpicture}
    =
    \begin{tikzpicture}[baseline=(base.south)]
        \node (base) at (0,-2.5) {};
        
        \node (Tdtop) at (0,0) {$T^\dual$};
        \node (Ttop) at (1,0) {$T$};

        \node (Tdbot) at (0,-5) {$T^\dual$};
        \node (Tbot) at (1,-5) {$T$};
        

        \draw (Tdtop) [downward,withbackarrow] to (Tdbot);
        \draw (Ttop) [downward,witharrow] to (Tbot);
    \end{tikzpicture}
\end{equation*}

Here in the first equation we have introduced the definitions of $r$ and $i$. In the second equation we have isotoped, so that in the third equation we can use the definition of a twisted-trivial braiding to exchange $t$ for a crossing. In the fourth equation, we have isotoped, leaving $\id_{T^\dual \wedge T}$ as desired.

We now show that and that $ir \wedge \id_T = \id_T$, from which the stability equation $ir \wedge \id_{T^\dual \wedge T} = \id_{T^\dual \wedge T}$ follows.

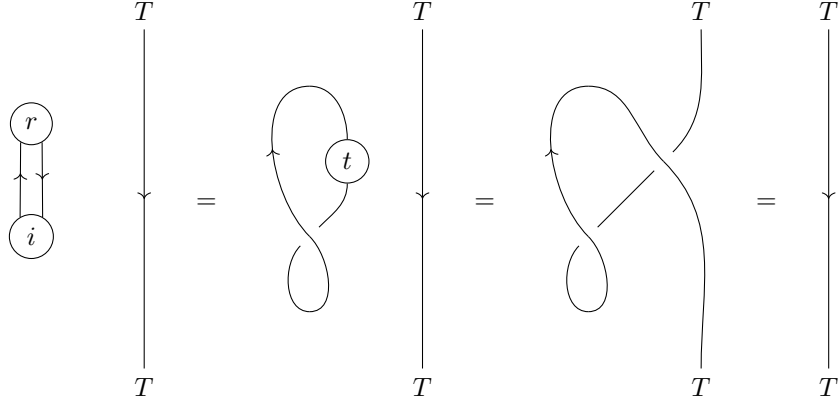
\begin{equation*}
    \begin{tikzpicture}[baseline=(base.south)]
        \node (base) at (0,-2.5) {};
        
        \node (Ttop) at (2,0) {$T$};

        \node (Tbot) at (2,-5) {$T$};
        
        \node (R) [circlemor] at (.5,-1.5) {$r$};
        \node (I) [circlemor] at (.5,-3) {$i$};
        
        \draw (Ttop) [downward,witharrow] to (Tbot);
        \draw (R.-120) [downward,withbackarrow] to (I.120);
        \draw (R.-60) [downward,witharrow] to (I.60);
    \end{tikzpicture}
    \quad = \quad
    \begin{tikzpicture}[baseline=(base.south)]
        \node (base) at (0,-2.5) {};
        
        \node (Ttop) at (2,0) {$T$};

        \node (Tbot) at (2,-5) {$T$};
        
        \node (Cap) [mycap] at (.5,-1) {};
        \node (Twist) [circlemor] at (1,-2) {$t$};
        \node (C0) [crossing] at (.5,-3) {};
        \node (Cup) [mycup] at (.5,-4) {};
        
        \draw (Ttop) [downward,witharrow] to (Tbot);
        \draw (Cap.center) [capedown] to (Twist);
        \draw (Cap.center) [w_nw,withbackarrow] to (C0.center);
        \draw (Twist) [downcrossne] to (C0.north east);
        \draw (C0.center) [se_e] to (Cup.center);
        \draw (C0.south west) [sw_w] to (Cup.center);
    \end{tikzpicture}
    \quad = \quad
    \begin{tikzpicture}[baseline=(base.south)]
        \node (base) at (0,-2.5) {};
        
        \node (Ttop) at (2,0) {$T$};

        \node (Tbot) at (2,-5) {$T$};
        
        \node (Cap) [mycap] at (.5,-1) {};
        \node (C1) [crossing] at (1.5,-2) {};
        \node (C0) [crossing] at (.5,-3) {};
        \node (Cup) [mycup] at (.5,-4) {};
        
        \draw (Ttop) [downcrossne] to (C1.north east);
        \draw (Cap.center) [e_nw] to (C1.center);
        \draw (Cap.center) [w_nw,withbackarrow] to (C0.center);
        \draw (C1.south west) [sw_ne] to (C0.north east);
        \draw (C1.center) [crosssedown] to (Tbot);
        \draw (C0.center) [se_e] to (Cup.center);
        \draw (C0.south west) [sw_w] to (Cup.center);
    \end{tikzpicture}
    \quad = \quad
    \begin{tikzpicture}[baseline=(base.south)]
        \node (base) at (0,-2.5) {};
        
        \node (Ttop) at (0,0) {$T$};

        \node (Tbot) at (0,-5) {$T$};
        
        
        \draw (Ttop) [downward,witharrow] to (Tbot);
    \end{tikzpicture}
\end{equation*}

Here in the first equation, we have introduced the definitions of $r$ and $i$. In the second equation we have used the definition of twisted-trivial braiding to trade $t$ for a crossing. In the third equation, we have isotoped, leaving us with $\id_T$ as desired.
    
(4):The first statement is part of Proposition \ref{prop:smashing-localization}. For the second, let $T$ have twisted-trivial braiding. If $X$ is $T^\dual \wedge T$-stable, then by definition $r \wedge \id_X$ is an isomoprhism $X \cong T^\dual \wedge T \wedge X$. Conversely, from the proof of (\ref{prop:clopens.item:untwist}) we have $ri \wedge \id_T = \id_{T^\dual \wedge T \wedge T}$ and $ir \wedge \id_T = \id_T$, so that $T$ is $T^\dual \wedge T$-stable, and thus for any $Z$, $T\wedge Z$ is also $T^\dual \wedge T$-stable.
\end{proof}

\subsection{Complementary smashing localizations}\label{sec:comp-smash}

In this section, we discuss smashing-cosmashing localizations which have a \defterm{complement} (\cref{def:clopen-comp}), and how they lead to global splittings of the categories they live in (\cref{prop:ordinary-split}).

\begin{Def}\label{def:clopen-comp}
Let $(\calC,\wedge,S)$ be a braided monoidal category with finite biproducts preserved by $\wedge$ in each variable. If $S = X \oplus Y$, where $X$ and $Y$ are clopen idempotents, we say that $X,Y$ are \defterm{complementary} clopen idempotents.
\end{Def}

\begin{prop}\label{prop:ordinary-split}
Let ($\calC,\wedge,S)$ be a braided monoidal semiadditive category, and let $E$ be a clopen idempotent. Then $E$ is a retract of $S$. If a complement $S/E$ to $E$ exists, then $S/E$ is also a clopen idempotent, and thus $E,S/E$ are complementary clopen idempotents. In this case, for all objects $X,Y$, the natural map $\calC(X,Y) \to \calC(E \wedge X, E \wedge Y) \times \calC(S/E \wedge X, S/E \wedge Y)$ is a bijection, and thus $\calC$ splits as a product of its full subcategories of $E$-stable and $E$-torsion objects. In particular, an object is $E$-stable iff it is $S/E$-torsion and vice versa.
\end{prop}
\begin{proof}
For the first statement, let $r,i$ be the maps exhibiting $E$ as a clopen idempotent, and let $s,j$ be maps exhibiting $S/E$ as a complementary idempotent to $E$. Then the splitting equation is satisfied by $s,j$, so we just need to verify the stability equation, which says that $\id_{S/E} \wedge js = \id_{S/E}$. It suffices to verify that $js \wedge js = js$. We may write $js = 1 - ip$, and then we must verify that $1 \wedge 1 - ip \wedge 1 - 1 \wedge ip + ip \wedge ip = 1 - ip \wedge 1$. We know that $ip \wedge 1$, $1 \wedge ip$, and $ip \wedge ip$ are all equal, so this is true.

We have $\calC(X,Y) = \calC((E \wedge X \oplus S/E \wedge X, E \wedge Y \oplus S/E \wedge Y) = \calC(E \wedge X, E \wedge Y) \times \calC(E \wedge X, S/E \wedge Y) \times \calC(S/E \wedge X, E \wedge Y) \times \calC(S/E \wedge X, S/E \wedge Y)$, so we want to show that the cross-terms vanish. Now, since $E$ is self-dual by \cref{prop:clopens}(\ref{prop:clopens.item:equiv}), we may calculate that $\calC(E \wedge X, S/E \wedge Y) = \calC(X, E \wedge S/E \wedge Y)$ and similarly since $S/E$ is self-dual we have $\calC(S/E \wedge X, E \wedge Y) = \calC(X, S/E \wedge E \wedge Y)$, so it suffices to show that $E \wedge S/E = 0$. Now, $E \wedge S/E$ splits the idempotent $ir \wedge js$, which coincides with $irjs \wedge 1 = 0 \wedge 1 = 0$, so indeed $E \wedge S/E = 0$.
\end{proof}

\section[Infinity-Categorical Preliminaries]{$\infty$-Categorical Preliminaries}\label{chap:oo-cat}

This chapter is technical in nature. It mostly consists of assembling some tools from \cite{ha} which will be used in \cref{chap:dual-colim} and \cref{chap:app} to investigate the structure of $\infty$-categories with certain duality and cocompleteness properties. \cref{sec:rev-op} reviews some basics of the theory of operads from \cite{ha}, which we shall use to formalize statements about symmetric monoidal $\infty$-categories in \cref{sec:colims} and the later chapters. \cref{sec:det-ho} discusses several ways to leverage the 1-categorical results of \cref{chap:1-cat} in the $\infty$-categorical context. Of particular interest is \cref{subsec:cob}, where the 1-dimensional cobordism hypothesis is used to justify defining a dualizable object in a symmetric monoidal $\infty$-category in terms of the homotopy category, and to show that there exists a universal way to adjoin duals to a symmetric monoidal $\infty$-category in a variety of contexts. \cref{sec:colims} discusses $\infty$-categories with certain colimits, with a particular eye toward their monoidal properties, as discussed in \cite{ha}. These properties are used in \cref{chap:dual-colim} to construct universal examples of symmetric monoidal $\infty$-cateogories with various duality, cocompleteness, and exactness properties.

\subsection{A quick review of Lurie's theory of operads}\label{sec:rev-op}

In this section, we review some of the basics of Lurie's theory of operads from \cite{ha}. This provides a language to speak rigorously about symmetric monoidal $\infty$-categories.

\begin{Def}[\cite{ha}]
Let $\Fin$ denote the 1-category of finite sets, and let $\Fin_\ast$ denote the 1-category of finite pointed sets. For $n \in \nats$, let $n_+ = \{0,1,\dots,n\}$ denote a set of cardinality $n+1$ with basepoint $0$. A morphism $f: m_+ \to n_+$ of $\Fin_\ast$ is said to be \defterm{active} if $f\inv(0) = \{0\}$, or equivalently if $f$ is in the image of the functor $(-)_+: \Fin \to \Fin_+$ which adds a disjoint basepoint. A morphism $f: m_+ \to n_+$ of $\Fin_\ast$ is said to be \defterm{inert} if $f\inv(i)$ is a singleton for $i \neq 0$, i.e. if $f$ ``only collapses to the basepoint". For $1 \leq i \leq n$, we denote by $\rho_i: m_+ \to 1_+$ the inert morphism which collapses every element to the basepoint except for $i$.
\end{Def}

\begin{Def}[\cite{ha}]
An \defterm{$\infty$-operad} $\calO$ comprises an $\infty$-category $\calO^\otimes$ and a functor $p: \calO^\otimes \to \Fin_\ast$. For $n_+ \in \Fin_\ast$, we denote by $\calO_{n_+}^\otimes$ the fiber of $p$ at $n_+$, and we sometimes abusively write $\calO$ for the \defterm{underlying $\infty$-category $\calO_{1_+}^\otimes$}. This data is required to satisfy the following two conditions:
\begin{enumerate}
    \item For every inert morphism $f: m_+ \to n_+$ of $\Fin_\ast$ and every object $X$ of $\calO^\otimes$ over $m_+$, there exists a cocartesian lift of $f$ emanating from $X$.
\end{enumerate}
In particular, reindexing along an inert morphism $f: m_+ \to n_+$ induces a functor $f_\ast: \calO^\otimes_{m_+} \to \calO^\otimes_{n_+}$. The second condition is:
\begin{enumerate}[resume]
    \item (Segal condition) For each $m_+ \in \Fin_\ast$, the functors $(\rho_1)_\ast,\dots, (\rho_m)_\ast$ induce an equivalence of categories $\calO^\otimes_{m_+} \to \prod_{i=1}^m \calO^\otimes_{1_+}$. 
\end{enumerate}
A \defterm{morphism of $\infty$-operads} $\calO \to \calP$ is a functor $\calO^\otimes \to \calP^\otimes$ over $\Fin_\ast$ which preserves cocartesian lifts of inert morphisms.
\end{Def}

\begin{rmk}
Let $\calO$ be an $\infty$-operad. Then the Segal condition allows us to view any object $X \in \calO^\otimes_{n_+}$ as an $n$-tuple of objects of $\calO = \calO_{1_+}$.
\end{rmk}

\begin{rmk}
Let $\calO$ be an $\infty$-operad such that $\calO^\otimes$ is a 1-category. Then $\calO$ corresponds to a \defterm{colored symmetric operad} in the usual sense, also known as a \defterm{symmetric multicategory}. The correspondence is as follows. Given a symmetric multicategory $\calC$, define $\calO^\otimes \to \Fin_\ast$ to be the category such that $\calO^\otimes_{n_+} = \calC^n$, and with homsets defined by $\calO^\otimes((C_1,\dots,C_m),(D_1,\dots,D_n)) := \prod_{j=1}^n \calC(C_1,\dots,C_m;D_j)$; composition and the functor to $\Fin_\ast$ are defined in the obvious way. $\calO^\otimes$ is called the \defterm{category of operators} of $\calC$, and is an $\infty$-operad. This construction is part of an equivalence of categories.
\end{rmk}

\begin{rmk}
Let $\calO$ be an $\infty$-operad. We see from the construction of the category of operators that reindexing along an inert morphism should be viewed as the operation of \emph{forgetting} certain objects from a tuple of objects (namely those objects whose index is mapped to the basepoint).
\end{rmk}

\begin{Def}[\cite{ha}]
A \defterm{symmetric monoidal $\infty$-category} $\calC$ comprises an $\infty$-category $\calC^\otimes$ and a functor $p: \calC^\otimes \to \Fin_\ast$. For $n_+ \in \Fin_\ast$, we denote by $\calC_{n_+}^\otimes$ the fiber of $p$ at $n_+$, and we sometimes abusively write $\calC$ for the \defterm{underlying $\infty$-category $\calC_{1_+}^\otimes$}. This data is required to satisfy the following two conditions:
\begin{enumerate}
    \item $p$ is a cocartesian fibration.
\end{enumerate}
In particular, reindexing along any morphism $f: m_+ \to n_+$ induces a functor $f_\ast: \calC^\otimes_{m_+} \to \calC^\otimes_{n_+}$. The second condition is:
\begin{enumerate}[resume]
    \item (Segal condition) For each $m_+ \in \Fin_\ast$, the functors $(\rho_1)_\ast,\dots, (\rho_m)_\ast$ induce an equivalence of categories $\calC^\otimes_{m_+} \to \prod_{i=1}^m \calC^\otimes_{1_+}$. 
\end{enumerate}
A \defterm{lax symmetric monoidal functor} between symmetric monoidal $\infty$-categories is a morphism of the underlying operads. A \defterm{strong symmetric monoidal functor} between symmetric monoidal $\infty$-categories is functor over $\Fin_\ast$ which preserves cocartesian lifts of all morphisms. We let $\SMC$ denote the $\infty$-category of symmetric monoidal $\infty$-categories and strong symmetric monoidal functors.
\end{Def}

\begin{rmk}
Let $\calC$ be a symmetric monoidal $\infty$-category. Then in particular, $\calC$ is an $\infty$-operad. In addition to reindexing along inert morphisms (corresponding to forgetting objects from a tuple), $\calC^\otimes$ permits reindexing along \emph{active} morphisms. The interpretation of these reindexing functors is as follows. Reindexing along the unique active morphism $2_+ \to 1_+$ corresponds to the tensor product functor $\otimes: \calC \times \calC \to \calC$. Reindexing along the unique morphism $0_+ \to 1_+$ corresponds to the the unit object, viewed as a functor $\ast = \calC^0 \to \calC$. Reindexing along a general surjective active morphism $f: m_+ \to n_+$ corresponds to tensoring together everything in each fiber of $f$, and reindexing along a general injective active morphism $f: m_+ \to n_+$ corresponds to inserting copies of the unit at all tuple indices not in the image of $f$. Reindexing along a general active morphism $f : m_+ \to n_+$ is a composite of these two operations.
\end{rmk}

\begin{thm}[\cite{ha}]
There is an equivalence of $\infty$-categories between symmetric monoidal $\infty$-categories and strong symmetric monoidal functors on the one hand, and algebras in $\Cat_\infty$ for the $E_\infty$ operad on the other.
\end{thm}

\begin{Def}
Let $\calD$ be an $\infty$-operad. We say that $\calD$ is \defterm{unital} if the nullary operation is representable, i.e. if $\calD^\otimes(\ast,-): \calD^\otimes_{1_+} \to \Top$ is represntable, where $\ast \in \calD^\otimes_{0_+}$ is the unique nullary object. We say that $\calD$ is \defterm{exponentially closed} if, for every $D \in \calD$, the functor $D \otimes (-): \calD \to \calD$ has a right adjoint $F(D,-): \calD \to \calD$. Then $F(D,D')$ is called the \defterm{internal hom} of $D,D' \in \calD$.

Let $\calC$ be a symmetric monoidal $\infty$-category and $\calD \subseteq \calC$ a full suboperad. We say that $\calD$ is a \defterm{$\otimes$-ideal} in $\calC$ if $C \in \calC, D \in \calD \Rightarrow C \otimes D \in \calD$, where $\otimes$ is taken in $\calC$. More generally, we say that $\calD$ is \defterm{closed under $\otimes$} in $\calC$ if $D,D' \in \calD \Rightarrow D \otimes D' \in \calD$. If $\calC$ is exponentially closed with internal hom $F$, then we say that $\calD \subseteq \calC$ is a \defterm{exponential ideal} if $F(C,D) \in \calD$ for all $C \in \calC, D \in \calD$, and \defterm{closed under exponentiation} if $F(D,D') \in \calD$ for all $D,D' \in \calD$.
\end{Def}

\begin{lem}\label{lem:sub-mon}
Let $\calC$ be a symmetric monoidal $\infty$-category which is exponentially closed and let $\calD \subseteq \calC$ be a full suboperad. Suppose that $\calD$ closed under exponentiation and closed under $\otimes$ in $\calC$, and that $\calD$ is unital. Then $\calD$ is a symmetric monoidal $\infty$-category, and the inclusion $\calD \to \calC$ is a lax monoidal functor preserving $\otimes$.
\end{lem}
\begin{proof}
Because $\calD \subseteq \calC$ is a full suboperad, the functor $\calD^\otimes \to \Fin_\ast$ is an isofibration. It will suffice to verify that cocartesian lifts exist for all morphisms of $\Fin_\ast$. Because cocartesian morphsims are closed under composition, it will suffice to show that $\calD$ has cocartesian lifts over inert morphisms, active surjections, and injections in $\Fin_\ast$ separately. Cocartesian lifts over inert morphisms exist by virtue of $\calD$ being an operad (and they are computed as in $\calC^\otimes$ by virtue of $\calD \subseteq \calC$ being a full suboperad inclusion). Because $\calD \subseteq \calC$ is closed under $\otimes$, it follows that cocartesian lifts over active surjective morphisms of $\Fin_\ast$ also exist, computed as in $\calC^\otimes$. Because cocartesian morphisms are closed under composition, it will now suffice to exhibit cocartesian lifts over the injective map $i_{m+1}: m_+ \to (m+1)_+$ which misses $m+1$, for each $m \in \nats$. Let $I \in \calC, J \in \calD$ be the respective units and let $F$ be the exponential in $\calC$. We claim that a cocartesian lift of $D_1,\dots,D_m$ over $i_{m+1}$ is given by $D_1,\dots,D_m,J$. That is, for $f: m+1 \to n$, we claim that the map $\calD^\otimes_f(D_1,\dots,D_m,J;\vec D') \to \calD^\otimes_{fi_{m+1}}(D_1,\dots,D_m;\vec D')$ is an isomorphism. By the Segal condition, it will suffice to consider the case where $n=1$ and $f$ is active. In this case, we have

\begin{align*}
    \calD^\otimes_f(D_1,\dots,D_m,J;D') &= \calC^\otimes_f(D_1,\dots,D_m,J;D') \\
    &= \calC(D_1 \otimes \cdots \otimes D_m \otimes J, D') \\
    &= \calC(J,F(D_1\otimes \cdots \otimes D_m,D')) \\
    &= \calD(J,F(D_1\otimes \cdots \otimes D_m,D')) \\
    &= \calD^\otimes_{0}(\ast;F(D_1\otimes \cdots \otimes D_m,D')) \\
    &= \calC^\otimes_{0}(\ast;F(D_1\otimes \cdots \otimes D_m,D')) \\
    &= \calC(I,F(D_1\otimes \cdots \otimes D_m,D')) \\
    &= \calC(D_1\otimes \cdots \otimes D_m \otimes I,D') \\
    &= \calC(D_1\otimes \cdots \otimes D_m,D') \\
    &= \calC^\otimes_{fi_{m+1}}(D_1,\dots,D_m;D') \\
    &= \calD^\otimes_{fi_{m+1}}(D_1,\dots,D_m;D')
\end{align*}
as desired. Here $\ast$ denotes the unique object of $\calC^\otimes_{0_+} =\calD^\otimes_{0_+}$ and $0$ is the unique morphism $0_+ \to 1_+$ in $\Fin_\ast$. In this argument, we have used that $\calD^\otimes \subseteq \calC^\otimes$ is a full subcategory, the existence of cocartesian lifts for $\calC^\otimes \to \Fin_\ast$, the exponential closedness of $\calC$, the fully faithfulness of $\calD \to \calC$ and the fact that $\calD$ is closed under $\otimes$ and exponentiation, the universal proeprty of $J$, the fully faithfulness of $\calD^\otimes \to \calC^\otimes$, the unitality of $\calC$, the exponential closedness of $\calC$, unitality of $\calC$, the definition of $D_1 \otimes \cdots \otimes D_m$, and the fully faithfulness of $\calD^\otimes \to \calC^\otimes$ again.
\end{proof}

\begin{Def}
Let $\calC$ be an $\infty$-category and $\calD$ a full subcategory. A functor $L: \calC \to \calD$ is said to be a \defterm{localization functor} if it is left adjoint to the inclusion $i: \calD \to \calC$ and moreover $Li$ is naturally equivalent to the identity on $\calD$. The composite $iL: \calC \to \calC$ may also be referred to as a \defterm{localization functor}. An \defterm{$L$-local morphism} is a morphism in $\calC$ whose image under $L$ is an equivalence. An \defterm{$L$-local functor} out of $\calC$ is a functor taking $L$-local morphisms to equivalences; these span a full subcategory $\Fun_{L\textrm{-loc}}(\calC,\calE)$ of $\Fun(\calC,\calE)$. 

Let $\calO$ be an $\infty$-category, $\calC \to \calO$ a category over $\calO$, and $\calD$ a full subcategory of $\calC$ over $\calO$. An adjunction functor $L: \calC^\to_\leftarrow \calD:i$ is said to be a \defterm{localization over $\calO$} if it is a localization and $L,i$ are over $\calO$. We do not require that the unit and counit of the adjunction $L \dashv i$ be vertical. We may call the endofunctor $iL$ a \defterm{localization functor over $\calO$}, and an \defterm{$L$-local functor over $\calO$} out of $\calC$ is a functor over $\calO$ which is $L$-local; these span a full subcategory $\Fun_{L\textrm{-loc},\calO}(\calC,\calE)$ of $\Fun_\calO(\calC,\calE)$.

Let $\calO$ be an $\infty$-category, $\calC \to \calO$ a cocartesian fibration, and $\calD$ a full subcategory of $\calC$ over $\calO$ which is cocartesian over $\calO$. We do not require that the inclusion $\calD \to \calC$ preserve cocartesian edges. A functor $L: \calC \to \calD$ is said to be a \defterm{cocartesian localization functor over $\calO$} if is a localization functor, it is over $\calO$, and it preserves cocartesian edges. We may refer to $iL$ as a \defterm{cocartesian localization functor over $\calO$}. An \defterm{$L$-local cocartesian functor over $\calO$} is an $L$-local functor which is cocartesian over $\calO$; these span a full subcategory $\Fun^{\cocart}_{L\textrm{-loc},\calO}(\calC,\calE)$ of $\Fun_\calO^\cocart(\calC,\calE)$.

In the case where $\calO$ is an $\infty$-operad, recall that the cocartesian fibrations $\calC, \calD$ over $\calO$ are called \defterm{$\calO$-monoidal $\infty$-categories}. In this case, we correspondingly refer to a cocartesian localization functor over $\calO$ as an \defterm{$\calO$-monoidal localization}, and we refer to a cocartesian localization functor over $\calO$ as an \defterm{$L$-local $\calO$-monoidal functor}; these span a full subcategory $\Fun_{L,\calO}^\cocart(\calC,\calE)$ of $\Fun_\calO(\calC,\calE)$.
\end{Def}

\begin{rmk}
The notion of a symmetric monoidal $\infty$-category is recovered as an $\calO$-monoidal $\infty$-category when $\calO^\otimes = \Fin_\ast$ is the $E_\infty$-operad.
\end{rmk}

\begin{prop}
\begin{enumerate}
    \item If $L :\calC \to \calD$ is a localization functor, it has the following universal property. For any $\infty$-category $\calE$, precomposing $L$ induces an equivalence between $\Fun(\calD, \calE)$ and $\Fun_{L\textrm{-loc}}(\calC,\calE)$, with inverse given by restriction to $\calD$.
    \item If $L: \calC^\to_\leftarrow \calD:i$ is a localization functor over $\calO$, it has the following universal property. For any $\infty$-category $\calE$ over $\calO$, the precomposing $L$ induces an equivalence between $\Fun_\calO(\calD, \calE)$ and $\Fun_{L\textrm{-loc},\calO}(\calC,\calE)$, with inverse given by restricting to $\calD$.
    \item If $L: \calC \to \calD$ is a cocartesian localization functor over $\calO$, it has the following universal property. For any $\infty$-category $\calE$ cocartesian over $\calO$, precomposing $L$ induces an equivalence between $\Fun^\cocart_\calO(\calD, \calE)$ and $\Fun^\cocart_{L,\calO}(\calC,\calE)$, with inverse given by restriction to $\calD$.
    
    In particular, an $\calO$-monoidal localization functor $L: \calC \to \calD$ induces an equivalence between $\Fun_\calO(\calD,\calE)$ and $\Fun_{L,\calO}(\calC,\calE)$, with inverse given by restriction to $\calD$.
\end{enumerate}
\end{prop}
\begin{proof}
$(1)$ is \cite[Proposition 5.2.7.12]{htt}.

$(2)$ follows from $(1)$ as soon as we note that since $L$ and the inclusion $i: \calD \to \calC$ are over $\calO$, precomposing them preserves the property of being over $\calO$.

$(3)$ follows from $(2)$ as soon as we check the following. Let $F: \calD \to \calE$ be a functor over $\calO$. We must show that if $FL: \calC \to \calE$ preserves cocartesian morphisms, then $F$ preserves cocartesian morphisms. (Because we do not assume that $i$ preserves cocartesian morphisms, this does not immediately follow from the fact that $F \simeq FLi$.) To see this, it suffices to show that any cocartesian morphism in $\calD$ is isomorphic to the image of a cocartesian morphism in $\calC$. This is straightforward: let $f$ be a morphism in $\calO$, let $D \in \calD$ lie over the domain of $f$, and let $f_\ast^D$ be a cocartesian lift of $f$ out of $D$. Let $f_\ast^{iD}$ be a cocartesian lift of $f$ out of $iD \in \calC$. We have a factorization $if_\ast^D = \overline{if_\ast^D}f_\ast^{iD}$ where $\overline{if_\ast^D}$ is vertical. Applying $L$, we have $Lif_\ast^D = L \overline{if_\ast^D}Lf_\ast^{iD}$, with $L \overline{if_\ast^D}$ vertical. Now $Li f_\ast^D$ is isomorphic to $f_\ast^D$ and hence cocartesian, while $Lf_\ast^{iD}$ is cocartesian because $L$ preserves cocartesian morphisms. Hence the vertical factorization $L \overline{if_\ast^D}$ must be an isomorphism. So $f_\ast^D$ is isomorphic to $Lif_\ast^D$, which in turn isomorphic to $Lf_\ast^{iD}$, which is the image under $L$ of a cocartesian morphism, completing the proof.
\end{proof}

\begin{lem}\label{lem:comm-init}
Let $(\calC,\wedge,S)$ be a symmetric monoidal $\infty$-category. Then $S$ is naturally (in fact, uniquely) an $E_\infty$-algebra in $\calC$, and as such is the initial object of the $\infty$-category of $E_\infty$-algebras in $\calC$.
\end{lem}
\begin{proof}
See \cite[Corollary 3.2.1.9]{ha}.
\end{proof}

\subsection{Structures detected in the homotopy category}\label{sec:det-ho}

In this section, we discuss certain properties of $\infty$-categories which may be detected by passing to the homotopy category. Notably, in \cref{subsec:cob}, the 1-dimensional cobordism hypothesis is used to show that duals may be freely adjoined to the objects of a symmetric monoidal $\infty$-category in a wide variety of contexts. Moreover, in \cref{subsec:fun-split} the category-wide splittings arising from dualizable objects with twisted-trivial braiding (\cref{sec:tt-dual} and \cref{sec:comp-smash}) are lifted to the $\infty$-categorical setting.

\subsubsection{The Cobordism Hypothesis}\label{subsec:cob}

In this subsection, we use the 1-dimensional cobordism hypothesis to justify defining a dualizable object in a symmetric monoidal $\infty$-category to be an object whose image in the homotopy category is dualizable. Moreover, the 1-dimensional cobordism hypothesis is used to verify the hypotheses of the Adjoint Functor Theorem and conclude that in many contexts, duals maybe freely adjoined to all objects of a symmetric monoidal $\infty$-category (\cref{cor:freecats}).

Let $\SMD$ denote the full sub-$\infty$-category of $\SMC$ spanned by the small symmetric monoidal $\infty$-categories which have duals for objects. We would like to show that the inclusion $\SMD \to \SMC$, and variants where the categories involved have various colimits, have left adjoints. This is not difficult once we avail ourselves of a form of the 1-dimensional cobordism hypothesis.

\begin{thm}\label{thm:cobhyp}[{\cite{lurie-cob,harpaz-cob}}]
Let $\calC$ be a symmetric monoidal $\infty$-category. Then the homotopy category of $\calC$ has duals if and only if $\calC$ is right orthogonal to the inclusion $\FinBij \to \Bord_1^\fr$ of the symmetric monoidal category of finite sets and bijections (under disjoint union) into the symmetric monoidal $\infty$-category of framed 1-dimensional bordisms. Moreover, an object $C \in \calC$ is dualizable if and only if the symmetric monoidal functor $\FinBij \to \calC$ classifying $C$ extends to a functor $\Bord_1^\fr \to \calC$, in which case the extension is unique up to a contractible space of choices.
\end{thm}

\begin{cor}\label{cor:freecats}
The inclusion $\SMD \to \SMC$ has a left adjoint $\freeduals$ exhibiting $\SMD$ as a localization of $\SMC$. If $\calK \to \SMC$ is a right-adjoint functor from a presentable $\infty$-category, then $\SMD \times_\SMC \calK \to \calK$ has a left adjoint exhibiting $\SMD \times_\SMC \calK$ as a localization of $\calK$.
\end{cor}
\begin{proof}
The second statement follows from the first because limits of presentable categories and right adjoint functors are computed as for the underlying categories.

The first statement holds because by Theorem \ref{thm:cobhyp}, $\SMD$ is characterized as the full subcategory right orthogonal to a small (in fact a singleton) set of morphisms $\{\FinBij \to \Bord_1^\fr\}$, and therefore is an accessible localization of the presentable $\infty$-category $\SMC$.
\end{proof}

\begin{cor}\label{cor:pushout}
If $\calC$ is a symmetric monoidal $\infty$-category and $S \subseteq \calC$ is a set of objects, then there is a universal symmetric monoidal $\infty$-category $\freeduals_S \calC$ receiving a symmetric monoidal functor from $\calC$ carrying the object of $S$ to dualizable objects. If $\calK$ is a cocomplete $\infty$-category admitting a functor $U : \calK \to \SMC$, and if $P\FinBij, P\Bord_1^\fr \in \calK$ are objects such that that $\calK(P\FinBij, -) \cong \SMC(\FinBij , U-)$ and $\calK(P\Bord_1^\fr, -) \cong \SMC(\Bord_1^\fr, U-)$, then for any $K \in \calK$ and $S \subseteq UK$, there is likewise a universal object $\freeduals_S^\calK K$ of $\calK$ admitting a morphism $f : K \to \freeduals_S^\calK K$ such that $Uf : U\calK \to U(\freeduals_C^\calK K)$ carries each $C \in S$ to a dualizable object.
\end{cor}
\begin{proof}
For the first statement, $\freeduals_C \calC$ is computed by a pushout in $\SMC$: $\freeduals_C \calC = \calC \cup_{\amalg_{s \in S}\FinBij} \amalg_{s \in S} \Bord_1^\fr$, where $\amalg_{s \in S}$ denotes a coproduct in $\SMC$. (Note that  pushouts in $\SMC$ are computed via a bar construction in $\Cat_\infty$.) Likewise, for the second statement we have $\freeduals_S^\calK K = K \cup_{\amalg_{s \in S} P\FinBij} \amalg_{s \in S} P\Bord_1^\fr$. (Again, if $\calK$ is an $\infty$-category of commutative algebra objects in some symmetric monoidal $\infty$-category $\calL$, then pushouts in $\calK$ are computed via a bar construction in $\calL$).
\end{proof}

\subsubsection{Certain limits and colimits}\label{subsec:lim-colim-ho}

In this subsection, we discuss certain very basic limit and colimit constructions which may be detected in the homotopy category.

\begin{lem}\label{lem:coprod-hocat}
Let $\calC$ be an $\infty$-category with $K$-colimits, where $K$ is discrete. Then the homotopy category $\ho \calC$ has $K$-colimits, and the canonical functor $\calC \to \ho \calC$ preserves them.
\end{lem}
\begin{proof}
This follows from the fact that $\pi_0: \Top \to \Set$ commutes with products.
\end{proof}

\begin{cor}\label{cor:pt-semiadd-hocat}
Let $\calC$ be an $\infty$-category with an initial object and a terminal object. Then $\calC$ is pointed if and only if $\ho\calC$ is pointed.

Let $\calC$ be a pointed $\infty$-category with finite products and finite coproducts. Then $\calC$ is semiadditive if and only if $\ho\calC$ is semiadditive.

Let $\calC$ be a semiadditive $\infty$-category. Then $\calC$ is additive if and only if $\ho \calC$ is additive.
\end{cor}
\begin{proof}
For the first statement, $\ho \calC$ has an initial and terminal object by \cref{lem:coprod-hocat} and its dual, and $\calC \to \ho\calC$ preserves initial and terminal objects. The map from the latter to the former is an isomorphism in $\calC$ iff it is an isomorphism in $\ho\calC$, verifying the claim.

The proof of the second statement is similar, using the canonical map $A \amalg B \to A \times B$ coming from the pointed structure.

The proof of the third statement is also similar, using the map $\begin{pmatrix}1 & 1 \\ 0 & 1 \end{pmatrix}: A \oplus A \to A \oplus A$.
\end{proof}

\begin{lem}\label{lem:split-cof}
Let $\calC$ be a pointed $\infty$-category, and let $i: E \to S$ be a split monomorphism in $\calC$. If $i$ has a cofiber $q: S \to S/E$, then $q$ is the cokernel of $i$ in the homotopy category $\ho\calC$.
\end{lem}
\begin{proof}
The cofiber sequence $E \to S \to S/E$ induces a fiber sequence $\calC(E,D) \leftarrow \calC(S,D) \leftarrow \calC(S/E,D)$ for any $D \in \calC$. Because $E \to S$ splits, the induced long exact sequence of homotopy groups includes in particular a short exact sequence $\pi_0 \calC(E,D) \leftarrow \pi_0 \calC(S,D) \leftarrow \pi_0 \calC(S/E,D)$, which is to say that $E \to S \to S/E$ is a cofiber sequence in $\ho \calC$ as desired.
\end{proof}

\begin{cor}\label{cor:split-cof}
Let $\calC$ be a pointed $\infty$-category, and suppose that $X^{\overset{r}{\leftarrow}}_{\underset{i}{\to}} Z$ is a retract. If $\calC$ is semiadditive and $X$ is a group object, then any cofiber of $i$ (or fiber of $p$) is a complement to $X$.
\end{cor}
\begin{proof}
This follows from \cref{lem:split-cof} and \cref{lem:cogroup}.
\end{proof}

\begin{lem}\label{lem:dual-stable}
\begin{enumerate}
    \item Let $\calC$ be a symmetric monoidal $\infty$-category. Then the unit object is self-dual.
    \item Let $\calC$ be a symmetric monoidal $\infty$-category. Then any $\otimes$-invertible object is dualizable (with dual given by its inverse).
    \item Let $\calC$ be a symmetric monoidal pointed $\infty$-category. Then the zero object is dualizable.
    \item Let $\calC$ be a symmetric monoidal semiadditive $\infty$-category. Then the dualizable objects are closed under direct sum.
    \item Let $\calC$ be a symmetric monoidal stable $\infty$-category. Then the dualizable objects in $\calC$ are closed under finite limits and finite colimits.
    \item Let $\calC$ be a symmetric monoidal $\infty$-category with split idempotents. Then the dualizable objects in $\calC$ are closed under retracts.
\end{enumerate}
\end{lem}
\begin{proof}
(1) is trivial. (2) is an instance of the fact that any equivalence of categories may be upgraded to an adjoint equivalence. For (3), if $(X,X^\dual,\eta,\varepsilon)$ and $(X',{X'}^\dual,\eta',\varepsilon')$ are duality data, then writint $(X\oplus X') \wedge (X^\dual \oplus {X'}^\dual) = (X \wedge X^\dual) \oplus (X \wedge {X'}^\dual) \oplus ({X'} \wedge X^\dual) \oplus ({X'} \wedge {X'}^\dual)$ and $(X^\dual \oplus {X'}^\dual) \wedge (X \oplus {X'}) = (X^\dual \wedge X) \oplus (X^\dual \wedge {X'}) \oplus ({X'}^\dual \wedge X) \oplus ({X'}^\dual \wedge {X'})$, then $(X \oplus X, X^\dual \oplus {X'}^\dual , [\eta,0,0, \eta'], [\varepsilon,0,0,\varepsilon'])$ is a duality datum.

For (4), it suffices in light of (3) to show that dualizable objects are closed under cofibers and desuspension. The latter follows from closure under tensor and the fact that $\Sigma\inv S = (\Sigma S)^\dual$ is dualizable. For the former, from a cofiber sequence $X \to Y \to Z$ with $X,Y$ dualizable, we obtain a fiber sequence $X^\dual \leftarrow Y^\dual \leftarrow F$, and we claim that $F$ is dual to $Z$. For any $A,B \in \calC$, we have $\calC(Z \otimes A, B) = \calC(Y \otimes A, B) \times_{\calC(X\otimes A, B)} \{0\} = \calC(A,Y^\dual \otimes B)\times_{\calC(A,X^\dual\otimes B)} \{0\} = \calC(A,F \otimes B)$. That is, we have a an adjunction $Z \otimes(-) \calC^\to_\leftarrow \calC: F \otimes(-)$, which by the Yoneda lemma implies that $F$ is right adjoint to $Z$ in $\catdeloop\calC$, i.e. that $F = Z^\dual$ as desired.

For (5), suppose that $A$ has a dual $A^\dual$. Suppose that $B$ is a retract of $A$, splitting the idempotent $e : A \to A$. Then there is a dual idempotent $e^\dual : A^\dual \to A^\dual$, and a splitting $B^\dual$ of this idempotent provides a dual for $B$.
\end{proof}

\subsubsection{Functorial Splitting}\label{subsec:fun-split}

In this subsection we discuss how to use the observations of \cref{chap:1-cat} \cref{subsec:lim-colim-ho} to obtain splittings of entire $\infty$-categories from the existence of dualizable objects with twisted-trivial braiding.

\begin{prop}\label{prop:monoidal-localization}
Let $(\calC, \wedge, S)$ be an $E_n$-monoidal $\infty$-category where $n \geq 2$. Suppose that $E$ is a clopen idempotent. Then
\begin{enumerate}
    \item The full subcategory $\calC_E$ of $E$-stable objects is canonically $E_n$-monoidal, with the localization functor $\calC \to \calC_E$ being an $E_n$-monoidal localization.
    \item If $\calC$ has (co)limits of shape $I$ preserved by $\wedge$ in each variable, then so does $\calC_E$. Moreover, the functors $\calC^\to_\leftarrow \calC_E$ preserve these colimits, and an $E_n$-monoidal functor $\calC_E \to \calD$ preserves $I$-(co)limits iff $\calC \to \calC_E \to \calD$ does.
    \item\label{prop:monoidal-localization.item:split} If $E$ has a complement $S/E$, then the localization functor $\calC \to \calC_E \times \calC_{S/E}$ is an $E_n$-monoidal equivalence.
\end{enumerate}
\end{prop}
\begin{proof}
For (1), by \cite[Proposition 2.2.1.9]{ha}, it suffices to show that $E$-local equivalences are stable under tensoring. This is clear: if $E \wedge f$ and $E \wedge g$ are equivalences, then $E \wedge (f \wedge g) = (E \wedge f) \wedge (E \wedge g)$ is an equivalence.

For (2), note that $\calC_E$ is closed in $\calC$ under $I$-(co)limits. For if we have an $I$-diagram in $\calC_E$ with (co)limit in $\calC$, we may smash the colimit diagram with $E$; by hypothesis this results in a new (co)limit diagram. The base of the original diagram is naturally isomorphic to the new diagram, so the old and new (co)limits are also isomorphic via the natural map, i.e. the (co)limit lies in $\calC_E$. So $i$ preserves and reflects $I$-(co)limits. Moreover $iL$ (which is given by smashing with $E$) also preserves $I$-(co)limits. It follows that $L$ preserves $I$-colimits. Since $i$ and $L$ preserve $I$-(co)limits, (2) follows.

For (3), the functor is a product of $E_n$-monoidal functors and hence also $E_n$-monoidal. That it is an equivalence can be checked at the level of homotopy categories, and follows from Proposition \ref{prop:ordinary-split}.
\end{proof}

\begin{thm}\label{prop:monoidal-splitting}
Let $(\calC,\wedge,S)$ be an $E_n$-monoidal $\infty$-category, where $n \geq 2$.
\begin{enumerate}
\item\label{prop:monoidal-splitting.item:1} If $E$ is a clopen idempotent, then the $E_n$-monoidal localization $\calC \to \calC_E$ is the universal $E_n$-monoidal functor inverting $E$ under $\wedge$.

\item\label{prop:monoidal-splitting.item:2} If $T$ has twisted-trivial braiding and a dual $T^\dual$, then the $E_n$-monoidal localization $\calC \to \calC_{T^\dual \wedge T}$ is the universal $E_n$-monoidal functor inverting $T$ under $\wedge$.

\item\label{prop:monoidal-splitting.item:3} Suppose that $\calC$ has finite biproducts preserved by $\wedge$ in each variable. If $E$ is a clopen idempotent with complement $S/E$, then $\calC \to \calC_{S/E}$ is the universal $E_n$-monoidal functor preserving finite biproducts and sending $E$ to 0.
\end{enumerate}
\end{thm}
\begin{proof}
For (1), $E$ is certainly $\wedge$-invertible in $\calC_E$ -- in fact it is the monoidal unit. By the universal property of $E_n$-monoidal localization, it suffices to verify that an $E_n$-monoidal functor $F: \calC \to \calD$ inverts $E$ under $\wedge$ iff it takes $E$-local morphisms to equivalences. Now, $F$ inverts $E$ under $\wedge$ iff it takes the unit and counit exhibiting $E$ as self-dual to isomorphisms. This unit and counit are none other than the canonical maps $E \to S$ and $S \to E$, so they are inverted iff $F$ takes $E \to S$ to an isomorphism. But becuase $F$ is $E_n$-monoidal, this happens iff $F$ takes $E \wedge X \to X$ to an isomorphism for every $X$, i.e. takes $E$-local morphisms to equivalences as claimed.

For (2), $T$ is certainly $\wedge$-invertible in $\calC_{T^\dual \wedge T}$, with inverse $T^\dual$ -- that is to say, $T^\dual \wedge T$ is the monoidal unit. A functor $F: \calC \to \calD$ inverts $T$ under $\wedge$ if and only if it sends the duality data for $T$ to equivalences, if and only if it sends the clopenness data for $T^\dual \wedge T$ to equivalences, if and only if (by (1)) it factors though $\calC \to \calC_{T^\dual \wedge T}$.

For (3), it suffices to verify that $F: \calC \to \calD$ inverts $S/E$ iff it takes $E$ to $0$. This follows from the fact that in a biproduct $Z = X \oplus Y$, we have $X = 0$ iff $Y \to Z$ is an isomorphism.
\end{proof}

\subsection[Infinity-categories with certain colimits]{$\infty$-Categories with certain colimits}\label{sec:colims}

In this section, we review some material from \cite{ha} on $\infty$-categories with certain colimits, with an eye toward monoidal properties. We also begin considering certain subcategories of such categories, starting with the pointed case. This study will continue in \cref{chap:dual-colim}.

\begin{lem}\label{lem:free-cocomp-k}
Let $\calK$ be a class of small $\infty$-categories and let $\calC$ be a small $\infty$-category. Let $\calR$ be a set of cocones on diagrams with shape as in $\calK$. Then there exists an $\infty$-category $\calP^\calK_\calR(\calC)$ with all colimits of shape in $\calK$, and a functor $\calC \to \calP^\calK_\calR(\calC)$ carrying the cocones of $\calR$ to colimiting cocones, which is universal with these properties.

Moreover, $\calP^\calK_\calR(\calC)$ may be constructed as follows. Let $\Psh(\calC)$ be the $\infty$-category of presheaves on $\calC$, let $S \subseteq \Mor \Psh(\calC)$ be the collection of morphisms $\varinjlim_{k \in K} R(k) \to R(\infty)$ for each $R: K^\triangleright \to \calC$ in $\calR$, and let $L: \Psh(\calC) \to S\inv \Psh(\calC)$ be the localization functor. Then $\calP^\calK_\calR(\calC)$ is the closure of $L(\calC) \subseteq S\inv \Psh(\calC)$ under $\calK$-shaped colimits.
\end{lem}
\begin{proof}
This is \cite[Proposition 5.3.6.2]{htt}. The explicit description is given in the proof.
\end{proof}

\begin{notation}
We will denote $\calP_\calR(\calC):= \calP^\calK_\calR(\calC)$ for $\calK$ the class of \emph{all} small colimits. We will denote $\calP^\calK(\calC) := \calP^\calK_\calR(\calC)$ for $\calR$ the empty class of cocones. Thus $\calP(\calC) = \Psh(\calC)$ denotes $\calP^\calK_\calR(\calC)$ when $\calK$ is the class of all small colimits and $\calR$ is the empty class of cocones.
\end{notation}

\begin{Def}\label{def:cat-k}
Let $\calK$ be a class of small categories. Let $\Cat_\calK$ be the symmetric monoidal $\infty$-category of $\infty$-categories with $\calK$-indexed colimits, where the symmetric monoidal structure is that of \cite[Corollary 4.8.1.4]{ha}.
\end{Def}

\begin{rmk}\label{rmk:cat-k}
According to \cite[Notation 4.8.1.2]{ha}, a morphism in $\Cat_\calK^\otimes$ from $\calC_1,\dots, \calC_n$ to $\calD$ is a functor $\calC_1 \times \cdots \times \calC_n \to \calD$ which preserves $\calK$-indexed colimits separately in each variable. In particular, when $n=0$, a nullary morphism $\ast \to \calD$ is simply an object of $\calD$.

According to \cite[Proposition 4.8.1.3]{ha}, the tensor product in $\Cat_\calK$ may be constructed by setting $\calC \otimes_\calK \calD = \calP^\calK_{\calK \boxtimes \calK}(\calC \times \calD)$, where $\calK \boxtimes \calK$ is the collection of diagrams described in \cite[Notation 4.8.1.7]{ha}. That is, $\calK \boxtimes \calK$ comprises those diagrams of the form $(C_k,D)_{k \in K}$ for $K \in \calK$, with vertex $(\varinjlim_{k \in K} C_k , D)$, as well as those diagrams of the form $(C,D_k)_{k \in K}$ for $K \in \calK$ with vertex $(C,\varinjlim_{k \in K} D_k)$. As recalled in \cref{lem:free-cocomp-k}, the proof of \cite[Proposition 5.3.6.2]{htt} shows this means that $\calC \otimes_\calK \calD$ is the smallest full subcategory of $(\calK \boxtimes \calK)\inv \Psh(\calC \times \calD)$ containing the image $L(\calC \times \calD)$ (where $L: \Psh(\calC \times \calD) \to (\calK \boxtimes \calK)\inv \Psh(\calC \times \calD)$ is the localization functor) and closed under $\calK$-shaped colimits.
\end{rmk}

\begin{notation}\label{notation:2-cat}
For $\calC_1,\dots,\calC_n,\calD \in \Cat_\calK$, we let $\Fun_\calK(\calC_1,\dots,\calC_n;\calD)$ be the $\infty$-category of functors $\calC_1 \times \cdots \times \calC_n \to \calD$ preserving $\calK$-colimits separately in each variable, so that the multi-hom-space $\Cat_\calK(\calC_1,\dots,\calC_n;\calD)$ is the underlying $\infty$-groupoid of $\Fun_\calK(\calC_1,\dots,\calC_n;\calD)$.
\end{notation}

\begin{rmk}\label{rmk:2-cat}
As noted in \cite[4.8.1.6]{ha}, the construction $\Fun_\calK$ of \cref{notation:2-cat} exhibits the symmetric monoidal $\infty$-category $\Cat_\calK$ as an exponentially closed symmetric monoidal $\infty$-category.
\end{rmk}

\begin{rmk}\label{rmk:ha-4819}
There is an equivalence between $E_\infty$-algebra objects in $\Cat_\calK$ and symmetric monoidal $\infty$-categories with $\calK$-indexed colimits over which $\otimes$ distributes.
More generally, there is an equivalence between $E_n$-monoidal $\infty$-categories with $\calK$-colimits over which $\otimes$ distributes, and $E_n$-algebra objects in $\Cat_\calK$. See \cite[4.8.1.9]{ha}.
\end{rmk}

\begin{Def}
Let $\calK$ be a class of small $\infty$-categories containing the empty category, so that any $\calK$-cocomplete category $\calC \in \Cat_\calK$ has an initial object. Let $\Cat_{\calK,\ast} \subset \Cat_\calK$ denote the full sub-operad of $\Cat_\calK$ comprising those $\calK$-cocomplete categories which are pointed (i.e. where the initial object is a zero object).
\end{Def}

\begin{rmk}
As in \cref{rmk:ha-4819},there is an equivalence between $E_\infty$-algebra objects in $\Cat_{\ast,\calK}$ and pointed symmetric monoidal $\infty$-categories with $\calK$-indexed colimits over which $\otimes$ distributes.
\end{rmk}

\begin{lem}\label{lem:pointed-free}
Let $\calK$ be a class of small $\infty$-categories containing the empty category. Then the free pointed $\calK$-cocomplete $\infty$-category on an object is given by $\Top_\ast^\calK$, the smallest full subcategory of $\Top_\ast$ which contains $S^0$ and is closed under $\calK$-colimits. That is, for any $\calD \in \Cat_{\ast,\calK}$, evaluation at $S^0$ determines an equivalence of categories $\Cat_{\ast,\calK}(\Top_\ast^\calK,\calD) \to \calD^\sim$.
\end{lem}
\begin{rmk}
In the following proof, we heavily use \cite[Section 4.4.5]{htt}, which Lurie has rewritten since the book was published. References in the following proof are to theorem numbers from the current version of \cite{htt} freely available from Lurie's website at \url{https://www.math.ias.edu/~lurie/}.
\end{rmk}
\begin{proof}[Proof of \cref{lem:pointed-free}]
Following \cite[Definition 4.4.5.2]{htt}, let $\Idem^+$ denote the walking split idempotent, which is a 1-category, with an initial object $Y$ and another object $X$ of which $Y$ is a retract. We will use \cite[Proposition 4.4.5.6]{htt}, which says that as an $\infty$-category, $\Idem^+$ is freely generated by the objects $Y,X$, the two morphisms $i: Y \to X$ and $r: X \to Y$, and the homotopy $ri \simeq \id_Y$. Let us contemplate the $\infty$-category $\calP^\calK_\calR(\Idem^+)$ where $\calR$ comprises the cocone on the empty diagram with vertex at $Y$. Unraveling the definitions, $\calP(\Idem^+)$ is the $\infty$-category of spaces $\underline X$ equipped with a retract $\underline Y$, and $\calP_\calR(\Idem^+) \subset \calP(\Idem^+)$ is the full subcategory where $\underline Y$ is contractible. Thus we have a canonical equivalence $\calP_\calR(\Idem^+) \simeq \Top_\ast$. Under this identification, the localization $\calP(\Idem^+) \to \Top_\ast$ carries the representable at $X$ to $S^0$ and $Y$ to $0$. So $\calP^\calK_\calR(\Idem^+) \subset \Top_\ast$ is identified with $\Top_\ast^\calK$.

By the universal property of $\Top_\ast^\calK = \calP^\calK_\calR(\Idem^+)$, evaluation at $S^0$ determines an equivalence of categories $\Cat_\calK(\calP^\calK_\calR(\Idem^+),\calD) \to \Cat_{\{\emptyset\}}(\Idem^+,\calD)$ for any $\calD \in \Cat_\calK$. That is, $\calK$-cocontinuous functors $\calP^\calK_\calR(\Idem^+) \to \calD$ are identified with functors $\Idem^+ \to \calD$ preserving the initial object. By \cite[Proposition 4.4.5.6]{htt}, this means that evaluation at $X$ determines a functor $\Cat_{\{\emptyset\}}(\Idem^+,\calD) \to \calD^\sim$ whose fiber over $D \in \calD$ is the space whose points comprise the data of a morphism $0 \to D$, a morphism $D \to 0$, and a homotopy between the composite and $\id_0$ in $\calD(0,0)$ (here $0 \in \calD$ is the initial object). If $\calD$ is pointed, then this space is contractible, so the functor $\Cat_{\{\emptyset\}}(\Idem^+,\calD) \to \calD$ is an equivalence. It follows that the functor $\Cat_\calK(\calP^\calK_\calR(\Idem^+),\calD) \to \calD^\sim$, is an equivalence.

The proof is completed as soon as we note that $\Top_\ast^\calK$ is itself pointed.
\end{proof}


\begin{lem}\label{lem:pointed-zero}
Let $\calK$ be a collection of small $\infty$-categories containing the empty category. Then $\Cat_{\ast,\calK}$ is a pointed $\infty$-category; its zero object is the terminal category $[0]$.
\end{lem}
\begin{proof}
The unique functor $\calC \to [0]$ is as in $\Cat_\infty$. The unique functor $[0] \to \calC$ picks out the initial object. This functor has no automorphisms because the initial object of $\calC$ has no automorphisms. It is $\calK$-colimit-preserving.
\end{proof}

\begin{lem}\label{lem:pointed-mon}
Let $\calK$ be a collection of small $\infty$-categories containing the empty category. The full suboperad $\Cat_{\ast,\calK} \subset \Cat_\calK$ is a $\otimes$-ideal and an exponential ideal, and is unital. 
The unit is the $\infty$-category $\Top_\ast^\calK$ of \cref{lem:pointed-free}.
\end{lem}
\begin{proof}
For any $\calC \in \Cat_\calK$, the functors $\calC \otimes (-)$ and $\Fun_\calK(\calC,-)$ are 2-functorial under the natural enrichment of $\Cat_\calK$ in itself via its exponentially closed structure. Moreover, these functors preserve the zero object of $\Cat_\calK$, i.e. the terminal category $[0]$. An $\infty$-category $\calD$ is pointed if and only if the unique functor $\calD \to [0]$ has a left and right adjoints which are isomorphic. Moreover, if $\calD \in \Cat_\calK$, then all the functors involved preserve $\calK$-colimits. Since adjunctions are preserved by any 2-functor, it follows that $\calC \otimes (-)$ and $\Fun_\calK(\calC,-)$ preserve pointed $\infty$-categories, i.e. $\Cat_{\ast,\calK}$ is a $\otimes$-ideal and an exponential ideal.

That $\Cat_{\ast,\calK}$ is unital follows from 
\cref{lem:pointed-free}.
\end{proof}

\begin{cor}\label{cor:free-on-objs}
Let $\calK$ be a collection of small $\infty$-categories containing the empty category. The full suboperad $\Cat_{\ast,\calK} \subset \Cat_\calK$ is symmetric monoidal, and the inclusion functor is lax symmetric monoidal and preserves $\otimes$. The unit is the $\infty$-category $\Top_\ast^\calK$ of \cref{lem:pointed-free}.
\end{cor}
\begin{proof}
This follows from \cref{lem:pointed-mon} and \cref{lem:sub-mon}.
\end{proof}




\begin{rmk}[The infinitary perspective]\label{rmk:ind}
Let $\calK \subseteq \calK'$ be an inclusion of classes of small $\infty$-categories. By \cite[Remark 4.8.1.8]{ha}, the functor $\calP^{\calK'}_\calK: \Cat_\calK \to \Cat_{\calK'}$ is left adjoint to the forgetful functor, and moreover $\calP^{\calK'}_\calK$ is strong symmetric monoidal. As noted in \cite[Proposition 4.8.1.10 and Corollary 4.8.1.14]{ha}, it is in particular the case that if $\calC$ is a symmetric mononoidal $\infty$-category with compatible finite colimits, then $\Ind(\calC)$ is also symmetric monoidal under Day convolution, and has the universal property that $\Fun^{\otimes,L}(\Ind(\calC),\calD) = \Fun^\otimes(\calC,\calD)$ when $\calD$ is a symmetric monoidal presentable $\infty$-category. Moreover, if $\calK$ is a presentably symmetric monoidal $\infty$-category and $T \in \calK$ is dualizable with twisted-trivial braiding, then in the functorial splitting of \cref{prop:monoidal-splitting}, the factors $\calK_{T\wedge T^\dual}$ and $\calK/T$ are presentably symmetric monoidal and the localization functors onto them are symmetric monoidal left adjoints.
\end{rmk}

%
%

\section[Duals and Colimits]{Symmetric Monoidal $\infty$-Categories with Duals and Certain Colimits}\label{chap:dual-colim}

In this section, the kernel of ideas from \cref{chap:1-cat} and the infrastructure from \cref{chap:oo-cat} are combined to study the structure of certain $\infty$-categories of symmetric monoidal $\infty$-categories with duals and certain colimits. In \cref{sec:can-split} we discuss the fact that any symmetric monoidal $\infty$-category with duals and certain colimits admits a splitting as a product of certain subcategories characterized by various properties. For example, one factor is stable, another is an additive 1-category, etc. The splittings all arise from certain objects with twisted-trivial braidings guaranteed to exist by virtue of having the appropriate colimits. The splittings are also functorial in nature, so any sutiably cocontinuous, symmetric monoidal functor between such $\infty$-categories respects them, and thus the $\infty$-category of all such $\infty$-categories itself splits according to these factors. In \cref{sec:init}, we contribute toward a preliminary understanding of some of these factors, by computing the free symmetric monoidal $\infty$-category with duals and $\calK$-colimits for various $\calK$.

\subsection{Canonical splittings}\label{sec:can-split}

In this subsection, we identify two canonical objects with twisted-trivial braiding which are guaranteed to exist in any symmetric monoidal $\infty$-category which is suitably cocomplete. The first (\cref{subsec:gp-agp}) corepresents grouplike elements in the semiadditive setting with sufficient cofibers. The second (\cref{subsec:stab-costab}) is the suspension of the unit, which exists in the pointed case with suspensions. In both cases, the yoga of \cref{subsec:fun-split} tells us that if in addition these objects are dualizable, we obtain a splitting of the entire $\infty$-category. These splittings are characterized in the respective sections, and their interaction is discussed in \cref{subsec:3-split}.

\subsubsection{The grouplike / anti-grouplike splitting}\label{subsec:gp-agp}

In this subsection, we describe the object corepresenting grouplike elements of hom-spaces in a semiadditive symmetric monoidal $\infty$-category with cofibers. We discuss the splitting which results from \cref{subsec:fun-split} if this object is dualizable. 

\begin{lem}\label{lem:gp-space}
Let $(X,\mu,\eta)$ be a homotopy commutative, homotopy associative, $H$-space. Let $i: X^\gp \to X$ denote the inclusion of the grouplike part of $X$, i.e. the disjoint union of connected components of $X$ which have inverses under the multiplication $\mu$. Note that $X^\gp$ is a grouplike $H$-space under the restriciton of the map $\mu$. Let $-1: X^\gp \to X^\gp$ be a map sending each point to a homotopy inverse under $\mu$. Then the following commutative square is a homotopy pullback:
\begin{center}
    \begin{tikzcd}
        X^\gp \ar[r,"{(i,i(-1))}"] \ar[d] \ar[dr,"\lrcorner",very near start,phantom] & X \times X \ar[d,"\mu"] \\
        \ast \ar[r,"\eta"] & X
    \end{tikzcd}
\end{center}
\end{lem}
\begin{proof}
We show this in several steps.

\textbf{Step 1:} First note that by \cref{rmk:cogp-cof-of-diag}, the lemma is true when $X$ is discrete.

Now let $F$ be the fiber of $\mu$, so that we have a natural map $X^\gp \to F$. 

\textbf{Step 2:} Because $\mu$ is split by $\id_X \times \eta$ (or by $\eta \times \id_X)$, we have short exact sequences $\pi_n(F) \to \pi_n(X \times X) \to \pi_n(X)$ for each $n \in \nats$. When $n=0$, this tells that $\pi_0(F) = \pi_0(X)^\gp$ by Step 1. Since $\pi_0(X)^\gp = \pi_0(X)^\gp$, this shows that the comparison map $\pi_0(X^\gp) \to \pi_0(F)$ is a bijection.

\textbf{Step 3:} For $n \geq 1$, we identify the short exact sequence $\pi_n(F) \to \pi_n(X \times X) \to \pi_n(X)$ with the short exact sequence $\pi_0(\Omega^n F) \to \pi_0(\Omega^n X \times \Omega^n X) \to \pi_0(\Omega^n X)$. By Step 2, the map $\pi_0((\Omega^n X)^\gp) \to \pi_0(\Omega^n F)$ is a bijection. Since $\Omega^n (X^\gp) = \Omega^n X = (\Omega^n X)^\gp$, we have that $\pi_n(X^\gp) = \pi_0(\Omega^n (X^\gp)) \to \pi_0(\Omega^n F) = \pi_n(F)$ is a bijection. So by Whitehead's theorem, $X^\gp \to F$ restricts to an equivalence on the connected component of the identity.

\textbf{Step 4:} $F$ inherits the structure of a commutative, associative $H$-space as the fiber of a map of such, and the map $X^\gp \to F$ is a map of $H$-spaces. Because $X^\gp$ is grouplike and $X^\gp \to F$ is $\pi_0$-surjective, $F$ is also grouplike. Therefore, because $X^\gp \to F$ is an equivalence at the connected component of the identity (Step 3) and a $\pi_0$-bijection (Step 1), it follows that $X^\gp \to F$ is an equivalence.
\end{proof}

\begin{Def}
Let $\calC$ be a semiadditive $\infty$-category. Denote by $\calC_\gp \subseteq\calC$ the full subcategory of those objects $C$ such that the $E_\infty$-space $\calC(C,D)$ is grouplike for all $D \in \calC$. If $\calC = \calC_\gp$, we say that $\calC$ is \defterm{additive}. If $\calC$ has cofibers, then for $C \in \calC$, let $C_\gp$ be the cofiber of the diagonal $\Delta: C \to C \oplus C$.

Denote by $\calC_\notgp \subseteq \calC$ the full subcategory of objects $C$ such that $\calC(C,D)$ has no nonzero grouplike elements, for all $D \in \calC$. If $\calC = \calC_\notgp$, we say that $\calC$ is \defterm{anti-additive}.
\end{Def}

\begin{cor}\label{cor:gp-semiadd}
Let $\calC$ be a semiadditive $\infty$-category with cofibers, and let $C \in \calC$. Then the composite map $r_C : C \xrightarrow{(\id_C,0)} C \oplus C \to C_\gp$ induces a map factors canonically through the inclusion $\calC(C,D)^\gp \to \calC(C,D)$, and the induced map $r_C^\ast: \calC(C_\gp,D) \to \calC(C,D)^\gp$ is an equivalence.
\end{cor}
\begin{proof}
The cofiber sequence $C \xrightarrow{\Delta} C\oplus C \to C_\gp$ induces a fiber sequence $\calC(C,D) \overset \mu \leftarrow \calC(C,D) \times \calC(C,D) \leftarrow \calC(C_\gp,D)$ where $\mu$ is the addition on $\calC(C,D)$. So this follows from \cref{lem:gp-space}.
\end{proof}

\begin{thm}\label{thm:gp-mon}
Let $(\calC,\wedge,S)$ be an $E_n$-monoidal $\infty$-category with biproducts and cofibers over which $\wedge$ distributes. Then the composite map $r = r_S: S \xrightarrow{(\id_S,0)} S \oplus S \to S_\gp$ of \cref{cor:gp-semiadd} is a closed idempotent. The $S_\gp$-stable objects are exactly the objects of $\calC_\gp$, and the $S_\gp$-torsion objects are exactly the objects of $\calC_\notgp$.
\end{thm}
\begin{proof}
Note that for any $C \in \calC$, the cofiber sequence $C \to C\oplus C \to C_\gp$ may be identified with the cofiber sequence $C \wedge S \to C \wedge (S \oplus S) \to C \wedge S_\gp$. By \cref{cor:gp-semiadd}, we have that $\calC(S_\gp \wedge C, D) = \calC(C,D)^\gp$. In particular, $\calC(S_\gp \wedge S_\gp, D) = \calC(S_\gp,D)^\gp = (\calC(S,D)^\gp)^\gp = \calC(S,D)^\gp = \calC(S_\gp, D)$. So by the Yoneda lemma, $S_\gp \to S_\gp \wedge S_\gp$ is an isomorphism, i.e. $S \to S_\gp$ is a closed idempotent. From the equivalence $\calC(S_\gp \wedge C, D) = \calC(C,D)^\gp$ of \cref{cor:gp-semiadd} also follows the description of $S_\gp$-stable and $S_\gp$ torsion objects.
\end{proof}


\begin{cor}\label{cor:gp-dual}
Let $(\calC, \wedge ,S)$ be an $E_n$-monoidal $\infty$-category with biproducts and cofibers over which $\wedge$ distributes. Suppose that the cofiber $S_\gp$ of the diagonal $\Delta: S \to S \oplus S$ is dualizable. Then $S_\gp$ is a clopen idempotent, with complement $S_\notgp$. The induced functor $\calC \to \calC_\gp \times \calC_\notgp$ is an equivalence.
\end{cor}
\begin{proof}
By \cref{thm:gp-mon} and \cref{prop:clopens}(\ref{prop:clopens.item:equiv}), $S_\gp$ is a clopen idempotent. Because $S_\gp$ is a cogroup object (\cref{thm:gp-mon}), it follows from \cref{cor:split-cof}, that the cofiber $S/E$ of $E \to S$ is a complementary clopen idempotent. So the equivalence $\calC \to \calC_\gp \times \calC_\notgp$ follows from \cref{thm:gp-mon} and \cref{prop:monoidal-localization}(\ref{prop:monoidal-localization.item:split}). 
\end{proof}

\begin{cor}\label{cor:gp-loc}
Let $(\calC, \wedge ,S)$ be an $E_n$-monoidal $\infty$-category $(n \geq 2)$ with biproducts and cofibers over which $\wedge$ distributes. Suppose that the cofiber $S_\gp$ of the diagonal $\Delta: S \to S \oplus S$ is dualizable. Then the localization $\calC \to \calC_\gp$ is the universal $E_n$-monoidal, biproduct-and-cofiber-preserving functor to an additive $E_n$-monoidal $\infty$-category with biproducts and cofibers over which $\wedge$ distributes. The localization $\calC \to \calC_\notgp$ is the universal $E_n$-monoidal, biproduct-and-cofiber-preserving functor to an anti-additive $E_n$-monoidal $\infty$-category with biproducts and cofibers over which $\wedge$ distributes. 
\end{cor}
\begin{proof}
This follows from \cref{cor:gp-dual} and \cref{prop:monoidal-splitting}.
\end{proof}

\subsubsection{The stable / trivial suspension splitting}\label{subsec:stab-costab}

In this subsection, we discuss the suspension of the unit object. When this object is dualizable, it induces (via \cref{subsec:fun-split}) a splitting of the $\infty$-category into a stable part and a part with a dual characterization: all objects have trivial suspension.

\begin{Def}
Let $\calC$ be a pointed $\infty$-category with suspension $\Sigma$ (that is, for every $C \in \calC$ the object $\Sigma C = 0 \cup_C 0$ exists). We say that $\calC$ is \defterm{weakly stable} if the suspension functor $\Sigma: \calC \to \calC$ is an equivalence of categories, and \defterm{stable} if in addition $\calC$ has finite colimits (equivalently, $\calC$ has finite limits). We say that $\calC$ \defterm{has trivial suspension} if the suspension functor $\Sigma$ is constant at $0$.
\end{Def}

\begin{rmk}\label{rmk:stab-cogp}
Recall from \cref{eg:cogp-susp} that because $S^1 \in \Top_\ast$ is a cogroup object, it follows that any suspension object, in an $\infty$-category with suspensions and finite coproducts, is a cogroup object. Therefore, if $\calC$ is a weakly stable semiadditive $\infty$-category, then $\calC$ is additive.
\end{rmk}
\begin{rmk}
In an additive $\infty$-category, the coequalizer of two maps $f,g$ may be computed as the cofiber of $f-g$. It follows that an additive $\infty$-category has finite colimits iff it has cofibers. A weakly stable $\infty$-category is stable if and only if it has finite coproducts and cofibers, if and only if it has finite colimits.
\end{rmk}

\begin{prop}\label{prop:stab-dual}
Let $(\calC,\wedge,S)$ be an $E_n$-monoidal $\infty$-category $(n \geq 2)$ with finite coproducts and cofibers over which $\wedge$ distributes. Suppose that the suspension of the unit $\Sigma S$ has a dual $\Sigma\inv S$. Then $S_\stab := \Sigma S \wedge \Sigma\inv S$ is a clopen idempotent, with a complement $S_\trivsusp$.
\end{prop}
\begin{proof}
Because $\Sigma S$ has twisted-trivial braiding (\cref{eg:twist-susp-unit}), $S_\stab$ is a clopen idempotent by \cref{prop:clopens}(\ref{prop:clopens.item:untwist}). Because $S_\stab$ is a suspension object, it is a cogroup object (\cref{rmk:stab-cogp}). So by \cref{cor:split-cof}, $S_\trivsusp$ is a complementary clopen idempotent.
\end{proof}

\begin{Def}
Let $(\calC,\wedge,S)$ be an $E_n$-monoidal $\infty$-category $(n \geq 2)$ with finite coproducts and cofibers over which $\wedge$ distributes. We write $\calC_\stab = \calC_{S_\stab}$ and $\calC_\trivsusp = \calC_{S_\trivsusp}$.
\end{Def}

\begin{thm}
Let $(\calC,\wedge,S)$ be an $E_n$-monoidal $\infty$-category $(n \geq 2)$ with finite coproducts and cofibers over which $\wedge$ distributes. Suppose that the suspension of the unit $\Sigma S$ has a dual $\Sigma\inv S$. Then the the canonical functor $\calC \to \calC_\stab \times \calC_\trivsusp$ is an equivalence. The localization $\calC \to \calC_\stab$ is the universal $E_n$-monoidal, finite-coproduct-and-cofiber-preserving functor to an $E_n$-monoidal stable $\infty$-category. The localization $\calC \to \calC_\trivsusp$ is the universal $E_n$-monoidal, finite-coproduct-and-cofiber-preserving functor to an $E_n$-monoidal $\infty$-category with trivial suspension and finite coproducts and cofibers over which $\wedge$ distributes.
\end{thm}
\begin{proof}
The equivalence $\calC \to \calC_\stab \times \calC_\trivsusp$ follows from \cref{prop:stab-dual} and \cref{prop:monoidal-localization}(\ref{prop:monoidal-localization.item:split}). By \cref{rmk:stab-cogp}, $\calC_\stab$ is stable. Any $E_n$-monoidal, finite-coproduct-and-cofiber-preserving functor $F$ to an $E_n$-monoidal stable $\infty$-category $(\calD,\wedge S)$ carries $\Sigma S$ to the $\wedge$-invertible object $\Sigma S$ (which has monoidal inverse $\Sigma \inv S$). Therefore by \cref{prop:monoidal-splitting}(\ref{prop:monoidal-splitting.item:1}), $F$ factors uniquely through $\calC_\stab$. Because $\calC \to \calC_\stab$ is a product projection, the induced functor $\calC_\stab \to \calD$ preserves any colimits which $F$ does. This establishes the universal property of $\calC_\stab$. Similarly, any $E_n$-monoidal, finite-coproduct-and-cofiber-preserving functor $F$ to an $E_n$-monoidal $\infty$-category $(\calD,\wedge S)$ with trivial suspension and finite coproducts and cofibers over which $\wedge$ distributes, carries $\Sigma S$ to $0$. Therefore by \cref{prop:monoidal-splitting}(\ref{prop:monoidal-splitting.item:3}), $F$ factors uniquely through $\calC_\trivsusp$, and as before the preservation of the appropriate colimits is automatic. This establishes the universal property of $\calC_\trivsusp$.
\end{proof}

\subsubsection{The 3-fold splitting}\label{subsec:3-split}

In this subsection, we describe the interaction between the splittings of the previous two subsections, yielding a 3-fold splitting of any symmetric monoidal $\infty$-category with duals, finite coproducts, and cofibers.

\begin{prop}\label{prop:3-split}
Let $(\calC, \wedge, S)$ be an $E_n$-monoidal $\infty$-category $(n \geq 2)$ with finite biproducts and cofibers over which $\wedge$ distributes. Suppose that $\Sigma S$ and the cofiber $S_\gp$ of the diagonal $S \to S \oplus S$ have duals. Then the $E_n$-monoidal localization functors of \cref{cor:gp-dual} and \cref{prop:stab-dual} assemble into an equivalence of categories 

\[\calC \to \calC_{\stab} \times \calC_{\gptrivsusp} \times \calC_{\notgp}\]

Here the $E_n$-localization functor $\calC \to \calC_\gptrivsusp = \calC_\gp \cap \calC_\trivsusp$ is the universal $E_n$-monoidal functor to an $E_n$-monoidal additive $\infty$-category with trivial suspension and cofibers over which $\wedge$ distributes.
\end{prop}
\begin{proof}
As in \cref{rmk:stab-cogp}, any weakly stable $\infty$-category is additive. Therefore, the clopen idempotent $S_\stab$ refines the clopen idempotent $S_\gp$, and dually the clopen idempotent $S_\notgp$ refines the clopen idempotent $S_\trivsusp$. Moreover, because $n \geq 2$, the two closed idempotent structures $S \to S_\gp \to S_\gp \wedge S_\trivsusp$ and $S \to S_\trivsusp \to S_\trivsusp \wedge S_\gp$ agree, i.e. the localizations $\calC \to \calC_\gp$ and $\calC \to \calC_\trivsusp$ commute. The result follows.
\end{proof}

We have already explored some of the properties of $\calC_\stab$ and $\calC_\notgp$. The following property of $\calC_\gptrivsusp$ may come as a surprise:

\begin{prop}\label{prop:add-1-cat}
Let $(\calC, \wedge, S)$ be an $E_n$-monoidal $\infty$-category $(n \geq 2)$ with finite biproducts and cofibers over which $\wedge$ distributes. Suppose that $\Sigma S$ and $S_\gp$ have duals. Then the localization functor $\calC \to \calC_\gptrivsusp$ is the universal $E_n$-monoidal, finite-coproduct-and-cofiber-preserving functor to an $E_n$-monoidal additive 1-category with cofibers over which $\wedge$ distributes.
\end{prop}
\begin{proof}
We must first show that $\calC_\gptrivsusp$ is a 1-category, i.e. has discrete hom-spaces. For any $C,D \in \calC_\gptrivsusp$, we have $\pi_n \calC_\gptrivsusp(C,D) = \pi_0 \calC_\gptrivsusp(\Sigma^n C, D)$, where the basepoint is the zero morphism. Because suspension is trivial it follows that the component of $\calC_\gptrivsusp(C,D)$ at the zero morphism is contractible. Moreover, $\calC_\gptrivsusp(C,D)$ is a grouplike $H$-space because $\calC_\gptrivsusp$ is additive. So all of its connected components are homotopy equivalent, and so $\calC_\gptrivsusp(C,D)$ is discrete.

It is clear that $\calC_\gptrivsusp$ is additive and has cofibers over which $\wedge$ distributes. Conversely, if $\calD$ is an $E_n$-monoidal additive 1-category with cofibers over which $\wedge$ distributes, then $\calD$ has trivial suspension. So if $F: \calC \to \calD$ is an $E_n$-monoidal, finite-coproduct-and-cofiber-preserving functor, it carries $\Sigma S$ to $0$ and carries $S \to S_\gp$ to an isomorphism, and therefore factors uniquely through $\calC_\gptrivsusp$ by \cref{prop:monoidal-splitting}. It automatically preserves any colimits which $F$ does.
\end{proof}

Additive 1-categories with appropriate duals admit some further splittings:

\begin{prop}\label{prop:prime-closed}
Let $(\calC,\wedge,S)$ be an $E_n$-monoidal semiadditive 1-category with cofibers over which $\wedge$ distributes. Let $m \in \nats$ be a number. Then the cofiber $S/m$ of the map $m: S \to S$ (i.e. the $m$-fold sum of the identity with itself) is a closed idempotent. We have $\calC(C\wedge S/m,D) = \calC(C,D)[m]$ (where $A[m]$ denotes the $m$-torsion subgroup of $A$) for all $C, D \in \calC$, so that the $S/m$-local objects are those objects $C \in \calC$ such that $m=0$ in $\calC(C,D)$ for all $D \in \calC$, and the $S/m$-torsion objects are those $C \in \calC$ such that $\calC(C,D)[m] = 0$ for all $D \in \calC$.
\end{prop}

\begin{warning}
\cref{prop:prime-closed} is false in a general $\infty$-category. even if $m$ is prime. For instance, in spectra, $S/m$ is the mod-$m$ Moore spectrum, and $S/m \not \simeq S/m \wedge S/m$.
\end{warning}

\begin{proof}[Proof of \cref{prop:prime-closed}]
The cofiber sequence $C \xrightarrow m C \to C \wedge S/m$ induces a fiber sequence $\calC(C,D) \overset{m}{\leftarrow} \calC(C,D) \leftarrow \calC(C\wedge S/m, D)$ for any $D \in \calD$. Because $\calC$ is a 1-category, this says exactly that $\calC(C \wedge S/m,D) = \calC(C,D)[m]$. So $\calC(S/m\wedge S/m, D) = \calC(S/m,D)[m] = (\calC(S,D)[m])[m] = \calC(S,D)[m] = \calC(S/m, D)$; by Yoneda $S/m$ is a closed idempotent, and the rest follows.
\end{proof}

\begin{Def}
Let $\calC$ be an additive 1-category and $m \in \ints$. We say that $C \in \calC$ is \defterm{of characteristic dividing $m$} if the endomorphism $m: C \to C$ (i.e. the $m$-fold sum of the identity) is the zero morphism. We write $\calC/m \subseteq \calC$ for the full subcategory of objects of characteristic dividing $m$. We say that $\calC$ \defterm{has characteristic dividing $m$} if $\calC = \calC/m$, i.e. if all hom-groups have characteristic dividing $m$.

We say that $C$ is \defterm{$m$-torsion-free} if $m: C \to C$ is a monomorphism and \defterm{co-$m$-torsion-free} if $m: C \to C$ is an epimorphism. We write $\calC{(m)} \subseteq \calC$ for the full subcategory of objects which are co-$m$-torsion-free. We say that $\calC$ is \defterm{$m$-torsion-free} if $\calC = \calC{(m)}$, i.e. if all hom-groups are $m$-torsion-free.
\end{Def}

\begin{rmk}
The notation $\calC{(m)}$ is perhaps misleading: $m$ need not act invertibly on $\calC{(m)}$, but only without torsion.
\end{rmk}

\begin{rmk}
The reader may have expected the co-$m$-torsion-free objects to be called ``$m$-divisible" by analogy to the category of abelian groups. We have opted against such terminology, finding that it is unhelpful in the present setting. For us, the important fact about co-$m$-torsion-free objects $C$ is that $\Hom(C,D)$ is always an $m$-torsion-free abelian group, so it seems best to have the term ``torsion-free" appear in the name. In fact, we will have little use for $m$-torsion-free objects, and have only introduced them to stand as foils for the co-$m$-torsion-free objects.
\end{rmk}

\begin{cor}\label{cor:prime-split}
Let $(\calC,\wedge,S)$ be an $E_n$-monoidal semiadditive 1-category with cofibers over which $\wedge$ distributes. Let $m \in \nats$ be a number, and suppose that $S/m$ has a dual. Then $S/m$ is a clopen idempotent, with complement $S{(m)}$. The induced functor $\calC \to \calC/m \times \calC{(m)}$ is an equivalence. The localization $\calC \to \calC/m$ is the universal $E_n$-monoidal, semiadditive, cofiber-preserving functor to an additive 1-category of characteristic dividing $m$ with cofibers over which $\wedge$ distributes.
The localization $\calC \to \calC{(m)}$ is the universal $E_n$-monoidal, semiadditive, cofiber-preserving functor to an additive 1-category with cofibers over which $\wedge$ distributes which is $m$-torsion-free.
\end{cor}

\begin{proof}[Proof of \cref{cor:prime-split}]
That $S/m$ is clopen follows from \cref{prop:clopens}(\ref{prop:clopens.item:equiv}). Because $\calC$ is additive, $S/m$ is a cogroup object, so has a complement by \cref{cor:split-cof}, and \cref{prop:ordinary-split}. So the equivalence $\calC \to \calC/m \times \calC{(m)}$ follows by \cref{prop:monoidal-splitting}.

We have seen that $\calC/m$ has characteristic $m$. If $F: \calC \to \calD$ is an $E_n$-monoidal, semiadditive, cofiber-preserving functor to an additive 1-category of characteristic dividing $m$ with cofibers over which $\wedge$ distributes, then $F$ carries $S \to S/m$ to an isomorphism. So $F$ factors uniquely through $\calC/m$ by \cref{prop:monoidal-splitting}, with colimit preservation being automatic. Likewise, we have seen that $\calC{(m)}$ is $m$-torsion-free. If $F: \calC \to \calD$ is an $E_n$-monoidal, semiadditive, cofiber-preserving functor to an additive 1-category with cofibers over which $\wedge$ distributes which is $m$-torsion-free, then $F$ carries $S/m$ to $0$. So $F$ factors through $\calC{(m)}$ by \cref{prop:monoidal-splitting}, with colimit preservation being automatic.
\end{proof}

\begin{rmk}
In the setting of \cref{cor:prime-split}, the localizations $\calC \to \calC/m$ for varying $m$ are compatible with one another whenever they exist. In particular, if $m = p_1^{e_1} \cdots p_r^{e_r}$ is the prime factorization, then $\calC_m = \calC/p_1 \times \cdots \times \calC/p_r$ so long as $S/p_1,\dots, S/p_r$ all have duals.
\end{rmk}

\subsection{Initial objects}\label{sec:init}

In this section, we begin a preliminary study of symmetric monoidal $\infty$-categories with duals and various finite colimits suggested by the taxonomy of \cref{subsec:3-split}. The goal, achieved in \cref{subsec:all}, is to compute the initial object of the $\infty$-category of symmetric monoidal $\infty$-categories with duals, finite coproducts, and cofibers. This is achieved by considering each of the factors from \cref{subsec:3-split} in turn.

\begin{Def}
Let $\SMC_{\amalg,\cof}$ denote the $\infty$-category of symmetric monoidal $\infty$-categories with finite coproducts and cofibers over which $\wedge$ distributes. Let
\[\SMD_{\amalg,\cof} = \SMD \times_\SMC \SMC_{\amalg,\cof}\]
denote the full subcategory of symmetric monoidal $\infty$-categories with finite coproducts, cofibers, and duals. Let $\SMD_\stab \subset \SMD_{\amalg,\cof}$ denote the full subcategory of symmetric monoidal stable $\infty$-categories with duals. Let $\SMD_{\gptrivsusp,\cof} \subset \SMD_{\amalg,\cof}$ denote the full subcategory of symmetric monoidal additive 1-categories with duals and cofibers. For each prime $p$, let $\SMD_{p,\trivsusp,\cof} \subset \SMD_{\amalg,\cof}$ denote the full subcategory of symmetric monoidal additive 1-categories with duals and cofibers of characteristic dividing $p$. Let $\SMD_{\notgp,\amalg,\cof} \subset \SMD_{\amalg,\cof}$ denote the full subcategory of symmetric monoidal anti-additive $\infty$-categories with duals and cofibers.
\end{Def}

\subsubsection{The stable case}\label{subsec:stab}

In this subsection, we compute the initial object of the $\infty$-category of symmetric monoidal $\infty$-categories with duals, finite coproducts, and cofibers which are moreover \emph{stable}. Unsurprisingly, it turns out to be the $\infty$-category of finite spectra, symmetric monoidal under smash product.

\begin{Def}
Let $\Cat_\rex$ be the symmetric monoidal $\infty$-category $\Cat_\calK$ (see \cref{def:cat-k}) when $\calK$ is the collection of finite categories. That is, $\Cat_\rex$ is the symmetric monoidal $\infty$-category of $\infty$-categories with finite colimits. Let $\Cat_\stab \subset \Cat_\rex$ be the full suboperad of $\Cat_\rex$ whose objects are the stable $\infty$-categories.
\end{Def}

\begin{lem}\label{lem:stable-free}
Let $\calD \in \Cat_\stab$ be a stable $\infty$-category. Then evaluation at the sphere $\bbS \in \Spt^\fin$ determines an equivalence of categories $\Cat_\stab(\Spt^\fin, \calD) \to \calD^\sim$.
\end{lem}
\begin{proof}
By \cref{cor:free-on-objs}, $\Top_\ast^\fin$ is the free pointed, right exact $\infty$-category on an object. The Spanier-Whitehead construction tells us that $\Spt^\fin = \varinjlim(\Top_\ast^\fin \xrightarrow \Sigma \Top_\ast^\fin \xrightarrow \Sigma \dots)$. Let $\calD$ be a stable $\infty$-category, and $F: \Top_\ast \to \calD$ the right-exact functor classifying an object $D \in \calD$. Because $F$ commutes with suspension, which is invertible on $\calD$, it easily follows that the space of extensions of $F$ along $\Sigma^\infty: \Top_\ast^\fin \to \Spt^\fin$ is contractible. As $\Spt^\fin$ is stable, the result follows.
\end{proof}

\begin{thm}\label{thm:stab-mon}
The full suboperad $\Cat_\stab \subset \Cat_\rex$ is a $\otimes$-ideal and an exponential ideal, and unital with unit the $\infty$-category $\Spt^\fin$ of finite spectra.
\end{thm}
\begin{proof}
The proof will follow the pattern of the proof of \cref{lem:pointed-mon}.
We first show that $\Cat_\stab \subset \Cat_\rex$ is a $\otimes$-ideal. Let $\calC \in \Cat_\stab$ and $\calD \in \Cat_\rex$; we will show that the tensor product $\calC \otimes \calD$ in $\Cat_\rex$ is stable. Because the suspension functor $\Sigma: \calC \to \calC$ is an equivalence, it follows by functoriality of $(-) \otimes \calD$ that $\Sigma \otimes \calD: \calC \otimes \calD \to \calC \otimes \calD$ is an equivalence. But this is nothing other than the suspension functor on $\calC \otimes \calD$. This follows from the fact that $0 \otimes \calD = 0$, where $0$ is a constant functor at the initial object, and the fact that $(F\cup_G H) \otimes \calD = (F \otimes \calD) \cup_{G \otimes \calD} (H \otimes \calD)$ for any functors $F,G,H$.

Unitality follows from \cref{lem:stable-free}.
\end{proof}

\begin{cor}\label{cor:stab-mon}
The full suboperad $\Cat_\stab \subset \Cat_\rex$ is symmetric monoidal, with unit $\Spt^\fin$. The inclusion functor is lax symmetric monoidal, preserving $\otimes$.
\end{cor}
\begin{proof}
This follows from \cref{thm:stab-mon} and \cref{lem:sub-mon}.
\end{proof}

\begin{prop}\label{prop:init-stab}
The initial object of $\SMC_\stab$ is the symmetric monoidal $\infty$-category $\Spt^\fin$ of finite spectra, which is also the initial object of $\SMD_\stab$.
\end{prop}
\begin{proof}
By \cref{lem:comm-init}, the first statement follows from \cref{thm:stab-mon}. So for the second statement, it suffices to observe that $\Spt^\fin$ has duals for all objects. This is guaranteed by the theory of Spanier-Whitehead duality \cite{dold-puppe}, which reduces to the observation that $\bbS \in \Spt^\fin$ is (tautologically) dualizable, and that $\Spt^\fin$ is stable and generated under finite colimits and desuspensions by $\bbS$, so that this follows from \cref{lem:dual-stable}.
\end{proof}

\subsubsection{The anti-additive case}\label{subsec:agp}

In this subsection, we compute the initial object in the $\infty$-category of symmetric monoidal $\infty$-categories with duals, finite coproducts, and cofibers which are moreover \emph{anti-additive} in the sense of \cref{subsec:gp-agp}. It turns out to be the $(2,1)$-category of spans of finite sets, symmetric monoidal under cartesian product.

\begin{Def}
Let $\Cat_\amalg$ be the symmetric monoidal $\infty$-category $\Cat_\calK$ (see \cref{def:cat-k}) when $\calK$ is the collection of finite discrete categories. That is, $\Cat_\amalg$ is the symmetric monoidal $\infty$-category of $\infty$-categories with finite coproducts. Let $\Cat_\oplus \subset \Cat_\amalg$ be the full subcategory of $\Cat_\amalg$ whose objects are the semiadditive $\infty$-categories.
\end{Def}

\begin{lem}
Let $\calK$ be a collection of small $\infty$-categories containing the finite discrete categories. Then $\Cat_\calK$ is semiadditive; its biproducts are given by cartesian product.
\end{lem}
\begin{proof}
By \cref{lem:pointed-zero}, $\Cat_\calK$ is pointed by the terminal category $[0]$. Assume for the moment that $\Cat_\calK$ is small. Then by abstract nonsense, $\Cat_\calK$ has finite coproducts. So by \cref{cor:pt-semiadd-hocat}, it will suffice to verify that the homotopy category $\ho\Cat_\calK$ is semiadditive. But this is clear: every object $\calC \in \ho\Cat_\calK$ has a commutative monoid structure given by finite coproducts, and every morphism is a morphism of monoids because the morphisms of $\Cat_\calK$ preserve finite coproducts.

If $\calK$ is not small, then write it as a directed union $\calK = \cup_i \calK_i$ where each $\calK_i$ is small. By the foregoing, each $\Cat_{\calK_i}$ is semiadditive, and the inclusions $\Cat_{\calK_i} \to \Cat_{\calK_j}$ preserve the semiadditive structure. Thus $\Cat_\calK$ is also semiadditive. 
\end{proof}

\begin{thm}\label{thm:semiadd-mon}
The full suboperad $\Cat_\oplus \subset \Cat_\amalg$ is a $\otimes$-ideal and an exponenetial ideal. It is also unital, with unit $\Span(\Fin)$.
\end{thm}
\begin{proof}
The proof will follow the pattern of the proof of \cref{lem:pointed-mon} and \cref{thm:stab-mon}. If $\calC \in \Cat_\amalg$, then $\calC \otimes_\amalg (-)$ and $\Fun_\amalg(\calC,-)$ are 2-functorial and preserve zero objects and biproducts. If $\calD$ is semiadditive, then the diagonal $\calD \to \calD \times \calD$ has a left and right adjoint which coincide, and these functors live in $\Cat_\amalg$. The adjunctions are preserved by $\calC \otimes_\amalg (-)$ and $\Fun_\amalg(\calC,-)$, and therfore $\calC \otimes_\amalg \calD$ and $\Fun_\amalg(\calC,\calD)$ are semiadditive. Thus $\Cat_\oplus$ is a $\otimes$-ideal and an exponential ideal.

Unitality follows from \cite[Theorem A.1]{glasman}.
\end{proof}

\begin{cor}\label{cor:semiadd-mon}
The full suboperad $\Cat_\oplus \subset \Cat_\amalg$ is symmetric monoidal with unit $\Span(\Fin)$. The inclusion functor is lax symmetric monoidal and preserves $\otimes$.
\end{cor}
\begin{proof}
This follows from \cref{thm:semiadd-mon} and \cref{lem:sub-mon}.
\end{proof}

\begin{cor}\label{cor:init-semiadd}
The symmetric monoidal $(2,1)$-category $\Span(\Fin)$ is the initial object in the $\infty$-category of symmetric monoidal semiadditive $\infty$-categories.
\end{cor}
\begin{proof}
By \cref{lem:comm-init}, this follows from \cref{thm:semiadd-mon}.
\end{proof}

\begin{cor}\label{cor:init-antiadd}
The symmetric monoidal $(2,1)$-category $\Span(\Fin)$ is the initial object in $\SMD_{\notgp,\amalg,\cof}$.
\end{cor}
\begin{proof}
By \cref{cor:init-semiadd}, it will suffice to verify that $\Span(\Fin)$ has duals and cofibers, is anti-additive, and that the cofibers are preserved by any symmetric monoidal functor to an anti-additive symmetric monoidal $\infty$-category with duals and cofibers. Indeed, coproducts in $\Span(\Fin)$ are disjoint union. Since every object is a finite disjoint union of copies of the unit, it follows from \cref{lem:dual-stable} that every object is dualizable (and indeed, self-dual). Since sums in hom-spaces are given by disjoint union, $\Span(\Fin)$ is clearly anti-additive, and in particular has trivial suspension.

Let us compute cofibers in $\Span(\Fin)$. Any morphism in $\Span(\Fin)$ factors as a backwards morphism in $\Fin$ followed by a forward morphism in $\Fin$. We may further decompose the backward arrow as a coproduct of a bijection with maps $1 \leftarrow 0$ followed by a coproduct of a bijection with maps $1 \leftarrow n$ with $n \geq 2$. The cofiber of the former are suspensions and hence zero. The cofiber of the latter likewise vanish by antiadditivity (vanishing of the cofiber of $1 \leftarrow 2$ is the \emph{definition} of antiadditivity). So the cofiber of a backward map is zero, and this cofiber is preserved by any semiadditive functor to an antiadditive $\infty$-category. A forward map may be factored as a coproduct of a bijection with maps $n \to 1$ followed by a coproduct of a bijection with maps $0 \to 1$. The cofiber of a map $0 \to 1$ is $1$, and is preserved by all functors. It follows that the cofiber of an inclusion map $i_0: 1 \to 2$ is $1$, preserved by all semiadditive functors (the relevant morphism $2 \to 1$ is the backward map $i_1: 2 \leftarrow 1$). Since $i_0$ is a section of the map $2 \to 1$, and the cofiber of the identity is $0$, it follows by pasting properties of pushout squares that the pushout of $i_1: 2 \leftarrow 1$ along the forward map $2 \to 1$ is $0$, and this pushout is preserved by all semiadditive functors. Then because the cofiber of $1 \to 0$ is $0$ by the triviality of suspension, it follows by pasting pushout squares that the cofiber of $2 \to 1$ is $0$, and this cofiber is preserved by any semiadditive functor to a semiadditive $\infty$-category with trivial suspension. Since we have already computed the cofiber of the map $0 \to 1$ to be $0$ and seen that this cofiber is preserved by all semiadditive functors, it follows that the cofiber of any morphism of $\Span(\Fin)$ exists and is preserved by any semiadditive functor to an antiadditive $\infty$-category, as desired.
\end{proof}

\subsubsection{The additive 1-category case}

In this subsection, we compute the initial object in the $\infty$-category of symmetric monoidal $\infty$-categories with duals, finite coproducts, and cofibers which are moreover \emph{additive}. It turns out to be a curious sort of restricted product of categories of finite-dimensional vector spaces over prime fields.

\begin{Def}
Let $\Cat_{\rex,\oplus} \subset \Cat_\rex$ be the full sub-operad of semiadditive $\infty$-categories with finite colimits. Let $\Cat_{\rex,\add} \subset \Cat_{\rex,\add}$ be the full suboperad of additive $\infty$-categories with finite colimits. Let $\Cat_{\rex,\add,\trivsusp} \subset \Cat_{\rex,\add}$ be the full subcategory of additive $\infty$-categories with finite colimits and trivial suspension (i.e, by \cref{prop:add-1-cat}, additive 1-categories with finite colimits).
\end{Def}

\begin{prop}\label{prop:semiadd-mon-rex}
The sub-operad $\Cat_{\rex,\oplus} \subset \Cat_\rex$ is a $\otimes$-ideal.
\end{prop}
\begin{proof}
The proof is the same as \cref{thm:semiadd-mon}.
\end{proof}

\begin{prop}
The sub-operad $\Cat_{\rex,\add} \subset \Cat_{\rex,\oplus}$ is a $\otimes$-ideal.
\end{prop}
\begin{proof}
Let $\calC \in \Cat_{\rex,\add}$ and $\calD \in \Cat_{\rex,\oplus}$. Then $\calC \otimes \calD = \calP^\rex_{\rex \boxtimes \rex}(\calC \times \calD)$. The localization of a representable $L(C,D)$ is a cogroup object because $C$ is a cogroup object and $L(0,D) = 0, L(C \oplus C,D) = L(C,D) \oplus L(C,D)$. Cogroup objects are closed under direct sums, so it remains to check that cogroup objects are closed under cofibers. This is true because being an object $X$ in a semiadditive category is a cogroup object if and only if the map $\begin{pmatrix} \id_X & \id_X \\ 0 & \id_X \end{pmatrix}: X \oplus X \to X \oplus X$ is an isomorphism (\cref{rmk:semiadd-cogp}). This map is natural with respect to all maps, so if $X \to Y \to Z$ is a cofiber sequence and this map is an isomorphism for $X$ and for $Y$, then it is also an isomorphism for $Z$.
\end{proof}

\begin{prop}
The sub-operad $\Cat_{\rex,\add,\trivsusp} \subset \Cat_{\rex,\add}$ is a $\otimes$-ideal.
\end{prop}
\begin{proof}
Let $\calC \in \Cat_{\rex,\add,\trivsusp}$ and $\calD \in \Cat_{\rex,\add}$; we wish to show that $\calC \otimes \calD$ has trivial suspension. But this follows from the fact that $(-) \otimes \calD$ is a 2-functor, and locally commutes with zero objects and finite colimits (and hence with suspensions).
%
\end{proof}

\begin{Def}
Let $\Cat_{\rex,\add,\trivsusp,\Split}$ be the following $\infty$-operad. An object is an object $\calC \in \Cat_{\rex,\add,\trivsusp}$ equipped with, for each prime $p$, a section of the map $X \to X/p$ natural in $X \in \calC$. The hom-space $\Cat_{\rex,\add,\trivsusp,\Split}(\calC_1,\dots,\calC_n;\calD)$ is the subgroupoid of $\Cat_{\rex,\add,\trivsusp}(\calC_1,\dots,\calC_n;\calD)$ comprising those functors $\calC_1\times \cdots \times \calC_n \to \calD$ which commute with these natural splittings separately in each variable. Composition is as in $\Cat_{\rex,\add,\trivsusp}$.
\end{Def}


\begin{Def}\label{def:evconst}
The additive 1-category $\eventuallyconstantcat$ is defined as follows. An object $(V,F)$ comprises an infinite tuple $(V_2,V_3,V_5,\dots)$ of finite-dimensional $\field_p$-vector spaces, one for each prime $p$, a finitely-generated free abelian group $F$, and isomorphisms $F/p \cong V_p$ for all but finitely-many primes. A morphism $f: (V,F) \to (W,G)$ comprises a tuple of $\field_p$-linear morphisms $f_p: V_p \to W_p$ (one for each prime $p$) and a group homomorphism $f_0: F \to G$ such that for all but finitely-many primes $p$ we have $f_p = f_0/p$ under the canonical identifications $V_p = F/p$ and $W_p = G/p$. Composition and $\otimes$ are defined in the evident ``componentwise" manner. 
\end{Def}

\begin{notation}\label{notation:ev-const}
The category $\eventuallyconstantcat$ is equipped with a canonical object, which we denote $S$, corresponding to the free abelian group $\ints$ with its reduction mod $p$ at each prime $p$. For a prime $p$, we denote by $S/p$ the object of $\eventuallyconstantcat$ corresponding to $\ints/p$ at the prime $p$ and $0$ at all other primes. For $n \in \ints$, we denote by $S(m)$ the object of $\eventuallyconstantcat$ whose free part is $\ints$, with its standard at all $\ell$ not dividing $m$, but whose component at the prime $p$ is $0$ for all $p$ dividing $m$. Note that every object is canonically of the form $(S/p_1)^{\oplus e_1}\oplus \cdots \oplus (S/p_r)^{\oplus e_r} \oplus S(p_1\cdots p_r)^{\oplus f}$, and that if $p$ is coprime to $m$, there is a canonical isomorphism $S(m) \cong S(p) \oplus S(mp)$.
\end{notation}

\begin{lem}\label{lem:evconst-split}
The category $\eventuallyconstantcat$ admits in a canonical way the structure of an object of $\Cat_{\rex,\add,\trivsusp,\Split}$.
\end{lem}
\begin{proof}
Clearly $\eventuallyconstantcat$ is an additive 1-category.

Let us show that $\eventuallyconstantcat$ has cofibers. Observe that $\Hom(S/p,S/\ell) = 0$ for $p \neq \ell$, $\Hom(S/p,S/p) = \ints/p$, $\Hom(S/p,S(p)) = \Hom(S(p),S/p) = 0$, and $\Hom(S(m),S(m))$ contains a copy of $\ints$ corresponding to those maps which are the reduction mod $\ell$ of a fixed map for every $\ell \notdivides m$. So in light of the biproduct decompositions noted in \cref{notation:ev-const} and the triviality of suspension, it will suffice to construct cofibers of maps $\phi: (S/p)^{\oplus e} \to (S/p)^{\oplus e'}$ and $\psi: S(m)^{\oplus f} \to S(m)^{\oplus f'}$ where $\psi$ is the reduction of an integral map at all primes not dividing $m$. In the former case, cofibers are computed as in $\Vect_{\field_p}$; this clearly works because if $(S/p)^{\oplus e'} \to X$ is a morphism in $\eventuallyconstantcat$, it must factor through the inclusion $X/p \to X$. In the latter case, let $Q$ be the cokernel (in the category of abelian groups) of the corresponding map $\psi: \ints^f \to \ints^{f'}$. Let $T \subseteq Q$ be the torsion subgroup and $F = Q / T$ the torsion-free quotient.
Let $\overline Q$ be the object of $\eventuallyconstantcat$ whose component at $p$ is $Q/p$ for $p \notdivides m$ and $0$ for $p \divides m$; the free part is $F$. To see that $\overline Q$ is the cokernel of $\psi$, consider first a map $\theta: S(m)^{\oplus f'} \to S/p$. Then $\theta$ must kill $S(mp)^{\oplus f'}$. If $p \divides m$, this means that $\theta = 0$, and so $\theta \psi = 0$ and $\theta$ factors uniquely through $Q$ via the zero map. Otherwise, $\theta$ factors through $S(m)^{\oplus f'} \to (S/p)^{\oplus f'}$. Thus $\theta \psi = 0$ if and only if $\theta (\psi/p) = 0$, if and only if $\theta$ factors uniquely through $Q/p$. As any factorization through $\overline Q$ must kill $\overline Q(p)$, this implies that the factorization through $\overline Q$ is unique. Next, consider a map $\theta: S(m)^{\oplus f'} \to S(mn)$, where $n$ is the torsion exponent of $Q$. Then $\theta \psi = 0$ if and only if $\theta$ factors uniquely through $F$. This gives a unique factorization through $\overline Q$.

Finally, let us construct our splittings. The splitting of $S/\ell \to (S/\ell)/p$ is the identity when $p=\ell$ and zero when $p \neq \ell$. The splitting of $S(m) \to S(m)/p$ is zero when $p \divides m$ and the canonical inclusion $S/p \to S(m) = S/p \oplus S(mp)$ when $p \notdivides m$. We extend these definitions by taking direct sums. In doing this, we must check that our definitions are consistent with the relation $S(m) = S/p \oplus S(mp)$ when $p \notdivides m$. 
They are. Now we check that these splitting maps are natural. They are.
\end{proof}

\begin{prop}\label{prop:evconst-univ}
Evaluation at $S \in \eventuallyconstantcat$ determines an equivalence of categories $\Cat_{\rex,\add,\trivsusp,\Split}(\eventuallyconstantcat, \calD) \simeq \calD^\sim$ for any $\calD \in \Cat_{\rex,\add,\trivsusp,\Split}$.
\end{prop}
\begin{proof}
Let $D \in \calD$, and let $F: \eventuallyconstantcat \to \calD$ be a morphism of \\ $\Cat_{\rex,\add,\trivsusp,\Split}$ with $F(S) = D$. Then $\id_S$ must be carried to $\id_D$, and integral multiples of $\id_S$ must be carried to the corresponding integral multiples of $\id_D$. $F$ must carry $S/p$ to $D/p$ and $S(p)$ to the cokernel of the canonical splitting $D/p \to D$, which we denote $D(p)$. Note that $(D/p)/\ell = 0$ for $\ell \neq p$. Thus the idempotents defining the various $D(p)$'s commute (with their product being $0$ for $p \neq \ell$), and it follows that for any integer $m$, $S(m)$ is carried to the intersection of the $D(p)$'s for $p \divides m$; we denote this object $D(m)$. Because $F$ commutes with direct sums, its behavior on objects is now entirely forced. A nonzero morphism $S/p \to S/\ell$ exists only if $p=\ell$, in which case it is the reduction mod $p$ of a morphism $S \to S$, so its image in $\calD$ is forced. There are no nonzero morphisms $S/p \to S(m)$ or $S(m) \to S/p$ for $p \notdivides m$. Any morphism $S(m) \to S(m)$ which is the reduction of a morphism $S \to S$ has image in $\calD$ forced as well. Now, any morphism in $\eventuallyconstantcat$ is a sum of direct sums of the morphisms considered so far, and thus its image in $\calD$ is determined by additivity.

So if such an $F$ exists, it is unique. Let us verify that such an $F$ in fact exists. We have seen that we must have $F((S/p_1)^{\oplus e_1}\oplus \cdots \oplus (S/p_r)^{\oplus e_r} \oplus S(p_1\cdots p_r)^{\oplus f} = (D/p_1)^{\oplus e_1}\oplus \cdots \oplus (D/p_r)^{\oplus e_r} \oplus D(p_1\cdots p_r)^{\oplus f}$. This is well-defined on objects because $D(m) = D(pm)\oplus D/p$ for $p \notdivides m$. To see that $F$ is well-defined on morphisms, we must check first that $\calD(D/p,D/p)$ is $p$-torsion. This is indeed the case by the universal property of $D/p$. We must also check, for $p \notdivides m$, that the identity on $D(m)$ agrees with the sum of the identity on $D/p$ and the identity on $D(mp)$, which it does. It is clear that $F$ commutes with addition on hom-sets. To see that $F$ is functorial is straightforward. It is clear that $F$ commutes with the canonical splittings objects $X \to X/p$ -- for instance, when $X = S$, this holds by definition. $F$ commutes with finite coproducts, so it remains only to check that $F$ commutes with cofibers. This is clear from the description given in \cref{lem:evconst-split}.
\end{proof}

\begin{thm}
The $\infty$-operad $\Cat_{\rex,\add,\trivsusp,\Split}$ is a symmetric monoidal $\infty$-cate\-gory, and the inclusion $\Cat_{\rex,\add,\trivsusp,\Split} \to \Cat_{\rex,\add,\trivsusp}$ preserves binary tensor products up to equivalence. The unit object is $\eventuallyconstantcat$.
\end{thm}
\begin{proof}
We have seen in \cref{prop:evconst-univ} that the unit is representable. So it will suffice to show that the tensor product is representable and that the associativity constraints are isomorphisms.

Let $\calC,\calD \in \Cat_{\rex,\add,\trivsusp,\Split}$, and let $U\calC,U\calD$ be the underlying objects in $\Cat_{\rex,\add,\trivsusp}$. We start by showing that $U\calC \otimes U\calD$ admits splittings making it an object of $\Cat_{\rex,\add,\trivsusp,\Split}$. In fact, the data of such splittings on $\calE \in \Cat_{\rex,\add,\trivsusp}$ is equivalent
to the data of splittings $\calE = \calE/p \times \calE(p)$ for each prime $p$, where $\calE/p$ has the property that its hom-spaces are $p$-torsion and $\calE(p)$ has the property that the endomorphism $p:E \to E$ is an epimorphism for any $E \in \calE$. Moreover, the tensor product on $\Cat_{\rex}$, and hence on $\Cat_{\rex,\add,\trivsusp}$, preserves finite products (which are also finite coproducts in $\Cat_\rex$ and hence also in $\Cat_{\rex,\add,\trivsusp}$) separately in each variable: one way to see this is that the category $\Cat_{\rex}$ is semiadditive, with the addition operation on hom-spaces being $\oplus$; since $\otimes$ preserves this addition in each variable separately, it must preserve direct sums of objects of $\Cat_\rex$ separately. So if $\calC = \calC/p \oplus \calC(p)$ and $\calD = \calD/p \oplus \calD(p)$, then $\calC \otimes \calD = (\calC/p \otimes \calD/p) \oplus (\calC/p \otimes \calD(p)) \oplus (\calC(p) \otimes \calD/p) \oplus (\calC(p) \otimes \calD(p))$. The middle two terms have the property that $p=0$ and $p$ is an epimorphism on each object, so they vanish. We are left with the first term, which has the property that $p=0$ on each object, and the last term, which has the property that $p$ is an epimorphism on each object. That is, we have the desired splitting.

Moreover, from this description we see that if $F: \calC \times \calD \to \calE$ preserves these splitting separately in each variable, then the induced functor $\overline F: \calC \otimes \calD \to \calE$ preserves the designated splittings as well. Thus $\Cat_{\rex,\add,\trivsusp,\Split} \subset \Cat_{\rex,\add,\trivsusp}$ is closed under binary tensor products.
\end{proof}

\begin{cor}\label{cor:evconst-mon-init}
The initial object of $\SMC_{\rex,\add,\trivsusp,\Split}$ is $\eventuallyconstantcat$.
\end{cor}
\begin{proof}
The unit of a symmetric monoidal category, with its unique $E_\infty$ structure, is always the initial object of the category of $E_\infty$-algebras.
\end{proof}

\begin{prop}\label{prop:init-add}
The initial object of $\SMD_{\gptrivsusp,\cof}$ is $\eventuallyconstantcat$.
\end{prop}
\begin{proof}
There is a fully faithful forgetful functor $\SMD_{\gptrivsusp,\cof} \to \SMC_{\rex,\add,\trivsusp,\Split}$, where the splittings come from the clopen idempotent structure of $S/p$. So in order to verify that $\eventuallyconstantcat$ is the initial object of $\SMD_{\gptrivsusp,\cof}$, it will suffice by \cref{cor:evconst-mon-init} to verify that $\eventuallyconstantcat$ has duals for objects. But it is clear that every object is self-dual.
\end{proof}

\subsubsection{All together}\label{subsec:all}

In this subsection, we product together the initial objects from the rest of \cref{sec:init} to describe the initial object of the $\infty$-category of symmetric monoidal $\infty$-categories with duals, finite coproducts, and cofibers.

\begin{thm}\label{thm:freerigcopcof}
The initial object of $\SMD^{\amalg,\cof}$ is $\Spt^\fin \times \eventuallyconstantcat \times \Span(\Fin)$.
\end{thm}
\begin{proof}
By \cref{prop:3-split}, the initial object $\calI$ of $\SMD^{\amalg,\cof}$ splits as $\calI = \calI_\stab \times \calI_\gptrivsusp \times \calI_\notgp$. Moreover, the localization functors $\calI \to \calI_\stab$, $\calI \to \calI_\gptrivsusp$, $\calI \to \calI_\notgp$ exhibit their codomains respectively as the initial objects of $\SMD_\stab$, $\SMD_{\gptrivsusp,\cof}$, and $\SMD_{\notgp,\cof}$ respectively. So this follows from \cref{prop:init-stab}, \cref{prop:init-add}, and \cref{cor:init-antiadd}.
\end{proof}

\begin{thm}
The initial object $\calI$ of $\SMD^\rex$ is $\Spt^\fin \times \eventuallyconstantcat \times \calI_\notgp$.
\end{thm}
\begin{proof}
We have $\SMD^\rex_\stab = \SMD_\stab^{\amalg,\cof}$ and $\SMD^\rex_{\gptrivsusp} = \SMD^{\amalg,\cof}_\gptrivsusp$. So these two factors agree with the ones from \cref{thm:freerigcopcof}. The final factor will not be exactly $\Span(\Fin)$.
\end{proof}

\begin{rmk}
We do not know precisely what the category $\calI_\notgp$ is. 
However, $\calI_\notgp$ receives a functor from $\Span(\Fin)$ by \cref{thm:freerigcopcof}.
\end{rmk}

\section{Application to Equivariant Homotopy Theory}\label{chap:app}

In this chapter, we apply the foregoing results to prove \cref{cor:main-result}, which gives a new universal property for equivariant stable homotopy theory. We give several further, essentially equivalent universal properties in \cref{cor:maincor}. There are a number of ways to formulate this universal property, which was suggested in a preliminary form by Charles Rezk (\cite{rezk}).

\subsection{Equivariant homotopy theory}
Let $G$ be a compact Lie group, and let $\Tope{G}$ be the $\infty$-category of $G$-spaces, considered as symmetric monoidal under cartesian product. Let $\Tope{G}^\fin \subset \Tope{G}$ be the full symmetric monoidal subcategory of $G$-spaces with finitely many cells. We also have pointed versions $\Tope{G}_\ast, \Tope{G}^\fin_\ast$, considered as symmetric monoidal under smash product. These are to be contrasted to the $\infty$-categories $\Top^{BG}, \Top^{BG}_\ast$ of \defterm{Borel} $G$-spaces (pointed and unpointed respectively). Consider also the symmetric monoidal $\infty$-category $\Spte G$ of genuine $G$-spectra, and the full subcategory $\Spte G^\fin$ of finite genuine $G$-spectra.

Let $\SMC_{\rex}$ denote the $\infty$-category of symmetric monoidal $\infty$-categories with finite colimits. Let $\SMC_{\ast,\rex}$ denote the $\infty$-category of pointed symmetric monoidal $\infty$-categories with finite colimits. Let $\SMD_\rex$ denote the $\infty$-category of symmetric monoidal $\infty$-categories with finite colimits and duals for objects. Let $\freeduals_{\rex}: \SMC_{\rex} \to \SMD_{\rex}$ denote the left adjoint to the inclusion, and let $\freeduals_\rex: \SMC_{\ast,\rex} \to \SMD_\rex$ also denote the left adjoint to that inclusion. Let also $\SMP$ denote the $\infty$-category of presentably symmetric monoidal $\infty$-categories and left adjoint symmetric monoidal functors.

Before going further, let us remind ourselves of some basic facts about equivariant homotopy theory:

\begin{lem}\label{lem:induce-finite}
Let $G$ be a compact Lie group and $H \subseteq G$ a closed subgroup. Then induction from $H$-spaces to $G$-spaces carries finite $H$-CW-complexes to finite $G$-CW-complexes.
\end{lem}
\begin{proof}
Induction carries $H$-orbits to $G$-orbits, so this follows by induction on cells, since induction preserves homotopy colimits.
\end{proof}

\begin{lem}\label{lem:manifold-finite}
Any compact $G$-manifold $M$ is $G$-homotopy equivalent to a finite $G$-CW-complex. In particular, if $S^V$ is a $G$-representation sphere, then it is $G$-homotopy equivalent to a finite $G$-CW-complex.
\end{lem}
\begin{proof}
We will prove the theorem by induction on the dimension $d$ of $M$. By \cite[Corollary 4.11]{wasserman}, there exists an equivariant Morse function on $M$, which decomposes $M$ via a finite filtration $M_0 \subset M_1 \subset \cdots \subset M_n$. We have that $M_{i+1} = M_i \cup_{S(V)}D(v)$, where $V$ is an equivariant bundle over a $G$-orbit $G/H$, and $S(V) \subset D(V)$ are its associated sphere and disc bundles. Of course, $D(V)$ is $G$-homotopy equivalent to the orbit $G/H$, which is finite cellular by definition. The dimension of $D(V)$ is at most the dimension $d$ of $M$, so the dimension of $S(V)$ is strictly less. By induction $S(V)$ is also finite cellular. Thus, inducting on $i$ we obtain that $M$ is finite cellular as desired.
\end{proof}

\begin{lem}\label{lem:spt-g-duals}
The symmetric monoidal $\infty$-category $\Spte G^\fin$ of genuine finite $G$-spectra has all duals. Likewise, the symmetric monoidal $\infty$-category $\Spte{G}^\fd$ of finitely-dominated $G$-spectra has all duals.
\end{lem}
\begin{proof}
The second statement follows from the first because dualizable objects are closed under retracts (\cref{lem:dual-stable}). For the first statement, we may assume by induction that this is true for all proper closed subgroups $H \subset G$.
Following \cite[Theorem 4.10, Theorem 4.17]{greenlees-may}, it follows from the Wirthm\" uller isomorphism that the dual of an orbit $G/H_+$ is given by inducing up $S^{-L(H)}$ from $H$-spectra, where $L(H)$ is the tangent space of the identity of $G/H$. By \cref{lem:manifold-finite}, $S^{L(H)}$ is a a finite $H$-spectrum. If $H = G$, then $L(H) = 0$ and $S^{L(H)} = S^0$ is trivially a finite $H$-spectrum. Otherwise, $S^{-L(H)} = (S^{L(H)})^\dual$ is an $H$-spectrum by induction on the subgroups of $G$. Therefore by \cref{lem:induce-finite}, the dual of $(G/H)_+$ is a finite $G$-spectrum. Since $\Spte G^\fin$ is by definition generated under finite colimits and desuspensions by the orbits, it now follows from \cref{lem:dual-stable} that $\Spte G^\fin$ has all duals.
\end{proof}


We noted in Example \ref{eg:twist-G-sphere} that for any $G$-representation $V$, the 1-point compactification $S^V$ has twisted-trivial braiding in $\Tope{G}_\ast$. Let $S \subset \Tope{G}$ denote the set of finite-dimensional $G$-representation spheres. For each $S^V \in S$, \cref{prop:monoidal-splitting} tells us that the $\infty$-category $\freeduals_S^\SMP \Tope{G}_\ast$ decomposes as $\freeduals_S^\SMP \Tope{G}_\ast = \freeduals_S^\SMP \Tope{G}_\ast[{S^V}\inv] \times \freeduals^\SMP_S \Tope{G}_\ast / S^V$ where in the first factor $S^V$ is $\wedge$-invertible, while in the other it is trivial.

Our present goal is to show that the second factor vanishes, i.e.

\begin{thm}\label{thm:factor-killed}
Let $G$ be a compact Lie group, let $S \subset \Tope{G}$ denote the set of finite-dimensional $G$-representation spheres, and let $S^V \in S$. Then $\freeduals^\SMP_S \Tope{G}[(S^1)\inv]/S^V = 0$.
\end{thm}

The proof of \cref{thm:factor-killed} is deferred to the next section, \cref{sec:proof}. In the meantime, let us pause to deduce a major result of this thesis:


\begin{cor}\label{cor:main-result}
Let $G$ be a compact Lie group, and let $S \subset \Tope{G}_\ast$ denote the set of finite-dimensional representation spheres. Then $\freeduals_S^\SMP \Tope{G}_\ast[(S^1)\inv] = \Spte{G}$. That is, $\Spte{G}$ is the free presentably symmetric monoidal stable $\infty$-category on $\Tope{G}_\ast$ where the representation spheres become dualizable.
\end{cor}
\begin{proof}
Let $V$ be a finite-dimensional $G$-representation. By \cref{thm:factor-killed}, \[\freeduals_S^\SMP \Tope{G}_\ast[(S^1)\inv] / S^V = 0\]
By the product decomposition, this implies that
\[\freeduals_S^\SMP \Tope{G}_\ast[(S^1)\inv] = \freeduals_S^\SMP \Tope{G}^\fin[(S^1)\inv, (S^V)\inv]\]
i.e. $S^V$ is invertible in $\freeduals_S^\SMP \Tope{G}_\ast[(S^1)\inv]$. Now, $\Spte{G}$ is freely obtained from $\Tope{G}_\ast$ by inverting all the representation spheres $S^V$ (see \cite[Theorem A.2]{cmnn} for the case where $G$ is finite, and \cite[Corollary C.7]{gepner-meier} in general, in both cases the result is deduced from results of \cite{robalo}). By this universal property, we obtain a symmetric monoidal, left adjoint functor $\Spte{G} \to \freeduals_S^\SMP \Tope{G}_\ast[(S^1)\inv]$. The presentably symmetric monoidal $\infty$-category $\Spte{G}$ is stable, and by \cref{lem:spt-g-duals}, the representation spheres are dualizable in $\Spte{G}$, so we obtain a functor in the other direction as well by the universal property of $\freeduals_S^\SMP \Tope{G}_\ast$. The composite of the two functors in either order looks like the identity functor after precomposing with the canonical functor from $\Tope{G}_\ast$, and so by the universal property this composite is in fact the identity. Thus we have a symmetric monoidal equivalence under $\Tope{G}_\ast$ between $\Spte{G}$ and $\freeduals_S^\SMP \Tope{G}_\ast[(S^1)\inv]$ as claimed.
%
\end{proof}

We may deduce from \cref{cor:main-result} several related universal properties:

\begin{cor}\label{cor:maincor}
Let $G$ be a compact Lie group. Let $\Sigma^\infty_G : \Tope{G}_\ast \to \Spte{G}$ denote the equivariant suspension functor. Then
\begin{enumerate}
    \item\label{item:mc1} The functor $\Sigma^\infty_G$ is the universal presentably symmetric monoidal functor from $\Tope{G}_\ast$ to a presentably symetic monoidal stable $\infty$-category carrying each representation sphere to a dualizable object.
    \item\label{item:mc2} The functor $\Sigma^\infty_G$ is the universal presentably symmetric monoidal functor from $\Tope{G}_\ast$ to a presentably symmetric monoidal stable $\infty$-category carrying each object of $\Tope{G}^\fin_\ast$ to a dualizable object.
    \item\label{item:mc3} The functor $\Sigma^\infty_G$ is the universal presentably symmetric monoidal functor from $\Tope{G}_\ast$ to a presentably stable symmetric monoidal $\infty$-category carrying each compact object of $\Tope{G}_\ast$ to a dualizable object.
    \item\label{item:mc4} The functor $\Sigma^\infty_G$ is the universal compactly-generated symmetric monoidal functor from $\Tope{G}_\ast$ to a compactly-generated, stable symmetric monoidal $\infty$-category with duals for all compact objects.
\end{enumerate}  
Now let $\Sigma^{\infty,\fin}_G : \Tope{G}_\ast^\fin \to \Spte{G}^\fin$ denote the restriction / corestriction of $\Sigma^\infty_G$ to finite $G$-spaces / finite $G$-spectra.
\begin{enumerate}[resume]
    \item\label{item:mc5} The functor $\Sigma^{\infty,\fin}_G$ is the universal right exact symmetric monoidal functor from $\Tope{G}^\fin_\ast$ to a stable symmetric monoidal $\infty$-category carrying each representation sphere to a dualizable object.
    \item\label{item:mc6} The functor $\Sigma^{\infty,\fin}_G$ is the universal right exact symmetric monoidal functor from $\Tope{G}^\fin_\ast$ to a stable symmetric monoidal $\infty$-category carrying each object to a dualizable object.
    \item\label{item:mc7} The functor $\Sigma^{\infty,\fin}_G$ is the universal right exact symmetric monoidal functor from $\Tope{G}^\fin_\ast$ to a stable symmetric monoidal $\infty$-category with all objects dualizable.
\end{enumerate}
Finally, let $\Sigma^{\infty,\fd}_G : \Tope{G}_\ast^\fd \to \Spte{G}^\fd$ denote the restriction / corestriction of $\Sigma^\infty_G$ to finitely-dominated $G$-spaces / finitely-dominated $G$-spectra.
\begin{enumerate}[resume]
    \item\label{item:mc8} The functor $\Sigma^{\infty,\fd}_G$ is the universal right exact symmetric monoidal functor from $\Tope{G}^\fd_\ast$ to a right exact, idempotent-complete stable $\infty$-category carrying each representation sphere to a dualizable object.
    \item\label{item:mc9} The functor $\Sigma^{\infty,\fd}_G$ is the universal right exact symmetric monoidal functor from $\Tope{G}^\fd_\ast$ to a right exact, idempotent-complete, stable symmetric monoidal $\infty$-category carrying each object to a dualizable object.
    \item\label{item:mc10} The functor $\Sigma^{\infty,\fd}_G$ is the universal right exact symmetric monoidal functor from $\Tope{G}^\fd_\ast$ to a right exact, idempotent-complete, stable symmetric monoidal $\infty$-category with all objects dualizable.
\end{enumerate}
\end{cor}
\begin{proof}
\cref{item:mc1} is the statement of \cref{cor:main-result}. \cref{item:mc2} then follows because every finite $G$-spectrum is dualizable (\cref{lem:spt-g-duals}), and likewise \cref{item:mc3} follows because every finitely-dominated $G$-spectrum is dualizable (again \cref{lem:spt-g-duals}). Recall now that a symmetric monoidal left adjoint between compactly-generated symmetric monoidal $\infty$-categories is said to be \defterm{compactly-generated} if it preserves compact objects (or equivalently, if its right adjoint preserves filtered colimits). So for \cref{item:mc4}, it suffices to verify that $\Spte{G}$ has duals for all compact objects (which is \cref{lem:spt-g-duals}), that the functor $\Sigma^\infty_G$ preserves compact objects (which it does), and that if $\calK$ is any compactly-generated symmetric monoidal $\infty$-category with duals for compact objects, and if $F : \Tope{G}_\ast \to \calK$ is a compactly-generated symmetric monoidal functor, then the functor $\tilde F: \Spte{G} \to \calK$ induced by \cref{item:mc3} also preserves compact objects. This last point follows because the compact objects $\Spte{G}^\fd \subset \Spte{G}$ are contained in (in fact, coincide with) the closure of the image of $\Sigma^{\infty,\fd}_G$ under idempotent splitting; by hypothesis, $\tilde F$ carries the objects in the image of $\Sigma^{\infty,\fd}_G$ to dualizable objects, and dualizable objects are closed under retracts (\cref{lem:dual-stable}), so the result follows.

For \cref{item:mc5}, let $F: \Tope{G}^\fin_\ast \to \calC$ be a right exact symmetric monoidal functor to a stable symmetric monoidal $\infty$-category carrying each representation sphere to a dualizable object. We may assume without loss of generality that $\calC$ is small. There is an induced symmetric monoidal left adjoint $\Ind(F) : \Tope{G}_\ast \to \Ind(\calC)$, which continues to carry each representation sphere to a dualizable object. Thus by \cref{item:mc1} we obtain an esssentially unique extension $\widetilde{\Ind(F)} : \Spte{G} \to \Ind(\calC)$. The restriction $\tilde F : \Spte{G}^\fin \to \Ind(\calC)$ is right exact symmetric monoidal. Moreover, every object of $\Spte{G}^\fin$ is in the image of $\Sigma^{\infty,\fin}_G$ and therefore is contained in $\calC$. Thus $\tilde F$ corestricts to an extension $\bar F : \Spte{G}^\fin \to \calC$ of $F$. We have $\Ind(\bar F) = \widetilde{\Ind(F)}$; because $\Ind$ is a fully faithful functor, the essential uniqueness of $\bar F$ follows from the essential uniqueness of $\widetilde{\Ind(F)}$. For \cref{item:mc6}, the argument is similar, using \cref{item:mc2} instead of \cref{item:mc1}. \cref{item:mc7} follows because every object of $\Spte{G}^\fin$ is in fact dualizable (\cref{lem:spt-g-duals}).

\cref{item:mc8} follows from \cref{item:mc5}, using that $\Spte{G}^\fd$ is the idempotent completion of $\Spte{G}^\fin$. Similarly, \cref{item:mc9} follows from \cref{item:mc6}, and \cref{item:mc10} follows from \cref{item:mc7}.
\end{proof}

\begin{rmk}
Via the equivalence provided by the $\Ind$ functor (known as $\infty$-categorical Gabriel-Ulmer duality), \cref{item:mc8,item:mc9,item:mc10} also have equivalent statements in terms of compactly-generated symmetric monoidal $\infty$-categories. For example, under Gabriel-Ulmer duality the statement equivalent to \cref{item:mc10} is \cref{item:mc4}; we have not recorded the other two equivalent statements explicitly.
\end{rmk}

\subsection{The proof of \cref{thm:factor-killed}}\label{sec:proof}

This section is devoted to the proof of \cref{thm:factor-killed}. The proof is by induction on the structure of the orbit category $\calO_G$. More precisely, \cref{def:conj-sphere} associates to each \defterm{interval} $I$ in the lattice of conjugacy classes of closed subgroups $H \subseteq G$, a representation sphere $S^I$. Starting from the assumption of a symmetric monoidal left adjoint $F :\Tope{G}_\ast \to \calC$ with $\calC$ stable presentably symmetric monoidal with compact unit and $F(S^V) = 0$ for some representation sphere $V$, a geometric argument first shows (\cref{lem:conc-killed}) that if $I = \overline{\{H\}}$ is a ``singleton" interval, then $F(S^{\overline{\{H\}}}) = 0$. Then comes (\cref{thm:factor-killed-bis}) an induction on the size of the interval $I$, showing that $F(S^I) = 0$ for any interval $I$ (it is in this step that we use that hypothesis that $\calC$ have compact unit -- this is one delicate point). Taking $I=G$, we have $F(S^0) = 0$, which implies that $\calC = 0$. In a final step, the hypothesis that $\calC$ have a compact unit is removed, and \cref{thm:factor-killed} follows by the defining universal property of $\freeduals^\SMP_S \Tope{G}_\ast[(S^1)\inv]/S^V$.

The argument involves threading the needle between the infinitary setting of $\freeduals^\SMP_S \Tope{G}_\ast[(S^1)\inv]/S^V$ and the finitary setting of $\freeduals^\rex_S \Tope{G}^\fin_\ast[(S^1)\inv]/S^V$. The main argument of \cref{thm:factor-killed-bis} is carried out in the infinitary setting but with a compactness assumption. On the one hand, the infinitary setting is necessary in order to exploit a certain infinite colimit (\cref{lem:w+action}), while on the other hand the compactness is necessary in the inductive argument of \cref{thm:factor-killed-bis}\footnote{In fact, if $G$ is finite the argument will go through without the compactness hypothesis.} as alluded to above. Luckily, the universal properties involved here are rather robust (cf. \cref{cor:maincor}), and it is possible to pass back and forth between these settings rather easily.

\begin{lem}[Untwisting Lemma]\label{lem:untwist}
For any $X \in \Tope{G}_\ast$, there is a $G$-homeomorphism between $G_+ \wedge X$ with the diagonal action, and $G_+ \wedge X$ with the action on $G_+$ (forgetting the original action on $X$).
\end{lem}
\begin{proof}
We define a map $\phi: G \times X \to G \times X$ by $(g,x) \mapsto (g,gx)$. This is equivariant from the left action $\cdot^l$ to the diagonal action $\cdot^d$: $\phi(h\cdot^l (g,x)) = \phi((hg,x)) = (hg, hgx)$ while $h \cdot^d\phi(g,x)) = h\cdot^d(g,gx) = (hg,hgx)$. In the other direction, define $\psi(g,x) = (g,g\inv x)$. Then $\phi$ and $\psi$ are inverse to each other (and so $\psi$ is also equivariant). Moreover, these maps descend: we have $G_+ \wedge X = G \times X / G \times \ast$, and $\phi(g,\ast) = (g,g\ast) = (g,\ast)$ while $\psi(g,\ast) = (g,g\inv \ast) = (g,\ast)$.
\end{proof}

\begin{lem}
For any $H \subseteq G$, and any finite-dimensional $G$-representation $V$, the underlying space of $(S^V)^H$ is a sphere.
\end{lem}
\begin{proof}
The $H$-fixed points $V^H$ in the representation $V$ are a linear subspace. So $(S^V)^H = S^{V^H}$ is a sphere.
\end{proof}

\begin{Def}\label{def:conj-sphere}
Recall that we have fixed a compact Lie group $G$. An \defterm{upset} is an upwards-closed set in the lattice of conjugacy classes of subgroups of $G$. A \defterm{downset} is a downwards-closed set in the lattice of conjugacy classes of subgroups of $G$, i.e. the complement of an upset. A \defterm{interval} is the intersection of an upset and a downset. If $S$ is a subset of the poset of closed subgroups of $G$, write $\overline S$ for the closure of $S$ under conjugacy in $G$. In particular, if $H \subseteq G$ is a closed subgroup, then $\overline{\{H\}}$ denotes the conjugacy closure of the singleton set $\{H\}$, and $\overline{\downarrow H}$ denotes the collection of subgroups conjugate to a subgroup of $H$.

For any interval $I$, define $S^I$ to be the following pointed $G$-space, viewed as a presheaf on the orbit category: we have $(S^I)^H = S^0$ for $H \in I$ and $(S^I)^H = 0$ for $H \not \in I$. All transition maps between $S^0$'s are identities, and all other transition maps are zero as they must be.
\end{Def}

\begin{lem}
For any interval $I$, we have $S^I \wedge S^I  = S^I$ canonically. More generally, for any intervals $I,J$, we have $S^I \wedge S^J = S^{I \cap J}$, canonically.
\end{lem}

\begin{lem}\label{lem:cofseq}
If $D$ is a downset and $U$ is the complementary upset, then there is a natural cofiber sequence $S^D \to S^0 \to S^U$.
\end{lem}

\begin{lem}\label{lem:conc-borel}
For any conjucacy closure of a singleton $\overline{\{H\}}$, the full subcategory of objects of the form $X = S^{\overline{\{H\}}} \wedge Y$ consists of those pointed presheaves $X$ on $\calO_G$ which send all orbits other than $G/H$ to $0$ (equivalently, pointed $G$-spaces $X$ such that $X^K$ is contractible for $K$ not conjugate to $H$). This is monoidally equivalent to the $\infty$-category $\Top_\ast^W$ of pointed Borel $W$-spaces, where $W = W_G(H) = N_G(H)/H$ is the Weyl group of $H$ in $G$.
\end{lem}
This category is referred to as the category of \defterm{pointed $G$-spaces concentrated at $H$}.
\begin{proof}
By Elemendorf's theorem, $\Tope{G}_\ast$ is equivalent to the category of $\Top_\ast$-valued presheaves on the orbit category $\calO_G$; the equivalence sends a $G$-space $X$ to the functor $G/H \mapsto X^H$. The smash product is computed levelwise in this presheaf category. Thus smashing with $S^{\overline{\{H\}}}$ kills the $K$-fixed points for $K$ not conjugate to $H$, and leaves them unchanged for $K=H$, verifying the first statement. From this we see that the category of pointed $G$-spaces concentrated at $H$ is equivalent to the category of pointed presheaves on the full subcategory of the orbit category on the orbit $G/H$. Recall that as a $G$-space, every endomorphism of $G/H$ is an automorphism, and its automorphism group is the Weyl group $W = W_G(H) = N_G(H) / H$. Thus this full subcategory is a connected $\infty$-groupoid equivalent to the 1-object $\infty$-groupoid $BW$ whose automorphism group is $W$. As presheaves on $BW$ are naturally identified with Borel $W$-spaces, the second statement follows. 
\end{proof}

\begin{lem}\label{lem:w+corr}
Under the equivalence of \cref{lem:conc-borel}, $S^{\overline{\{H\}}} \wedge (G/H)_+$ corresponds to $W_+ \in \Top_\ast^W$. Moreover, $S^{\overline{\{H\}}}$ corresponds to $S^0$, and the equivalence of course respects suspension.
\end{lem}
\begin{proof}
The equivalence functor from $H$-concentrated spaces to $W$-spaces is given by taking $H$-fixed points. And indeed, the $H$-fixed points of $(G/H)_+$ are $W_+$ and the $H$-fixed points of $S^{\overline{\{H\}}}$ are $S^0$. The final statement is true of any equivalence.
\end{proof}

\begin{lem}\label{lem:w+action}
For any compact Lie group $W$, the Borel space $W \in \Top^{BW}$ admits a free action of $W$, whose homotopy orbits are the terminal object $\pt \in \Top^{BW}$. Thus, $W_+ \in \Top_\ast^W$ admits a $W$-action whose pointed homotopy orbits are $S^0 \in \Top_\ast^W$.
\end{lem}
\begin{proof}
The second statement follows from the first because the functor $(-)_+: \Top^{BW} \to \Top_\ast^W$, which adds a disjoint basepoint, preserves colimits. If our $W$-spaces are left $W$-spaces, then the action in question comes from $W$ acting on itself on the right. Because $\pt \in \Top^{BW}$ is terminal, there is a map $W_{hW} \to \pt$. Because we are in the Borel setting, to see this is an equivalence it suffices to check on underlying spaces. Moreover, the ``underlying space" functor $\Top^{BW} \to \Top$ preserves colimits, so we may compute the map $W_{hW} \to \pt$ at the level of underlying spaces. Since the $W$-action on $W$ is free, its homotopy orbits are its orbits, namely $\pt$. Thus the map $W_{hW} \to \pt$ is an equivalence as claimed.
\end{proof}

\begin{lem}\label{lem:conc-killed}
Let $F : \Tope{G}_\ast \to \calC$ be a symmetric monoidal left adjoint to a stable presentably symmetric monoidal $\infty$-category $\calC$. Let $V$ be a finite-dimensional $G$-representation, and assume that $F(S^V) = 0$. Then for any conjugacy closure of a singleton $\overline{\{H\}}$, $F(S^{\overline{\{H\}}})=0$.
\end{lem}
\begin{proof}
We have that $F(S^{\overline{\{H\}}} \wedge S^V \wedge (G/H)_+)=0$, and under the equivalence of \cref{lem:conc-borel}, the space $S^{\overline{\{H\}}} \wedge S^V \wedge (G/H)_+$ corresponds to a pointed Borel $W$-space of the form $S^{V^H} \wedge W_+$ by \cref{lem:w+corr}. By the untwisting lemma (\cref{lem:untwist}), we have in turn that $S^{V^H} \wedge W_+ \simeq S^n \wedge W_+$ where $S^n$ has the trivial action (and $n = \dim(V^H)$). By \cref{lem:w+corr} again, we have that $F(S^{\overline{\{H\}}} \wedge S^n \wedge (G/H)_+) = 0$. Now, $F$ commutes with suspension and $\calC$ is stable, so it follows that $F(S^{\overline{\{H\}}} \wedge (G/H)_+) = 0$. Using \cref{lem:w+corr} again, \cref{lem:w+action} tells us that there is a $W$-action on $S^{\overline{\{H\}}} \wedge (G/H)_+$ whose homotopy orbits are $S^{\overline{\{H\}}}$. Because $F$ preserves homotopy orbits (in fact it preserves all colimits), it follows that $F(S^{\overline{\{H\}}})=0$.
\end{proof}

\begin{thm}\label{thm:factor-killed-bis}
Let $G$ be a compact Lie group. Let $F : \Tope{G}_\ast \to \calC$ be a symmetric monoidal left adjoint to a stable presentably symmetric monoidal $\infty$-category $\calC$ with a compact unit. Let $V$ be a finite-dimensional $G$-representation, and assume that $F(S^V) = 0$. Then $\calC = 0$.
\end{thm}

\begin{proof}
Say that a map in $\Tope{G}_\ast$ is \defterm{local} if it is carried to an equivalence by $F$. Let $\calU$ be the collection of all upsets $U$ such that $S^0 \to S^U$ is local. Then $\calU$ is downward-closed. For if $U \subseteq {U'}$ and $S^0 \to S^U$ is local, then because this map factors through $S^0 \to S^{U'}$ we have that $F(S^0)$ is a retract of $F(S^{U'})$. Therefore, we have that $F(S^{(U')^c}) = 0$, so that $S^0 \to S^{U'}$ is local and so ${U'} \in \calU$. Moreover, $\calU$ is closed under codirected intersections. For $S^{\cap_i U_i} = \varinjlim_i S^{U_i}$, and if the maps $S^0 \to S^{U_i}$ are local, then so is the map from $S^0$ to the colimit, as the universal functor $F$ preserves colimits.

Since $\calU$ is nonempty (containing the top element as $S^0 \to S^0$ is local), by the previous paragraph it satisfies the dual hypotheses of Zorn's lemma. So let $U$ be a minimal element of $\calU$. Then $U$ itself is closed under codirected intersections (note that the lattice of conjugacy classes of closed subgroups of $G$ is has codirected intersections for any topological group $G$). For if we suppose otherwise, then we have $H = \cap_i H_i$ with $H_i \in U$, $H \not \in U$, then $U = \cup_i (U \setminus \overline{\downarrow H_i})$, and so $S^U = \varinjlim_i S^{U \setminus \overline{\downarrow H_i}}$. As $S^0 \to S^U$ is local and $F(S^0)$ (the unit of $\calC$) is compact by hypothesis, it follows that $F(S^0)$ is a retract of some $F(S^{U \setminus \overline{\downarrow H_i}})$. Then as before, we see that $F(S^{(U\setminus \overline{\downarrow H_i)^c}}) = 0$, so that $S^0 \to S^{U \setminus \overline{\downarrow H_i}}$ is local, contradicting the minimality of $U$.

Suppose for contradiction that $U$ is nonempty. Then by the previous paragraph the dual hypotheses of Zorn's lemma apply and $U$ has a minimal element $H$. Then $S^{\overline{\downarrow H}} \to S^{\overline{\downarrow H}} \wedge S^U$ is local. But $F(S^{\overline{\downarrow H}} \wedge S^U) = F(S^{\overline{\{H\}}}) = 0$ by Lemma \ref{lem:conc-killed}. Therefore $F(S^{\overline{\downarrow H}}) = 0$, and so $S^0 \to S^{(G \setminus \overline{\downarrow H})}$ is local by Lemma \ref{lem:cofseq}. Therefore $S^0 = S^0 \wedge S^0 \to S^U \wedge S^{(G \setminus \overline{\downarrow H})} = S^{U \setminus \overline{\{H\}}}$ is local, contradicting the minimality of $U$.

Therefore $U$ must be empty, and so in fact $S^0 \to S^\emptyset = 0$ is local, i.e. $F(S^0) = 0$. But $F$ is strong monoidal and in particular $F(S^0)$ is the unit object of $\calC$. Thus $\calC = 0$.
\end{proof}

\begin{rmk}
In the beginning of the section, the proof of \cref{thm:factor-killed-bis} was characterized as an induction on the structure of intervals $I$, but in the actual proof this induction is phrased in terms of Zorn's lemma, obfuscating the content a bit. In the case where $G$ is finite, then there are finitely many upsets on $G$ and moreover it is immediate that any upset $U$ in $G$ has a minimal element $H$. So in this case, the two uses of Zorn's lemma may be avoided and replaced with a straightforward induction on the poset of upsets $U$ in $G$.
\end{rmk}

\begin{proof}[Proof of Theorem \ref{thm:factor-killed}]
In \cref{thm:factor-killed-bis}, we may take $\calC = \Ind(\freeduals^\rex_S \Tope{G}_\ast^\fin[(S^1)\inv]/S^V)$ and $F$ to be induced by the identity functor on $\Tope{G}_\ast$ to conclude that \[\Ind(\freeduals^\rex_S \Tope{G}_\ast^\fin[(S^1)\inv]/S^V) = 0\]
It follows that $\freeduals^\rex_S \Tope{G}_\ast^\fin[(S^1)\inv]/S^V=0$. Now let $F' : \freeduals^\rex_S \Tope{G}_\ast^\fin[(S^1)\inv]/S^V \to \freeduals^\SMP_S \Tope{G}_\ast[(S^1)\inv]/S^V$ be universally induced by inclusion $\freeduals^\rex_S \Tope{G}_\ast^\fin \subset \Tope{G}_\ast$. Because this functor is symmetric monoidal and its domain is $0$, it follows that the codomain $\freeduals^\SMP_S \Tope{G}_\ast[(S^1)\inv]/S^V = 0$ as desired.
\end{proof}



%
%



%
%


\end{document}